\documentclass[11pt]{amsart}
\usepackage[utf8]{inputenc}
\usepackage{rotating}
\usepackage{pifont}
\usepackage{fdsymbol}

\let\centerdot\undefine
\let\doteq\undefine

\usepackage{amssymb}
\usepackage{amsthm}
\usepackage{stmaryrd}
\usepackage{mathrsfs}
\usepackage[svgnames,x11names,dvipsnames]{xcolor}
\usepackage{hyperref}
\hypersetup{colorlinks=true,linkcolor=Teal,citecolor=Teal}
\usepackage[all]{xy}
\SelectTips{cm}{10}  % Better arrowheads
\usepackage{float}
\usepackage[top = 1in, bottom = 1in, left = 1.2in, right=1.2in, marginparwidth=1.0in]{geometry}
\usepackage{fancyhdr}
\usepackage{mathtools}
\usepackage{scalefnt}
\usepackage{microtype}
\usepackage{pifont}
\usepackage{ifthen}
\usepackage{marginnote}
\usepackage{graphicx}
\usepackage{soul}
\usepackage{enumitem}
\usepackage{mdframed}
\usepackage{tikz}
\usetikzlibrary{cd}
\usepackage{devanagari}

\usepackage[nolabel]{showlabels} % uncomment (and comment out the next line) to remove \label{} labels
\newcommand{\too}[1]{\overset{#1}\to}

\newcommand{\isom}{\cong}
\renewcommand{\dsum}{\oplus}
\newcommand{\dsums}{\bigoplus}
\newcommand{\tensor}{\otimes}

\renewcommand{\bar}[1]{{\overline{#1}}}

\newcommand{\hteq}{\simeq}
\newcommand{\sm}{\wedge}

\newcommand{\im}{\operatorname{im}}

\newcommand{\Fun}{\operatorname{Fun}}

\newcommand{\Map}{\operatorname{Map}}

\newcommand{\Mod}{\operatorname{Mod}}
\newcommand{\Tot}{\operatorname{Tot}}

\DeclareMathOperator*{\colim}{colim}

\newcommand{\Spec}{\operatorname{Spec}}
\newcommand{\Ext}{\operatorname{Ext}}

\newcommand{\Tor}{\operatorname{Tor}}

\newcommand{\R}{\mathbb{R}}

\newcommand{\Z}{\mathbb{Z}}

\newcommand{\br}[1]{\left[ {#1} \right]}

\renewcommand{\epsilon}{\varepsilon}

\renewcommand{\phi}{{\mathchoice{\raisebox{2pt}{\ensuremath\varphi}}{\raisebox{2pt}{\!\! \ensuremath\varphi}}{\raisebox{1pt}{\scriptsize$\varphi$}}{\varphi}}}

\newcommand{\attop}[1]{{\let\textstyle\scriptstyle\let\scriptstyle\scriptscriptstyle\substack{#1}}}
\renewcommand{\atop}[1]{{\let\scriptstyle\textstyle\let\scriptscriptstyle\scriptstyle\substack{#1}}}

\makeatletter
\newcommand{\switchmargin}{
\if@reversemargin
\normalmarginpar
\else
\reversemarginpar
\fi
}
\makeatother

\numberwithin{equation}{subsection} % Number equations by section, like (1.1) instead of (1)

\theoremstyle{definition}
\newtheorem{theorem}[equation]{Theorem}
\newtheorem{theorem-definition}[equation]{Theorem--Definition}
\newtheorem{corollary}[equation]{Corollary}

\newtheorem{lemma}[equation]{Lemma} 
\newtheorem{proposition}[equation]{Proposition}

\newtheorem{definition}[equation]{Definition}

\newtheorem{remark}[equation]{Remark}
\newtheorem{rmk}[equation]{Remark}
\newtheorem{notation}[equation]{Notation}
\newtheorem{construction}[equation]{Construction}
\newtheorem{assumption}[equation]{Assumption}

\newtheorem*{theoremstar}{Theorem}
\newtheorem*{defstar}{Definition}

\newcommand{\argforcustom}{}
\theoremstyle{newplain}
\newtheorem{helperforcustom}[equation]{\argforcustom}
\newtheorem*{helperforcustomstar}{\argforcustom}

\usepackage{longtable}

\newcommand{\Sp}{\mathrm{Sp}}
\newcommand{\SH}{\mathrm{SH}}
\newcommand{\Be}{\mathrm{Be}}
\newcommand{\ta}{\text{{\dn t}}}
\newcommand{\Ccal}{\mathcal{C}}
\newcommand{\Fil}{{\operatorname{Fil}}}
\newcommand{\Sdef}[1]{\mathbb{S}^{\mathrm{def}}_{#1}} 
\newcommand{\Dec}{\mathrm{D\acute{e}c}}

\newcommand{\cofib}{\operatorname{cofib}}
\newcommand{\fiber}{\operatorname{fiber}}
\newcommand{\ANSS}{\mathrm{ANSS}}
\newcommand{\AHSS}{\mathrm{AHSS}}
\newcommand{\hfpss}{\mathrm{hfpss}}
\newcommand{\fil}{{\mathrm{fil}}}
\newcommand{\mot}{{\mathrm{mot}}}
\newcommand{\II}{I\!I}
\newcommand{\abs}[1]{\lvert #1\rvert}

\usepackage[mathlines]{lineno}

\DeclareMathOperator{\evslice}{\textup{ev.slice}}
\DeclareMathOperator{\post}{\textup{post}}
\DeclareMathOperator{\evpost}{\textup{ev.post}}

\DeclareMathOperator{\acc}{\textup{acc}}
\DeclareMathOperator{\alg}{\textup{alg}}
\DeclareMathOperator{\id}{\operatorname{id}}

\DeclareMathOperator{\op}{\operatorname{op}}

\DeclareMathOperator{\map}{\operatorname{map}}
\DeclareMathOperator{\CAlg}{\operatorname{CAlg}}

\usepackage{mathtools}
\usepackage{bbm}
\usepackage{bbding}
\usepackage{stmaryrd}
\usepackage{comment}
\newcommand\spoke{{\scaleobj{1.2}{\Yright}}}
\usepackage{scalerel}[2016/12/29]
\newcommand\sbsq{\scaleobj{0.65}{\blacksquare}}
\newcommand{\mc}[1]{\mathcal{#1}}
\newcommand{\mb}[1]{\mathbb{#1}}
\newcommand{\mr}[1]{\mathrm{#1}}

\newcommand{\mit}[1]{\mathit{#1}}
\newcommand{\bra}[1]{\langle {#1} \rangle}
\newcommand{\CC}{\mb{C}}
\newcommand{\RR}{\mb{R}}
\newcommand{\FF}{\mb{F}}
\newcommand{\ZZ}{\mb{Z}}
\newcommand{\E}[2]{\prescript{#1}{#2}{E}}
\newcommand{\tmf}{\mr{tmf}}
\newcommand{\Tmf}{\mr{Tmf}}
\newcommand{\TMF}{\mr{TMF}}
\newcommand{\sstar}{\medwhitestar}
\newcommand{\doteq}{\overset{\scriptstyle{\centerdot}}{=}}

\mathchardef\mh="2D

\usepackage{graphicx}
\graphicspath{{./EU-symbol.jpg}} 

\usepackage{lipsum}

\usepackage[capitalise]{cleveref}

\newcommand{\highlighteva}[1]{\ifmmode{\text{\sethlcolor{Lavender}\hl{$#1$}}}\else{\sethlcolor{Lavender}\hl{#1}}\fi}

\newcommand{\highlightgabe}[1]{\ifmmode{\text{\sethlcolor{SeaGreen}\hl{$#1$}}}\else{\sethlcolor{SeaGreen}\hl{#1}}\fi}

\definecolor{violet}{rgb}{0.56, 0.0, 1.0}

\definecolor{imadeupthiscolor}{rgb}{1, 0.4, 0.4}

\newcommand{\tslice}[1]{P_{\geq #1}}
\newcommand{\tpost}[1]{\tau_{\geq #1}}

\newcommand{\subslice}[2]{\Sp^{\textup{slice}}_{#1 ,\geq #2}}
\newcommand{\subpost}[2]{\Sp^{\post}_{#1,\geq #2}}

\newcommand{\bR}{{\mathbb{R}}}
\newcommand{\MANSS}{\mathrm{\textup{MANSS}}}

\newcommand{\Gr}{\mathrm{Gr}}

\title{A deformation of Borel equivariant homotopy}
\date{\today}
\author{Gabriel Angelini-Knoll}
\address{Department of Mathematics, Applied Mathematics, and Statistics, Case
Western Reserve University, Cleveland, OH, USA}
\email{gabriel.angelini-knoll@case.edu}

\author{Mark Behrens}
\address{Department of Mathematics, University of Notre Dame, Notre Dame, IN, USA}
\email{mbehren1@nd.edu}

\author{Eva Belmont}
\address{Department of Mathematics, Applied Mathematics, and Statistics, Case
Western Reserve University, Cleveland, OH, USA}
\email{eva.belmont@case.edu}

\author{Hana Jia Kong}
\address{School of Mathematical Sciences, Zhejiang University, Hangzhou, Zhejiang, China}
\email{hana.jia.kong@gmail.com}

\begin{document}
\maketitle

\begin{abstract}
We describe a deformation of the $\infty$-category of Borel $G$-spectra for a finite group $G$. This provides a new presentation of the $a$-complete real Artin--Tate motivic stable homotopy category when $G=C_2$ and gives a new interpretation of the $a$-completed $C_2$-effective slice spectral sequence. As a new computational tool, we present a modified Adams--Novikov spectral sequence which computes the $RO(G)$-graded Mackey functor valued homotopy of Borel $G$-spectra. 
\end{abstract}

\setcounter{tocdepth}{2}
\tableofcontents

\section{Introduction}
%\color{Teal}
\subsection{Overview}
Motivic homotopy theory, developed by Voevodsky and Morel (eg. \cite{MorelVoevodsky}), 
was famously used as a
tool to attack
the Milnor conjecture and the Bloch--Kato conjecture, as well as other deep problems in number theory and algebraic geometry.
It became useful in stable homotopy theory for a different reason: the ``Betti
realization'' functor from the $\CC$-motivic homotopy category to the stable homotopy category makes computations of stable homotopy groups of spheres more accessible \cite{IWX}. The $\RR$-motivic homotopy category is also equipped with a Betti realization functor to the $C_2$-equivariant stable homotopy category; Guillou and Isaksen \cite{GI-C2} use this to aid computing $C_2$-equivariant stable homotopy groups. It is natural to ask whether there is a related category that
supports an analogous Betti realization functor to $G$-equivariant spectra for
other finite groups $G$. However, we do not expect such a comparison category to come
from motivic homotopy theory; the $C_2$ case works because $\CC$ is an extension
of $\RR$ with Galois group $C_2$, and the Artin--Schreier theorem implies that
$e$ and $C_2$ are the only finite groups that can occur as absolute Galois
groups. 

Instead, we turn to another perspective by building the desired category as a deformation of the equivariant category. Work of Gheorghe--Wang--Xu \cite{GheorgheWangXu} (see also \cite{GIKR}, \cite{Pstragowski}) provided a new
characterization of the category of ($p$-completed, cellular) $\CC$-motivic
spectra as a \emph{deformation} (see Definition \ref{def:deformation}) of
($p$-completed) ordinary spectra. The Tate twist element $\tau$ (see \cite[Lem.~23]{HuKrizOrmsby}) provides a deformation parameter, and Betti realization can be described as inverting $\tau$. Gheorghe--Isaksen--Krause--Ricka \cite{GIKR} provided a model for this deformation using the category of modules over a certain algebra in filtered spectra coming from the Adams--Novikov filtration.

The filtered spectra model can be further generalized: given a well-behaved filtration functor (see Definition \ref{def:tower}) in some appropriate category $\mathcal{C}$, the deformation perspective provides an appropriate new category $\mathcal{D}$ equipped with a comparison functor $\mathcal{D}\to \mathcal{C}$, analogous to the Betti realization functor. The construction comes with a ``deformation parameter'' $\tau$ and the Betti realization functor is given by inverting $\tau$. 
Using this perspective, Burklund--Hahn--Senger \cite{BHS} showed that a variant of the cellular
$p$-complete $\R$-motivic category (the Artin--Tate category) arises from a filtration functor on $C_2$-equivariant spectra.
However, this is not amenable to generalization to an arbitrary
finite group $G$, because this filtration hinges on the existence of the $C_2$-spectrum $MU_\R$ known as Real complex cobordism. No analogue of $MU_\R$ for other finite groups $G$ has been found despite over a decade of searching (see \cite{HillHopkinsRavenelodd, HahnSengerWilson}). 

In this paper, we overcome this difficulty in the setting of Borel $G$-spectra for a finite group $G$. We give an alternate construction of the filtration that doesn't involve $MU_\R$ and can be generalized to other groups. Our paper gives two interpretations of this filtration:
\begin{itemize} 
\item as a generalized Adams filtration in Borel $G$-spectra with an appropriate base (Definition \ref{def:main}), or
\item as a modified Adams--Novikov filtration which mixes the (non-equivariant) Adams--Novikov filtration and the homotopy fixed points filtration (see Section \ref{sec: MANSS}).
\end{itemize}
We show these are equivalent (Proposition \ref{prop:MANSS=NE-Adams})
and show that the $G=C_2$ case of our construction recovers the $a$-completion
of the motivic category studied in \cite{BHS}
(Theorem \ref{thm: MANSS Galois reconstruction}). The second construction produces a spectral sequence,
which we call the \emph{modified Adams--Novikov spectral sequence} (MANSS). We demonstrate that it is amenable to computations in Section \ref{sec:applications}.

In the case of $G=C_2$, our construction provides a new description of a completed form of the \emph{$C_2$-effective slice spectral
sequence} introduced by the fourth author in \cite{Kon20}. The $C_2$-effective slice spectral
sequence is defined using
Betti realization and Voevodsky's effective slice filtration on the $\R$-motivic
homotopy category \cite{Voev02}. 
It agrees with the Hill--Hopkins--Ravenel slice spectral sequence \cite{HHR16} for Real orientable $C_2$-spectra, and in the non-Real orientable case it is more tractable than the slice spectral sequence. 
Burklund--Hahn--Senger \cite{BHS} have shown that the $C_2$-effective slice spectral
sequence agrees with the d\'ecalage (a degree-shifting operation) of the $MU_\R$-based Adams--Novikov spectral
sequence, and hence (after completion) with the d\'ecalage of the MANSS. Moreover, for other
groups $G$, the MANSS can be viewed as a Borel $G$-equivariant analogue of the
$C_2$-effective spectral sequence up to d\'ecalage.

\subsection{The main constructions}
We will use the presentations expressed in the language of filtered spectra by \cite{GIKR} and \cite{BHS}. A \emph{filtered spectrum} $X_\bullet\in \mathcal{C}^\Z$ is simply a $\ZZ$-indexed sequence of objects $X_i$ and maps
\[ 
	\cdots \leftarrow X_{-2} \leftarrow X_{-1} \leftarrow X_0 \leftarrow X_{1} \leftarrow X_2 \leftarrow \cdots 
\]
in a category $\mathcal{C}$.
We regard this sequence as a decreasing filtration of $X_{-\infty} \coloneqq \colim_i X_{-i}$. There is a symmetric monoidal structure on filtered objects given by \emph{Day convolution} \eqref{eq:day-convolution}.

Given a filtered object $X_\bullet$ and a Postnikov truncation functor
$\tau_{\geq \bullet}$, one can apply a degree-shifting operation called \emph{d\'ecalage} (see
\S\ref{sec:decalage}) to form another filtered object $\Dec^{\tau}(X)_\bullet$. The central idea, as introduced in \cite{GIKR}, is to study the category of modules over $\Dec(X)$:
\begin{equation}\tag{*} \Mod(\mathcal{C}^\Z; \Dec^{\tau}(X)). \end{equation}

The filtrations relevant to this paper are
\begin{enumerate} 
\item the \emph{Adams--Novikov filtration} $\bar{MU_p}^{\tensor \bullet}\tensor X$ of a
spectrum $X$ (where $\bar{(-)}$ denotes the fiber of the unit map);
\item the \emph{Real Adams--Novikov filtration} $\bar{MU_{\R,2}}^{\tensor
\bullet}\tensor X$ of a $C_2$-spectrum $X$;
\item the \emph{$MU[G]$-based Adams filtration} $\bar{MU_p \sm G_+}^{\tensor \bullet}\tensor X$ of a (Borel) $G$-spectrum $X$;
\item the \emph{modified Adams--Novikov filtration}, defined to be the Day convolution in Borel $G$-spectra of
	\begin{enumerate} 
	\item the \emph{homotopy fixed points filtration} of a Borel $G$-spectrum
	$X^h = F(EG_+,X)$, given by filtering $EG$ by skeleta, and
	\item the (nonequivariant) Adams--Novikov filtration of the sphere.
	\end{enumerate}
\end{enumerate}

Filtrations (1) and (2) pertain to the motivating results; (3) and (4) are our constructions.

\begin{theorem}[\cite{GIKR}]
Applying the construction (*) to the Adams--Novikov filtration with $\tau$ as the Postnikov filtration recovers the cellular $\CC$-motivic category $\SH^{\text{cell}}(\CC)_{i2}$.
\end{theorem}
Here and henceforth, we use $(-)_{ip}$ to denote the category of modules over the $p$-complete unit.

\begin{theorem}[\cite{BHS}]
Applying the construction (*) to the Real Adams--Novikov filtration with $\tau$ as even slice d\'ecalage (Notation \ref{notation:decalage}) gives the Artin--Tate category $\SH(\R)^{AT}_{i2}$, generated by bigraded suspensions $\Sigma^{i,j}\mr{Spec}(k)_+$, where $k \in \{\RR,\CC\}$.
\end{theorem}

Note that $\Spec(\RR)_+$ is the sphere spectrum, so this is analogous to the
cellular category but with one more generator.

Our main construction is the following.

\begin{defstar}[Definition~\ref{def:main}]\label{maindef1}
Define the category of \emph{complete artificial motivic spectra}
$\widehat{\SH}_{G,p}^{\text{art}}$ to be the construction (*) applied to the
$MU[G]$-based Adams filtration, with $\tau$ as the even Postnikov filtration.
\end{defstar}

Our main results are the following three theorems, which justify the definition above by relating it to the construction by \cite{BHS} in the case $G=C_2$. The first theorem relates (2) and (3) in the case $G=C_2$. The third theorem says that that (3) has the same spectral sequence as (4).

\begin{theoremstar}[Theorem \ref{ESSS is MANSS}]
For $X \in \Sp_{C_2, i2}$, there is an equivalence
\[ 
	\Dec^{\textup{ev.post}} (\overline{MU_2\sm {C_2}_+}^{\tensor \bullet}\tensor X) \simeq \Dec^{\textup{ev.slice}}(\overline{MU_{\R,2}}^{\tensor \bullet}\tensor X)^{h}
\]
of filtered Borel $C_2$-spectra, where $(-)^{h} = F((EC_2)_+, -)$ denotes Borel completion, $\Dec^{\textup{ev.post}}$ denotes the even Postnikov d\'ecalage, and $\Dec^{\textup{ev.slice}}$ denotes the even slice d\'ecalage.
\end{theoremstar}

\begin{theoremstar}[Theorem~\ref{thm: MANSS Galois reconstruction}]\label{mainthm-reconstruction}
There is an equivalence
\[ 
\widehat{\SH}_{C_2,2}^{\text{art}}\hteq	(\SH(\R)^{\textup{AT}}_{i2})_{a}^{\wedge}
\]
where the completion of $\SH(\RR)^{\mathrm{AT}}_{i2}$ is with respect to the motivic lift $a$ of the equivariant Euler class $a_\sigma$ defined in \cite[Eq.~(1)]{BHS}.
\end{theoremstar}

\begin{theoremstar}[Theorems~\ref{thm:MANSS=MU/a-Adams}]\label{mainthm2}
For $G=C_2$,
the spectral sequence associated to the
modified Adams--Novikov filtration agrees with the $\mathrm{MU}_{p}[C_2]$-based Adams spectral sequence from the $E_2$-term onwards. In particular, the construction (*) applied to the modified Adams--Novikov filtration also recovers $(\SH(\R)^{\textup{AT}}_{i2})_{a}^{\wedge}$.
\end{theoremstar}

While the comparison theorems use filtration (3), we are interested in the modified Adams--Novikov filtration because it is amenable to computation.

\subsection{Computing the MANSS}
We study the \emph{modified Adams--Novikov spectral sequence} (MANSS), the spectral sequence associated to the modified Adams--Novikov filtration as discussed in the previous section. Computing  this spectral sequence is relevant because, in many cases such as ours, computing the homotopy groups of a construction of the form (*) boils down to computing the spectral sequence associated to the filtration in question; see \cite[Theorem A.1]{BHS-manifolds}.

The MANSS is an $RO(G)$-graded
spectral sequence which is a kind of mixture of the (non-equivariant) Adams--Novikov spectral
sequence (ANSS) and the homotopy fixed point spectral sequence (HFPSS), and as
 such is very amenable to computations.

We introduce the following tools for computing with the MANSS:

\begin{enumerate}
\item In Section \ref{sec: E2 MANSS}, we introduce a pair of double complex spectral sequences which converge to the $E_2$-term of the $p$-primary MANSS:
\begin{align*}
{}_{I}\underline{E}^{*,*,\star}_2(G/H) = &\Ext^{*,*}_{(MU_{p})_*MU_{p}}(MU_*, H^*(H; (MU_{p})_*\Sigma^{-\star} X))  \implies {}^{\MANSS}\underline{E}_2^{*,*}(X)(G/H), \\
{}_{\II}\underline{E}^{*,*,\star}_2 (G/H)= &H^*(H; \Ext^{*,*}_{(MU_{p})_*MU_{p}}((MU_{p})_*, (MU_{p})_*\Sigma^{-\star}X))  \implies {}^{\MANSS}\underline{E}_2^{*,*}(X)(G/H).
\end{align*}

\item \label{it: 2 of tools for compting MANSS} We show these double complex spectral sequences collapse in the case where $X$ is the Borel completion of a genuine $C_2$-spectrum which is $MU_\RR$-\emph{projective} (i.e. its $MU_\RR$-homology is free on a set of generators in degrees which are multiples of $\rho$).  We give a complete computation of the $E_2$-term of the $2$-primary MANSS for such spectra in terms of the ANSS $E_2$-term of the underlying non-equivariant spectrum $X^e$, in a manner that mirrors the computation of the $E_2$-term of the $C_2$-effective slice spectral sequence for $MU_{\RR}$-projective spectra \cite[Thm.~5.6]{BHS}.
\vspace{10pt}

\item  \label{it: 3 of tools for compting MANSS} In Section \ref{sec: triv action case} we show that these double complex spectral sequences also collapse in the case where $X$ is a Borel $C_2$-spectrum with trivial $C_2$-action and torsion-free $MU$-homology.  In this case, the MANSS converges to 
\[ 
	\pi^{C_2}_{i+j\sigma}F((EC_2)_+,X) \cong \pi_i(F(P_j^\infty, X)), 
\]
and we give a complete computation of the MANSS $E_2$-term of such spectra in terms of the ANSS $E_2$-term of $X$, which is different from \eqref{it: 2 of tools for compting MANSS}, where the generators are now parameterized by the cells of the associated stunted projective spectra.
\vspace{10pt}
\end{enumerate}

For $X$ as in \eqref{it: 3 of tools for compting MANSS}  with $MU$-homology concentrated in even degrees, we introduce a complex motivic lift of the MANSS (the \emph{motivic MANSS}) 
\[ 
	\E{\MANSS}{\mot}^{*,\star,*}_r(\nu_\CC X) 
\]
in which $X$ is replaced with the motivic lift $\nu_{\CC}(X) \in \SH(\CC)_{i2}$ and the stunted projective spaces $P_j^\infty$ are replaced with the motivic lifts $\nu_{\CC}(P_j^\infty)$.  We have
\[ 
	\E{\MANSS}{\mot}_r^{*,\star}(X) \cong \E{\MANSS}{\mot}^{*,\star,*}_r(\nu_\CC X)[\tau^{-1}]. 
\]
We prove the following theorem.

\begin{theoremstar}[Theorem~\ref{thm:AHSS}]
The motivic MANSS 
\[
	\E{\MANSS}{\mot}^{*,*+j\sigma,*}_r(\nu_\CC(X) \otimes C\tau)
\]
is isomorphic (after reindexing) to the accelerated algebraic Atiyah--Hirzebruch spectral sequence (AHSS)
\[ 
	\E{\AHSS}{\acc,\alg}_1^{k,s,t}(X^{P^{\infty}_j}) = \Ext^{s,t}_{MU_*MU}(MU_*, MU_*(X^{P^{2k}_j/P^{2k-2}_j})) \implies \lim_m \Ext^{s,t}_{MU_*MU}(MU_*, MU_*(X^{P_j^{2m}})) 
\]
obtained from filtering $P_j^\infty$ by its even skeleta.\footnote{Because the filtration is by even skeleta instead of skeleta, we refer to this spectral sequence as the \emph{accelerated algebraic Atiyah--Hirzebruch spectral sequence}.}
\end{theoremstar}

The significance of this theorem is that the differentials in the algebraic Atiyah--Hirzebruch spectral sequence are computable in terms of the attaching maps of $P_j^\infty$.  We can then employ the zig-zag of spectral sequences

\begin{equation}\label{eq:AHSSMANSS}
\xymatrix@C-1em@R-1em{
\E{\MANSS}{mot}^{*,*+j\sigma,*}_{2r+1}(\nu_\CC(X) \otimes C\tau) \ar@{=}[d] &  
\E{\MANSS}{\mot}^{*,*+j\sigma,*}_{2r+1}(\nu_\CC(X)) \ar[l] \ar[r]  & 
\E{\MANSS}{\mot}^{*,*+j\sigma,*}_{2r+1}(\nu_\CC(X))[\tau^{-1}] \ar@{=}[d]
\\
\E{\AHSS}{\acc,\alg}_r^{*,*,*}(X^{P^\infty}_j)
 && \E{\MANSS}{}^{*,*+j\sigma}_{2r+1}(X)
}
\end{equation}
to lift differentials from the algebraic AHSS to the MANSS.

For any $G$ and $X$ we observe that the restriction map 
\[ 
	\pi^G_{V}(X^h) \to \pi_{\abs{V}}(X) 
\]
induces a map of spectral sequences

\begin{equation}\label{eq:ANSSMANSS} 
\E{\MANSS}{}^{*,\star}_r(X) \to \E{\ANSS}{}^{*,\abs{\star}}_r(X) 
\end{equation}
from the MANSS to the Adams--Novikov spectral sequence of the underlying non-equivariant spectrum $X$ (Corollary~\ref{cor: MANSS to ANSS}).  This gives another means of determining MANSS differentials.

\subsubsection{Computational examples}
We end this paper with some explicit computations of the MANSS.

\begin{itemize}
\item In the case of $G = C_2$, $p = 2$, and $X = ko_{C_2}$, we show that all of the differentials and multiplicative structure in the MANSS for $ko_{C_2}$ can be deduced from (\ref{eq:AHSSMANSS}) and (\ref{eq:ANSSMANSS}).  This recovers the main generating differentials in the $C_2$-effective slice spectral sequence determined by the fourth author in \cite{Kon20}.

\item In the case of $G = C_3$, $p = 3$, and $X = \tmf(2)$, the MANSS agrees with the $RO(C_3)$-graded homotopy fixed point spectral sequence for $\tmf(2)$.  Following \cite{HillHopkinsRavenelodd,HahnSengerWilson}, we extend the $RO(C_3)$-grading to an $RO^\spoke(C_3)$-grading which includes the ``spoke sphere'' $S^\spoke$.  Wilson computed the slice spectral sequence for $\tmf(2)$ in \cite{Wilsonthesis}.  We record an observation of Senger,  that the $\mu_3$-orientation of $\tmf(2)$ \cite{HahnSengerWilson}, the Hill--Hopkins--Ravenel norm, and known differentials in the integer graded homotopy fixed point spectral sequence (which is closely related to the decent spectral sequence for $\tmf_{(3)}$) completely determine the structure of the $RO^\spoke(C_3)$-graded MANSS. 
\end{itemize}

\subsection{Future work}
In future work, we plan to investigate how to extend our deformation $\widehat{\SH}^{\mr{art}}_{G,p}$ of $\Sp^{BG}_{ip}$ to a deformation $\SH^\mr{art}_{G,p}$ of $\Sp_{G,ip}$. {Let $\Sdef{C_p,p}$ denote $\Dec^{\textup{ev.post}}(MU_p[G]^{\tensor \bullet+1})$}. 
In the case of $G = C_p$, the idea would be to define the desired filtered $C_p$-spectrum $\Gamma_{C_p}\mb{S}_p$ as a pullback

\[
\xymatrix{
\Gamma_{C_p}\mb{S}_p \ar[r] \ar[d] &
\Gamma^{\Phi}_{C_p} a_\lambda^{-1}\mb{S}_{p} \ar[d] \\
\Sdef{C_p,p} \ar[r] &
a_\lambda^{-1}\Sdef{C_p,p}
}
\]
where $a_\lambda$ is the Euler class of the $2$-dimensional representation of $C_p$ given by rotation by $2\pi/p$.  This would require the construction of a suitable filtration $\Gamma^{\Phi}_{C_p} a_\lambda^{-1}\mb{S}_p$ of the geometrically local sphere $a^{-1}_\lambda \mb{S}_p$ which maps to the localized filtration $a_\lambda^{-1}\Sdef{C_p,p}$.

In the case of $G = C_2$, this amounts to understanding the $a$-localization of the $\RR$-motivic sphere in the Artin--Tate category.  Unfortunately, this is the one object that the authors of \cite{BHS} flag as very mysterious, so the proposed program will encounter significant challenges.  Nevertheless, the completion of this program would yield a notable application.  The ``$G$-effective slice spectral sequence'' 
\[ 
	E_r^{*,\star}(\Gamma_G X) 
\]
associated to the filtered object $\Gamma_G X$ would potentially allow for much of the computational payoffs of the hypothetical odd primary analog $BP_{\mu_p}$ of $BP_\RR$ discussed in \cite{HillHopkinsRavenelodd}, while circumventing the need to construct $BP_{\mu_p}$ itself. Note that Hill--Hopkins--Ravenel \cite[p.1-5]{HillHopkinsRavenelodd} have suggested that if such a spectrum $BP_{\mu_p}$ were known to exist, the its norm 
\[N_{C_p}^{C_{p^2}}BP_{\mu_p}\] 
could be used to resolve the $3$-primary Kervaire invariant problem.  

\subsection{Organization of the paper}
In Section~\ref{sec:foundations}, we introduce the foundations of filtered spectra which will be used throughout the paper.  We end this section with our definition of complete artificial motivic spectra.

In Section~\ref{sec:two-sseq}, we recall the construction of the $C_2$-effective slice spectral sequence.  We then define, for $E$ an arbitrary ring spectrum, the $E$-modified Adams spectral sequence.  We prove that this modified Adams spectral sequence is actually the Adams spectral sequence based on the Borel $G$-spectrum $E[G]$.  The modified Adams--Novikov spectral sequence (MANSS) is the instance of this spectral sequence where $E = MU_p$.  We finish this section by proving that the Borel completion of the $C_2$-effective slice filtration is the even Postnikov d\'ecalage of the MANSS filtration.

In Section~\ref{sec: computations}, we introduce our main computational tools.  We define the double complex spectral sequences, and compute their $E_2$-terms.  We then consider for $G = C_2$ the special case when $X$ is the Borel completion of an $MU_\RR$-projective $C_2$-spectrum (see Definition \ref{defn:MURprojective}), and the special case when $X$ has a trivial $C_2$-action.  We define the motivic MANSS, and show that mod $\tau$, the motivic MANSS is isomorphic to the algebraic Atiyah--Hirzebruch spectral sequence.

In Section~\ref{sec:applications}, we apply our methods to completely compute the MANSS in the case of the $C_2$-spectrum $ko_{C_2}$, and the $C_3$-spectrum $\tmf(2)$.

\subsection{Conventions and notation}\label{conventions}
Throughout, we will use the theory of $\infty$-categories (quasi-categories) as developed by \cite{Joy02,HTT,HA}. Since the nerve functor from $1$-categories to $\infty$-categories is fully faithful, we will not include this functor in the notation. We will also refer to $\infty$-categories simply as categories from this point onwards. Throughout, $G$ will denote a finite group. We use the following notations:
\begin{longtable}{lp{10pt}l}
   $\mathbb{Z}$  && $(\mathbb{Z},\le )$ as a partially ordered set\\
   $\mathbb{Z}^{\delta}$ && $\mathbb{Z}$ as a discrete category\\
   $\mathcal{C}^{\Fil}$ && $\Fun(\mathbb{Z}^{\op},\mathcal{C})$\\
    $\mathcal{C}^{\Gr}$ && $\Fun(\mathbb{Z}^{\delta},\mathcal{C})$\\
	  $\Pr^{L}$ && the category of presentable categories and left adjoints\\
      $\Pr^{L,\text{st}}$ && the category of presentable stable categories and left adjoints\\
      $\Sp$, $\Sp_G$, $\Sp^{BG}$ && the category of spectra, genuine $G$-spectra, Borel $G$-spectra\\
   $\SH(k)$ && Morel--Voevodsky $\infty$-category of motivic spectra over $k$ \cite{MorelVoevodsky}\\
   $\mathcal{C}_{ip}$ && $\infty$-category of modules over the $p$-complete unit in $\mathcal{C}$\\
   $\tensor$ && symmetric monoidal product in the category at hand\\
   $\Ext^{s,t}((MU_{p})_*X)$ && $\Ext^{s,t}_{(MU_{p})_*MU_{p}}((MU_{p})_*, (MU_{p})_*X)$
\end{longtable}

Recall that a presentably symmetric monoidal, stable category is an object in the category $\CAlg(\Pr^{L,\textup{st}})$. We follow \cite[p.16]{BHS} and use the following conventions: Given $\mathcal{C}\in \CAlg(\Pr^{L,\textup{st}})$, let $\mathcal{C}_p$ denote the category of $p$-complete objects in
$\mathcal{C}$.
For $\mathcal{D}\in \{\mathcal{C},\mathcal{C}^{\Fil},\mathcal{C}^{\Gr}\}$ and $\mathcal{C} \in \{\Sp,\Sp_G,\Sp^{BG}\}$, let $X_p$ denote the $p$-completion of
$X\in \mathcal{D}$ and let $\mathcal{D}_{ip}=\Mod(\mathcal{D};\mathbbm{1}_p)$ denote the category of modules over $\mathbbm{1}_p$.  When working in one of these categories $\mathcal{D}_{ip}$, we will write $\otimes$ for the symmetric monoidal structure $\otimes_{\mathbbm{1}_p}$. We use this convention in Definition \ref{def:main}, Proposition \ref{prop: deformation} and throughout Section \ref{sec: esss}, Section \ref{sec: ESSS vs MANSS}, Section \ref{sec: computations}, and Section \ref{sec:applications}. Note that we still use the symbol $\otimes_{R}$ to mean the (underived) tensor product in the category of modules over a discrete ring and we use $\otimes$ to mean the usual tensor product of (graded) abelian groups in these sections. We expect that the meaning is clear from the context. We also write $X_E$ for Bousfield localization of a spectrum $X$ at a spectrum $E$. 

We write $\map_{\mathcal{C}}$ for the mapping spectrum of a stable category $\mathcal{C}$ or simply $\map$ when the context is clear. Similarly, we write $\Map_{\mathcal{C}}$ for the mapping space of a category $\mathcal{C}$ and $\Map$ when the category $\mathcal{C}$ is clear from context. If $\mathcal{C}$ is $\mathcal{D}$-linear for some $\mathcal{D}\in \CAlg(\Pr^{L,\text{st}})$ we write $\map_{\mathcal{C}}^{\mathcal{D}}(-,-)\colon \thinspace \mathcal{C}^{\op}\times \mathcal{C}\to \mathcal{D}$ for the mapping object in $\mathcal{D}$ of the category $\mathcal{C}$ or simply $\map^{\mathcal{D}}$ when $\mathcal{C}$ is clear from the context. 

For $X \in \Sp_G$, and $V \in RO(G)$, we write
\[ 
	\pi^H_V(X) \coloneqq [\mathbb{S}^V,X]^H = \pi_0((\Sigma^{-V}X)^H)
\]
for the $RO(G)$-graded $H$-equivariant homotopy groups of $X$.  As $H$-varies, these form a Mackey functor, which we denote
\[ 
	\underline{\pi}_V(X). 
\]
For $X \in \Sp^{BG}$, we will similarly define
\[ 
	\pi^H_V(X^h) \coloneqq \pi_0((\Sigma^{-V}X)^{hH}), 
\]
and we denote the associated Mackey functor
\[ 
	\underline{\pi}_V(X^h). 
\]
We shall write $\{ E^{f,n}_r\}$ to denote a spectral sequence valued in abelian groups, and 
$\{\underline{E}^{f,n}_r\}$ to denote a spectral sequence valued in $G$-Mackey functors.  Given a spectral sequence $\{\underline{E}^{f,n}_r\}$ valued in Mackey functors, we will write $\{ E^{f,n}_r \}$ (removing the underline) to denote the value of the spectral sequence of Mackey functors at $G/G$:
\[ 
	E_r^{f,n} \coloneqq \underline{E}^{f,n}_r(G/G). 
\]
 
\subsection{Acknowledgements}
This project was begun in May 2021 at an American Institute of Mathematics workshop \emph{Equivariant techniques in stable homotopy theory} and finalized during a SQuaRE at the American Institute of Mathematics. The authors thank AIM for providing a supportive and mathematically rich environment. This project has received funding from the European Union's Horizon 2020 research and innovation programme under the Marie Sk\l{}odowska-Curie grant agreement No 1010342555. \thinspace \includegraphics[scale=0.1]{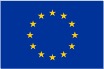} 
The second author completed part of this project as a visiting scholar at Northwestern University.
%The second author was partially supported by the NSF grant DMS-2005476, and he completed part of this project as a visiting scholar at Northwestern University. The third author was supported by the NSF grant DMS-2204357. The fourth author was supported by the NSF grant DMS-1926686.
This material is based upon work supported by the U.S. National Science Foundation under award No. DMS-2005476, DMS-2204357, DMS-1926686, and DMS-2506564. \thinspace \includegraphics[height=1em]{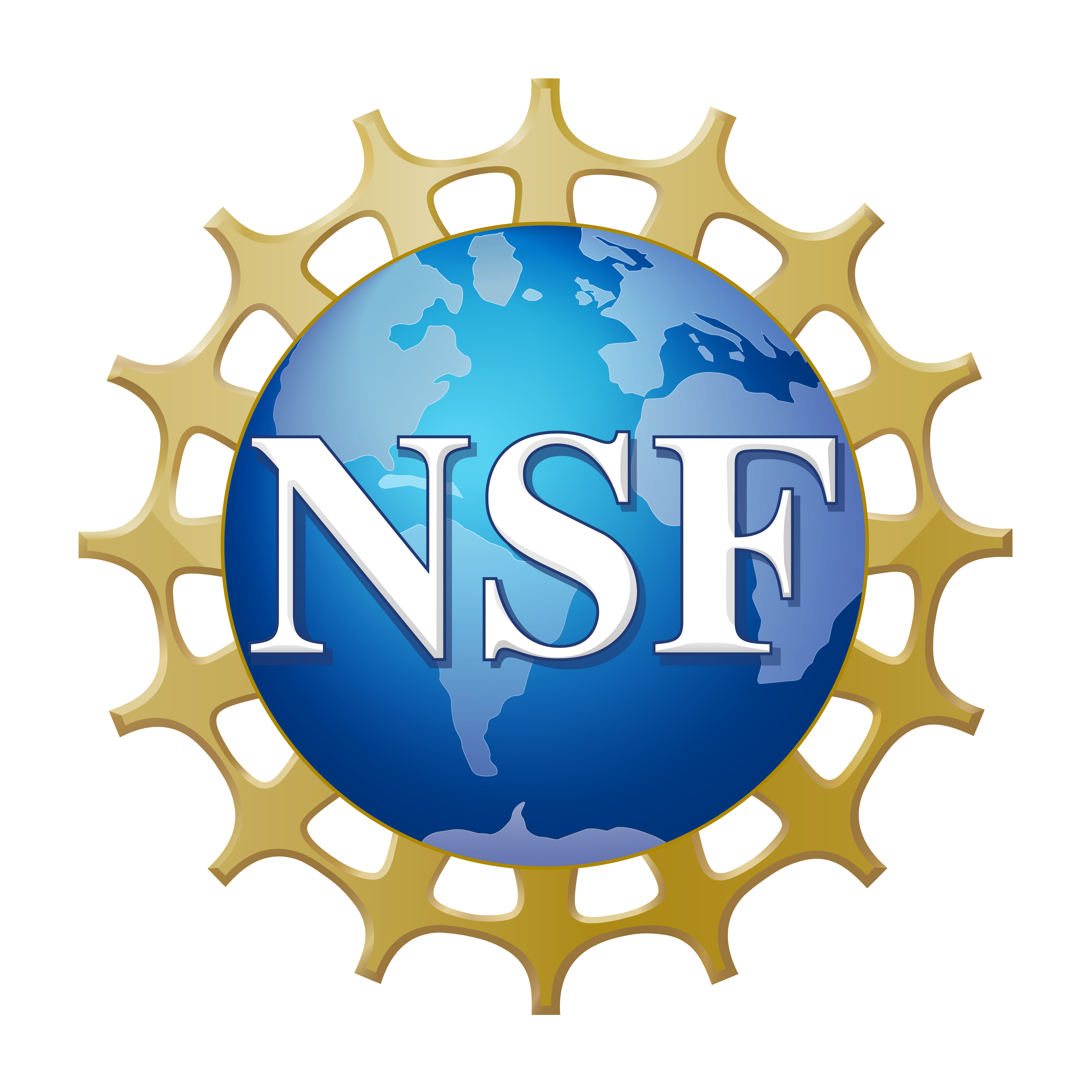} The fourth author was supported by ``the Fundamental Research Funds for the Central Universities'' No. 226-2025-00092 (China).
The authors thank Andrew Senger, for generously sharing his computations of the $RO^{\spoke}(C_p)$-graded homotopy fixed-point spectral sequence of $E_{k(p-1)}$.  The authors also benefited from conversations with Benjamin Antieau, Robert Burklund, Jeremy Hahn, Mike Hill, Max Johnson, and Noah Riggenbach. 

\section{Filtered objects and deformations}\label{sec:foundations}
Throughout this section, let 
$\mathcal{C}$ be a presentably symmetric monoidal, stable category and assume $\mathcal{C}$ is generated by a set $S$ of compact dualizable objects. In this situation, we say that the set $S$ \emph{rigidly generates} $\mathcal{C}$. 
\subsection{Filtered objects} \label{sec:filtered-objects}
Here we summarize results about filtered and graded objects in $\mathcal{C}$, denoted $\mathcal{C}^{\Fil}$ and $\mathcal{C}^{\Gr}$ (see Section \ref{conventions}). By our running assumptions on $\mathcal{C}$, we know that  $\mathcal{C}^{\Fil},\mathcal{C}^{\Gr}\in \CAlg(\Pr^{L,\text{st}})$  where we equip each of these categories with the Day convolution symmetric monoidal structure \cite{Day71,Gla16}. There is an associated graded functor 
\[ 
	\Gr^*\colon \thinspace \mathcal{C}^{\Fil}\to \mathcal{C}^{\Gr},
\] 
which is symmetric monoidal \cite[Theorem 1.13]{GP18}. 
Given a filtered spectrum $A_{\bullet}$ we define $A(-i)_{\bullet}$ such that $A(-i)_j=A_{j+i}$. 
Letting $\mathbbm{1} _{\mathcal{C}}$ denote the unit of the symmetric monoidal product on $\mathcal{C}$, we write $\Sigma^s\mathbbm{1} _{\mathcal{C}}(-t)$
for the $s$-th suspension of the symmetric monoidal unit shifted by $t$ and note that the homotopy groups in the category of graded objects in $\mathcal{C}$ are bigraded in the sense that for $A\in \mathcal{C}^{\Gr}$,
\[ \pi_{s,t}A\coloneqq \pi_{0}\text{map}(\Sigma^s \mathbbm{1}_{\mathcal{C}} (-t),A).\]
If $\mathcal{C}$ is $\Sp_{G}$-linear, we have a functor 
\[
	\text{map}^{\Sp_{G}}(-,-)=\mathcal{C}^{\text{op}}\times \mathcal{C}\to \Sp_{G},
\]
and $\mathcal{C}^{\Fil}$ and $\mathcal{C}^{\Gr}$ are also $\Sp_{G}$-linear.  
Note that since $\mathcal{C}$ is tensored over $\Sp_G$ we can can tensor with representation spheres to produce $\Sigma^VX\coloneqq \mathbb{S}^V\tensor X$ for any $X\in \mathcal{C}$. We can therefore consider the $RO(G)\times \mathbb{Z}$-graded Mackey functor valued homotopy of $A\in \mathcal{C}^{\Fil}$ with value
\[ \underline{\pi}_{V,t}A\coloneqq \underline{\pi}_{0}\text{map}^{\Sp_{G}}(\Sigma^{V}\mathbbm{1} _{\mathcal{C}}(-t),A)\]
at $(V,t)\in RO(G)\times \mathbb{Z}$. 

\begin{construction} \label{construction:fil-sseq}
Given an object $I\in\mathcal{C}^{\Fil}$, the spectral sequence associated to the filtered object $I$ is the
spectral sequence that comes from applying $\pi_*(-)\coloneqq \pi_*\map(\mathbbm{1} _{\mathcal{C}},-)$ to the tower of cofiber sequences
\[ 
    \xymatrix{
        \dots  & I_n\ar[l]\ar[d] & I_{n+1}\ar[l]\ar[d] & \dots\ar[l]
        \\ & I_n/I_{n+1} & I_{n+1}/I_{n+2}
    }
\]
where $I_n/I_{n+1}\coloneqq \cofib(I_{n+1}\to I_{n})$.
If $I_{\infty}=0$ then we say $I$ is complete and the associated conditionally convergent spectral sequence has signature 
\[E_1^{f,n}(I_\bullet) =\pi_{n}(\Gr^f I)\implies \pi_{n}I_{-\infty} \]
where $I_{-\infty}$ is the colimit.  The differentials take the form
\[ 
	d_r\colon \thinspace E_r^{f,n} \to E_r^{f+r,n-1}. 
\]

We will also use a variant of this construction for $\Sp_{G}$-linear categories.
Suppose $\mathcal{C}$ is $\Sp_{G}$-linear and $I\in \mathcal{C}^{\Fil}$. 
When $I_{\infty}=0$ then we say $I$ is complete and there is a conditionally convergent spectral sequence of Mackey functors with signature 
\[ 
	\underline{E}_1^{f,V}=\underline{\pi}_{V}(\Gr^f I)\implies \underline{\pi}_{V}I_{-\infty}.
\]
The differentials take the form
\[ 
	d_r\colon \thinspace \underline{E}_r^{f,V} \to \underline{E}_r^{f+r,V-1}. 
\]
\end{construction}
One source of filtered objects comes from cosimplicial objects.   Given a cosimplicial object
\[ 
	X^\bullet: \Delta \to \mc{C} 
\]
we get a filtered object $\{\fil_n X^{\bullet}\}$
by taking 
\begin{equation}\label{eq:cosimpfilt}
 \fil_n X^\bullet \coloneqq \mr{fiber}(\Tot X^\bullet \to \Tot_{n-1} X^\bullet)
\end{equation}
(where by convention, we define $\Tot_k X^\bullet$ to be trivial if 
$k$ is negative).  The associated spectral sequence is then the Bousfield--Kan spectral sequence, which we will denote:
\[ 
	E_r^{s,n}(X^\bullet) \coloneqq E_r^{s,n}(\fil_\bullet X^\bullet). 
\]

\begin{remark}
Our indexing convention differs from the standard indexing convention for the Bousfield--Kan spectral sequence.  To illustrate this, suppose that $E$ is a flat homotopy commutative ring spectrum, and that $X$ is a spectrum.  Then (with our indexing convention) the $E$-based Adams Spectral Sequence ($E$-ASS) for $X$ is the spectral sequence associated to the cosimplicial object $E^{\otimes \bullet + 1} \otimes X$ with
\[ 
	{}^{E\mh \textup{ASS}}E^{s,n}_2 \coloneqq E^{s,n}_2(E^{\otimes \bullet + 1} \otimes X) \cong \Ext^{s,s+n}_{E_*E}(E_*, E_*X) \implies \pi_n(X^{\wedge}_E)  
\]
where $X^{\wedge}_E \coloneqq \Tot E^{\otimes \bullet + 1} \otimes X$ is the $E$-nilpotent completion of $X$.\\
\end{remark}

Suppose $I$ and $J$ are filtered objects. Let $I\tensor J$ denote the Day convolution \cite{Gla16}; concretely, this is 
\begin{equation}\label{eq:day-convolution} (I\tensor J)_n = \colim_{i+j\geq n} I_i\tensor J_j \end{equation}
on objects.

\begin{lemma} \label{lemma:E1}
Let $d_1^I$ and $d_1^J$, respectively, denote the $d_1$ differentials in Construction \ref{construction:fil-sseq} associated to the filtered objects $I,J$.
Then we have 
\begin{equation}\label{eq:E1}
\Gr^n(I\tensor J) = \bigoplus_{i+j=n} \Gr^i I \tensor \Gr^j
J
\end{equation}
and whenever
\[ 
    x\otimes y\in \im\Big(\pi_s(\Gr^i I)\otimes \pi_s(\Gr^j
J)\longrightarrow \bigoplus_{i+j=n}\pi_s(\Gr^i I \tensor \Gr^jJ )
\Big)
\]
we have 
\[
	d_1^{I\otimes J}(x\otimes y)=d_1^{I}(x)\otimes y+(-1)^{|x|}x\otimes d_1^J(y).
\]
\end{lemma}

\begin{proof}
Equation \eqref{eq:E1} follows from the fact that the functor 
\[\Gr^{*} \colon \thinspace \mathcal{C}^\Fil\to \mathcal{C}^\Gr\]
is symmetric monoidal \cite[Theorem 1.13]{GP18} (see also \cite[Theorem
3.3.10]{AKS18}). For the statement about $d_1^{I\tensor J}$, we first claim that if $z\in \pi_s(\Gr^i I \tensor \Gr^jJ)$ then 
\begin{equation}
\label{eq:d1z-Gr} d_1(z)\in \pi_s(\Gr^{i+1} I \tensor \Gr^jJ) \oplus \pi_s(\Gr^i I \tensor \Gr^{j+1}J).
\end{equation}
To see this, let $P_{i,j}$ denote the pushout in the diagram
\[ 
	\xymatrix{
		I_{i+1}\tensor J_{j+1}\ar[r]\ar[d] & I_i\tensor J_{j+1}\ar[d]
		\\I_{i+1}\tensor J_j\ar[r] & P_{i,j}
		}.
\]
There is a map of cofiber sequences
\begin{align}\label{eq: d1 diagram}
 	\xymatrix{
		\Sigma^{-1}\Gr^i I \tensor \Gr^j J\ar[r]\ar[d] & P_{i,j}\ar[r]\ar[d] & I_i\tensor J_j\ar[d]  \\
		\Sigma^{-1}\underset{a+b=n}{\dsums}\Gr^a I \tensor \Gr^b J\ar[r]^-\partial & \underset{a+b\geq n+1}{\colim}I_a \tensor J_b\ar[r] & \underset{a+b\geq n}{\colim\thinspace}I_a \tensor J_b
}
\end{align}
for each pair of integers $i,j$ such that $i+j=n$ (see \cite[Lemma 3.3.1]{AKS18}).
The diagram \eqref{eq: d1 diagram} shows that the image of $P_{i,j}$
is in the image of the inclusion
\[ 
	(\Gr^{i+1}I\otimes \Gr^jJ)\oplus (\Gr^{i}I\otimes \Gr^{j+1}J) \to  Gr^{n+1}(I\tensor J)
\]
proving \eqref{eq:d1z-Gr}. 
 Moreover, the diagram
\[ 
	\xymatrix@C=50pt{
		\Sigma^{-1} \Gr^i I \tensor \Gr^j J\ar[rd]_-{\id \tensor \partial+\partial \tensor \id \ \ \ }\ar[r]
			& P_{i,j}\ar[d]\ar[r] & \underset{a+b\geq n+1}{\colim}I_a\tensor J_b\ar[d] \\
		& \atop{\Gr^i I \tensor J^{j+1} \\ \dsum\ I^{i+1}\tensor \Gr^jJ}\ar[r] & \Gr^{n+1}(I\tensor J)
		}
\]
commutes (see \cite[Lemma 4.7]{may-additivity}) and so, writing $\tilde{d}_1^{K}$ for the map of spectra inducing the map $d_1^K$, we have 
\[ 
	\tilde{d}_1^{I\otimes J}=\tilde{d}_1^{I} \tensor \id+\id \tensor \tilde{d}_1^{J}.
\]
Finally, naturality of the map
\[ 
   \pi_s(\Gr^i I)\otimes \pi_s(\Gr^j
J)\longrightarrow \bigoplus_{i+j=n}\pi_s(\Gr^i I \tensor \Gr^jJ )
\]
gives rise to the desired statement, where the sign comes from the natural isomorphism
\begin{align*}
	\Sigma^{-1} \Gr^i I \tensor \Gr^j J&\cong \Gr^i I \tensor \Sigma^{-1}\Gr^j. \qedhere
\end{align*}
\end{proof}

\subsection{Tower functors}
Fix a finite group $G$ and a category $\mathcal{C}\in \CAlg(\Pr^{L,\text{st}})$.

\begin{definition}[{\cite[p. 71]{BHS}}]\label{def:tower}
 A \emph{tower functor in $\mathcal{C}$} is a lax symmetric monoidal functor 
\[
	T\colon \thinspace \mathcal{C} \to  \mathcal{C}^{\Fil}.
\]
\end{definition}

Next we recall two examples of tower functors that come from different families of truncations of $\Sp_G$, namely those arising from the regular slice filtration and the Postnikov filtration.

\begin{definition}[\cite{Ull13}]
	Let $\subslice{G}{n}$ denote the coreflective subcategory of $\Sp_G$ generated by the set $\{{G}_+\wedge_{H}S^{\rho_H k}\}$ where $H$ is a subgroup of $G$, $\rho_H$ is the regular representation of $H$, and $|\rho_H|k\ge n$. 
	The associated functor 
	\[ 
		\tslice{\bullet}\colon \thinspace \Sp_{C_2}\to \Sp_{C_2}^{\Fil}
	\]
	is a tower functor (see \cite[Construction 4.2]{BHS}). Similarly, $\tslice{2\bullet}$ is a tower functor called the even slice tower.  
\end{definition}

\begin{definition}[\cite{HY18}]
	Let $\subpost{G}{n}$ denote the localizing subcategory generated by ${G}_+\wedge_{H} \mathbb{S}^k$ for $H$ a subgroup of $G$ and $k\geq n$.
		The associated functor 
	\[\tpost{\bullet}\colon \thinspace \Sp_{C_2}\to \Sp_{C_2}^{\Fil}\]
	is a tower functor.  Similarly, $\tpost{2\bullet}$ is a tower functor called the even Postnikov tower. This follows because it is a t-structure that is compatible with the symmetric monoidal structure (the full subcategory $\subpost{G}{0}$ clearly contains the symmetric monoidal unit and is closed under the symmetric monoidal product (see \cite[Remark 2.2.1.2]{HA})). 
\end{definition}

The following result will be useful later. 

\begin{lemma} 
\label{lem:post<slice}
There is an inclusion of subcategories
	$\subpost{C_2}{n} \subseteq \subslice{C_2}{n}$. Consequently, there is a natural transformation of tower functors
	\[
		\tpost{\bullet}\implies  \tslice{\bullet}.
	\]
	After postcomposing with the functor $(-)^e\colon \thinspace \Sp_{C_2}^{\Fil}\to \Sp^{\Fil}$, this natural transformation is an equivalence.
\end{lemma}

\begin{proof}
	We check the generators of $\subpost{C_2}{n}$ are in $\subslice{C_2}{n}$. 
	By definition, we have ${C_2}_+\tensor \mathbb{S}^k \in \subslice{C_2}{n},$ for all $k \geq n.$ 
	It suffices to show $\mathbb{S}^k\in \subslice{C_2}{n}$ for all $ k \geq n.$ Since a localizing subcategory is closed under suspension, it suffices to show $\mathbb{S}^n\in \subslice{C_2}{n}$. For $n=1$, this is implied by the cofiber sequence 
	$$C_{2+} \wedge S^1 \to S^1 \to S^\rho.$$
	The case for larger $n$ follows from the multiplicative relation \cite[Cor. 4.2]{Ull13}.
	
    	The underlying equivalence follows from the fact that both filtrations forget to the classical Postnikov filtration on underlying spectra. 
\end{proof}

\begin{remark}
Recall that there is a map $a \colon \thinspace \mathbb{S}^{-\sigma}\to \mathbb{S}$ in $\Sp_{C_2}$ called the Euler class with the property that completion at $a$
\[ 
	(-)_a^{\wedge}\colon \thinspace \Sp_{C_2}\longrightarrow \Sp_{C_2} 
\]
agrees with the cofree functor $F({EC_2}_+,-)$ sending a $C_2$-spectrum $X$ to its underlying Borel $C_2$-spectrum.  We also note that we can identify the essential image of the functor $(-)_a^{\wedge}$ with the category of Borel $C_2$-spectra $\Sp^{BC_2}$.  
\end{remark}

\begin{lemma}\label{lem:acomp_eq}
	Let $f\colon \thinspace X\to Y$ be a map in $\Sp_{C_2}$. If $f$ is an underlying equivalence, then $f$ induces an equivalence between the $a$-completions $f^\wedge_a\colon \thinspace X^\wedge_a \to Y^\wedge_a.$
\end{lemma}

\begin{proof}
Since $a$-completion is equivalent to the cofree construction
\[
	F(E{C_2}_+,-)\colon \thinspace \Sp_{C_2} \to \Sp_{C_2},
\]
it preserves weak equivalences (see, e.g. \cite[\S 2.2.1]{HM17}). The result follows since the essential image of the functor $F(E{C_2}_+,-)=(-)_a^{\wedge}$ is $\Sp^{BC_2}$ whose weak equivalences are detected by the underlying spectra. 
\end{proof}

We therefore have the following corollary of the last two lemmas. 

\begin{proposition}\label{a-complete slice is a-complete Postnikov}
There is a natural equivalence of tower functors
\[
	 (\tpost{\bullet}(-))_a^{\wedge}\overset{\simeq}{\longrightarrow}  (\tslice{\bullet}(-))_a^{\wedge}\colon \thinspace \Sp_{C_2}\to \left (\Sp^{BC_2}\right )^{\Fil}.
\]
\end{proposition}

\begin{proof}
By Lemma \ref{lem:post<slice} and Lemma \ref{lem:acomp_eq} it suffices to check that the functor $(-)_a^{\wedge}\colon \thinspace \Sp_{C_2}^{\Fil}\to \Sp_{C_2}^{\Fil}$ is lax symmetric monoidal. Since $(-)_a^{\wedge}= F(E{C_2}_+,-)\colon \thinspace \Sp_{C_2}^{\Fil}\to (\Sp^{BC_2})^{\Fil}$ is right adjoint to inclusion, it suffices to show that the inclusion functor is symmetric monoidal, which is evident. 
\end{proof}

\subsection{D\'ecalage} \label{sec:decalage}
We described \eqref{eq:cosimpfilt} how to associate a filtered object $\fil_\bullet (X^\bullet)$ to a cosimplicial object $X^\bullet$.  We now describe how a tower functor allows us to associate a different filtered object (the ``$T$-d\'ecalage'' filtration).  

\begin{definition}
Given a tower functor
$T\colon \thinspace \Ccal\to \Ccal^{\Fil}$, and a cosimplicial object $X^\bullet \in \mc{C}^{\Delta}$, the \emph{$T$-d\'ecalage} of $X^\bullet$ defined to be the filtered object given by
\[ 
	\Dec^T(X^\bullet) \coloneqq \Tot T(X^\bullet). 
\]
\end{definition}

For example, if $T=\tau_{\geq \bullet}$, then $\Dec^{T}(X^{\bullet})_w = \Tot_i\tau_{\geq w} X^{i}$.

\begin{notation}\label{notation:decalage}
We refer to the $P_{\ge \bullet}$-d\'ecalage as the \emph{slice d\'ecalage} ($\Dec^{\textup{slice}}$) and the $\tau_{\ge \bullet}$-d\'ecalage as the \emph{Postnikov d\'ecalage} ($\Dec^{\post}$).
We refer to the $P_{\ge 2\bullet}$-d\'ecalage as the \emph{even slice d\'ecalage} ($\Dec^{\evslice}$) and the $\tau_{\ge 2\bullet}$-d\'ecalage as the \emph{even Postnikov d\'ecalage} ($\Dec^{\evpost}$).
\end{notation}

\begin{remark}
The Postnikov d\'ecalage described above corresponds to ``classical d\'ecalage'' on the filtered object $\fil_\bullet X^\bullet$ via the Beilinson $t$-structure.  See \cite[Rmk.~3.7]{GIKR} for a discussion of this point.\\
\end{remark}

Postnikov d\'ecalage has the effect of shifting the pages of the associated spectral sequences.

\begin{theorem}[Hedenlund {\cite[Thm.~II.3.2]{Hedenlund}}]\label{thm:DecSS}
Suppose $X^\bullet$ is a cosimplicial spectrum.  There is a page-shifting isomorphism of spectral sequences
\[ 
	E^{s,n}_{r+1}(X^\bullet) \cong E^{s+n,n}_{r}(\Dec^{\post}(X^\bullet)). 
\]
\end{theorem}

Finally, suppose that $X^\bullet$ is a cosimplicial spectrum with $\pi_*X^s$ concentrated in even degrees.  Then $E^{s,n}_r(X^\bullet)$ is concentrated in degrees where $s+n$ is even, the associated graded of the Postnikov d\'ecalage of $X^\bullet$ is concentrated in even degrees, and we therefore have
\begin{equation}\label{eq:evDecSS}
	E^{s,n}_{2r+1}(X^\bullet) \cong E^{s+n,n}_{2r}(\Dec^{\post}(X^\bullet)) \cong E^{\frac{s+n}{2},n}_r(\Dec^{\textup{ev.post}}(X^\bullet)).
\end{equation}

\subsection{Deformations}
Here we present our deformation of Borel $G$-equivariant stable homotopy theory. We will see in Section \ref{sec: ESSS vs MANSS} that this gives a different presentation of the $a$-complete Artin--Tate real motivic stable homotopy category. We will also demonstrate in Section \ref{sec: computations} that this perspective is useful for computing $RO(G)$-graded homotopy groups of Borel $G$-spectra. First, we need to define a deformation of an object in $\CAlg(\Pr^{L,\text{st}})$.

\begin{notation}
Given a $\mathcal{C}^{\Fil}$-algebra $\mathcal{D}$ in $\CAlg(\Pr^{L,\text{st}})$, we write $\mathbbm{1} _{\mathcal{D}}(i)\coloneqq \mathbbm{1} _{\mathcal{D}}\tensor \mathbbm{1}_{\mathcal{C}}(i)$.
\end{notation}

\begin{definition}[{\cite[Definition C.14]{BHS}}] \label{def:deformation}
A \emph{deformation} of $\mathcal{C}$ is a diagram 
\[
    \begin{tikzcd}
                  &   \mathcal{C}^{\text{def}} \arrow{dr}{\text{Re}} &  \\
                  \mathcal{C} \arrow{ur}{c}\arrow{rr}{\id} && \mathcal{C}
    \end{tikzcd}
\]
in $\CAlg(\Pr^{L,\text{st}})$ such that 
\begin{enumerate}
    \item there is a map 
    \[
    	i\colon \thinspace \mathbb{Z}\to \ker\pi_0\text{Pic}(\text{Re})
    \] 
    of abelian groups sending $n$ to  $\mathbbm{1}_{\mathcal{C}^{\text{def}}}(n)$
    \item the map
    \[
    	\Map_{\mathcal{C}^{\text{def}}}(\mathbbm{1} _{\mathcal{C}^{\text{def}}} (i) , \mathbbm{1} _{\mathcal{C}^{\text{def}}}(j)) \to \Map(\mathbbm{1}_{\mathcal{C}},\mathbbm{1} _{\mathcal{C}})
    \]
    induced by $\text{Re}$ is an equivalence for all integers $i\le j$. 
    \item the category $\mathcal{C}^{\text{def}}$ is rigidly generated by a set $\{c(X_{\alpha})(i)\}$ where $\mathcal{C}$ is rigidly generated by the set $\{X_{\alpha}\}$.
\end{enumerate}
\end{definition}

Finally, we define our deformation of Borel $G$-spectra. Though this definition is currently unmotivated, our project in the next two sections will be to motivate this definition by showing that it recovers a suitable completion of the Artin--Tate $\mathbb{R}$-motivic stable homotopy category (see Section \ref{sec: ESSS vs MANSS}) and that it is useful for computations (see Section \ref{sec: computations}). First, we define an object in $\CAlg(\Sp_G^{\Fil})$. This will involve the $E$-based Amitsur complex
\begin{equation}
	\label{Amitsur complex}
	\begin{tikzcd}
	 E\ar[r,yshift=-0.5ex]\ar[r,yshift=0.5ex] &E^{\tensor 2}\ar[r,yshift=0.7ex]\ar[r,yshift=-0.7ex]\ar[r] &\dots 
\end{tikzcd}
\end{equation}
which is the standard cosimplicial spectrum whose totalization is the $E$-nilpotent completion $\mathbb{S}_{E}^{\wedge}$ and associated Bousfield--Kan spectral sequence can be identified with the $E$-based Adams spectral sequence. 
\begin{remark}\label{decalage as tower functor}
Given an associative monoid $E \in \Sp^{BG}$, then the functor $\Dec^{\post}(-\otimes E^{\otimes \bullet +1})$ given by the Postnikov d\'ecalage of the functor $-\otimes E^{\otimes \bullet +1}$ is again a tower functor (see \cite[Construction C.9]{BHS}).\\
\end{remark}
\begin{notation}
Given an associative monoid $E$ in $G$-spectra, we write $E[G]=E\tensor G_+$ for the group ring, which is again an associative monoid in $G$-spectra.\\ 
\end{notation}
\begin{definition}\label{def:main}
Let $\Sdef{G,p}\coloneqq \Dec^{\textup{ev.post}}(MU_p[G]^{\tensor \bullet+1})$ be the even Postnikov d\'ecalage of the $MU_p[G]$-based Amitsur complex \eqref{Amitsur complex} in the category $\Sp_{ip}^{BG}$.
Define the category of \emph{complete artificial motivic spectra} to be 
\[
	\widehat{\SH}_{G,p}^{\text{art}}\coloneqq \operatorname{Mod}((\Sp^{BG}_{ip})^{\Fil};\Sdef{G,p}).
\]
\end{definition}
\begin{proposition}\label{prop: deformation}
The category $\widehat{\SH}_{G,p}^{\text{art}}$ is a deformation of $\Sp_{ip}^{BG}$. 
\end{proposition}
\begin{proof}
This follows directly from \cite[Proposition C.2]{BHS} and Remark \ref{decalage as tower functor} using $F=\Dec^{\textup{ev.post}}(-\otimes MU_p[G]^{\otimes  \bullet+1})$. The fact that $\mathbb{S}_{G,p}^{\text{def}}$ is a commutative monoid in $(\Sp^{BG}_{ip})^{\Fil}$ follows from the proof of \cite[Theorem 5.64]{PP21}, which generalizes verbatim to $\Sp_{ip}^{BG}$. 
\end{proof}

\section{A tale of two spectral sequences} \label{sec:two-sseq}
Work of the fourth author \cite{Kon20} demonstrated the utility of the $C_2$-effective slice spectral sequence for computing $C_2$-equivariant homotopy groups. In this section, we give a new perspective on the $C_2$-effective slice spectral sequence. Specifically, we construct the modified Adams--Novikov spectral sequence (MANSS), which is useful for computing the homotopy groups of Borel $G$-spectra for any finite group $G$. In the case of $G = C_2$, we compare this to the $C_2$-effective slice spectral sequence.

One advantage of the modified Adams--Novikov spectral sequence is that the differentials can be computed from the classical Adams--Novikov spectral sequence for the underlying spectrum and the equivariant attaching maps for $EG$. When $X$ is a Borel $G$-spectrum with trivial action the differentials in the MANSS can often be deduced from the attaching maps of $BG$ together with the differentials in the Adams--Novikov spectral sequence for the underlying non-equivariant spectrum $X^e$. We illustrate this phenomenon in Section \ref{sec:applications}.

\subsection{The $C_2$-effective slice spectral sequence}\label{sec: esss}
Recall that we have the $C_2$-equivariant Betti realization functor from the $\R$-motivic stable homotopy category:
\[
    \text{Be}_{C_2}\colon \thinspace \SH(\mathbb R)_{i2}\to \Sp_{C_2,i2}.
\]
The $C_2$-effective slice spectral sequence \cite{Kon20} is the spectral sequence associated to the filtration we now define. Suppose $X=\text{Be}_{C_2}(\bar{X})$. There there is an effective slice filtration 
\[ 
	\dots \to \fil_{f+1}^{\text{eff}}\overline{X}\to \fil_{f}^{\text{eff}}\overline{X} \to \dots 
\]
in $\mathbb{R}$-motivic homotopy theory \cite{Voev02} with filtration quotients 
\[ 
	s_f(\overline{X})=\fil_{f}^{\text{eff}}\overline{X}/\fil_{f+1}^{\text{eff}}\overline{X}.
\]
We call the $C_2$-Betti realization of this filtration the \emph{$C_2$-effective slice filtration of $X$ associated to $\overline{X}$}.  The $C_2$-effective slice filtration gives a functor
\begin{align*} 
	i_* : \mr{SH}(\mb{R})_{i2} & \to \mr{Sp}^{\Fil}_{C_2,i2} \\
	\bar{X} & \mapsto \mr{Be}_{C_2}(\fil^\mr{eff}_\bullet X).
\end{align*}
Slice d\'ecalage of the $MU_{\mb{R},2}$-based Amitsur complex gives another functor valued in filtered $C_2$-spectra:
\begin{align*} 
	\Gamma_\mb{R}: \mr{Sp}_{C_2,i2} & \to \mr{Sp}^\mr{Fil}_{C_2,i2}\\
	X & \mapsto \Dec^{\evslice} \left(X\tensor MU_{\R,2}^{\tensor \bullet+1} \right).
\end{align*}

The main theorem of \cite{BHS} is the following.

\begin{theorem}[Galois reconstruction]\label{BHS Galois reconstruction}
The functor
\begin{equation}\label{eq:BHSGalRecon}
	 i_*: \mr{SH}(\mb{R})^{\mr{AT}}_{i2} \to \mr{Mod}(\mr{Sp}^{\mr{Fil}}_{C_2,i2};i_*\mb{S}_{\R,2}) 
 \end{equation}
is an equivalence, and there is an equivalence
\[ 
	i_*\mb{S}_{\R,2} \simeq \Gamma_\mb{R} \mb{S}_{C_2,2}. 
\]
\end{theorem}

It follows that when restricted to $\mr{Sp}_{C_2, i2}$, the functor $\Gamma_\mb{R}$ actually takes values in 
\[ 
	\mr{Mod}(\mr{Sp}^{\mr{Fil}}_{C_2,i2};i_*\mb{S}_{\R,2}) 
\]
and the composite
\[
	\nu_\mb{R} \coloneqq i_*^{-1} \circ \Gamma_\mb{R} : \Sp_{C_2,i2} \to \mr{SH}(\mb{R})^{\mr{AT}}_{i2} 
\]
is a section of the functor $\mr{Be}_{C_2}$.   Given $X \in \mr{Sp}_{C_2, i2}$, the $C_2$-effective slice filtration of $X$ relative to $\nu_\mb{R}(X)$ is the filtered spectrum $\Gamma_\mb{R}X$.  Therefore we shall refer to the filtration $\Gamma_\mb{R}X$ of $X$ simply as \emph{the $C_2$-effective slice filtration of $X$} (without reference to the motivic spectrum $\nu_\mb{R}(X)$ the $C_2$-effective slice filtration is taken with respect to).

The $C_2$-effective slice spectral sequence for $X$ associated to $\overline{X}$ is the $RO(C_2)\times \mathbb{Z}$-graded spectral sequence associated to the filtered $C_2$-spectrum $i_*\overline{X}$
and it has signature
\[
    {}^{C_2\mh \text{ESSS}}\underline{E}_1^{f,V}=\underline{\pi}_{V} \text{Be}_{C_2}(s_f(\bar{X})) \implies \underline{\pi}_{V}X.
\]
In the case where $\bar{X} = \nu_\R(X)$, we shall simply refer to this spectral sequence as the \emph{$C_2$-effective slice spectral sequence of $X$} ($C_2$-ESSS), and the $E_1$-term can be identified with 
\[ 
    \underline{\pi}_{\star}\Gr^\ast(\Gamma_\mb{R}(X)).
\]

One of the most attractive aspects of the $C_2$-ESSS is the computability of its $E_2$-term for the class of $MU_\mb{R}$-\emph{projective spectra}.  Recall the following definition from \cite{BHS}. 

\begin{definition}\label{defn:MURprojective}
A spectrum $X \in \Sp_{C_2,i2}$ is \emph{$MU_{\mb{R}}$-projective} if $MU_{\mb{R},2} \otimes X$ is a retract of a $C_2$-spectrum of the form $\bigoplus_i \Sigma^{n_i \rho} MU_{\mb{R},2}$.\\ 
\end{definition}
For example, the 2-complete sphere is $MU_{\R}$-projective.
The following theorem explains how the $E_2$-term of the $C_2$-ESSS for an $MU_\R$-projective $X \in \Sp_{i2}$ surprisingly only depends on the Adams--Novikov $E_2$-term for the underlying non-equivariant spectrum $X^e \in \Sp_{i2}$. 

\begin{theorem}[Burklund--Hahn--Senger {\cite[Thm.~5.6]{BHS}}]\label{thm:C2ESSSE2}
Suppose that $X \in \Sp_{C_2,i2}$ is $MU_\mb{R}$-projective.  Then the $E_1$-term of its $C_2$-ESSS is given by
\[ 
	{}^{C_2\mh \text{ESSS}}\underline{E}^{f,V}_1(X) = \bigoplus_{s} \Ext^{s,2f}((MU_{2})_*(X^e)\otimes \underline{\pi}_{V-f\rho+s} H\underline{\mb{Z}_{2}}). 
\]
\end{theorem}

\subsection{The modified $E$-based Adams spectral sequence}\label{sec: MANSS}
Given $X \in \Sp^{BG}$, and an associative monoid $E$ in $\Sp^{BG}$ we can form the $E$-based Adams spectral sequence of $X$. Let $\overline{E}$ denote the fiber of the unit map  $\overline{E}=\text{fib}(\mathbb{S}\to E)$. We can then construct the standard $E$-based Adams resolution 
\[ 
    \begin{tikzcd}
   	X \ar[d] & \overline{E} \otimes X \ar[d] \ar[l] &   \overline{E}^{\tensor 2} \otimes X \ar[d] \ar[l] &  \overline{E}^{\tensor 3} \otimes X \ar[l] \ar[d] & \ar[l] \dots \\
   	E \otimes X &   E\tensor \overline{E}  \otimes X&     E\tensor \overline{E}^{\tensor 2} \otimes X&  E\tensor \overline{E}^{\tensor 2} \otimes X & 
    \end{tikzcd}
\]
of $X$. We view this as a filtered spectrum 
\[ 
	\Gamma_{E}X \coloneqq \{ \dots \to \bar{E}^{\tensor 2} \otimes X \to \bar{E} \otimes X \to  X =  X = \dots  \} \in (\Sp^{BG})^{\mr{Fil}} 
\]
Let 
\[ 
	X^{\wedge}_E \coloneqq \Tot(E^{\bullet+1} \otimes X) 
\]
be the totalization of the $E$-based Amitsur complex tensored with $X$ \eqref{Amitsur complex}.
If the map $X \to X^\wedge_E$ 
is an equivalence, we have 
\[ 
	\lim_s \Gamma_E(X)_s=0. 
\]
Of course, we always have 
\[ 
	\colim_s \Gamma_E(X)_s = X.
\]

Next, we will set up the homotopy fixed point spectral sequence so that it also conditionally converges to the colimit. Define 
\[ 
	EG^{(k)}\coloneqq|\text{sk}_kB_{\bullet}(G,G,*)| 
\]
where $B_{\bullet}(G,G,*)$ is the usual two-sided bar construction whose geometric realization is $EG$. Let $EG^{(-1)}=\emptyset$ by convention. The homotopy fixed point spectral sequence of $X \in Sp^{BG}$ is the spectral sequence associated to the tower
\[ 
	F(EG^{(0)}_+,X) \leftarrow F(EG^{(1)}_+,X) \leftarrow F(EG^{(2)}_+,X) \leftarrow \cdots 
\]
whose inverse limit in $\Sp^{BG}$ is $X$.

We will explicitly describe the decreasing filtration of $X$ which gives rise to this tower.
Write $\widetilde{E}G^{(k)}$ for the cofiber of the collapse map 
\[ 
	EG^{(k)}_+\to S^0. 
\]
Then there is a cofiber sequence
\[ 
	F(\widetilde{E}G^{(s)},X) \rightarrow X \rightarrow F(EG^{(s)}_+,X). 
\]
It follows that the homotopy fixed points spectral sequence computing
$\underline{\pi}_\star(X)$ arises from the following decreasing filtration of $X$:
\begin{equation}\label{eq:hfpss-tower}
    \begin{tikzcd}
        X\ar[d] & \ar[l] F( \widetilde{E}G^{(0)}, X)\ar[d] &  \ar[l] F(\widetilde{E}G^{(1)},X) \ar[d] & \ar[l]  \dots \\ 
        F( EG^{(0)}_+,X) & F( EG^{(1)}_+/EG^{(0)}_+,X) & F(EG^{(2)}_+/EG^{(1)}_+,X)& 
    \end{tikzcd}
\end{equation}
Here, the filtration quotients in this diagram are computed using the diagram of cofiber sequences 
\begin{equation}\label{diagram of cof sequences for EG}
    \begin{tikzcd}
    EG_+^{(s-1)} \ar[r]\ar[d] & \mathbb{S} \ar[r] \ar[d] & \widetilde{E}G^{(s-1)} \ar[d] \\
     EG_+^{(s)} \ar[r] \ar[d] & \mathbb{S} \ar[r] \ar[d] & \widetilde{E}G^{(s)} \ar[d] \\
        EG_+^{(s)}/ EG_+^{(s-1)} \ar[r] & 0 \ar[r] &  \widetilde{E}G^{(s)}/\widetilde{E}G^{(s-1)}\\
    \end{tikzcd}
\end{equation}
which identifies $\Sigma^{-1} \widetilde{E}G^{(s)}/\widetilde{E}G^{(s-1)}$ with $EG^{(s)}_+/EG^{(s-1)}_+$ for $s\ge 0$.  We regard this construction as a functor 
\[
	\Gamma^{h}\colon \thinspace \mathbb{Z}^{\op}\to \Sp^{BG}. 
\]
Specifically, set 
\[ 
	\Gamma^{h}(X)_s \coloneqq\begin{cases} 
						X & \text{ for } s \le  0 \\   
						F(\widetilde{E}G_+^{(s-1)},X)&  \text{ for } s > 0.
					\end{cases}
\]
This has been arranged so that $\lim_s \Gamma^{h}(X)_s =0$ 
and $\colim_s \Gamma^h(X)_s=X$. To see this, note that since the map
\[ 
	EG_+ \longrightarrow S^0 
\]
is a non-equivariant equivalence, the space $\widetilde{E}G$ is non-equivariantly contractible.

Now we mix these two filtrations in the following way. By tensoring the Adams resolution for $\mathbb{S}$ with $X$, we produce the $E$-Adams resolution for $X$ in $\Sp^{BG}$. If we also filter $X$ by the filtration $\Gamma^h(X)$ as above then we produce an object 
\[ 
	X_{\bullet,\sbsq}\in \text{Fun}(\mathbb{Z}^{\op}\times \mathbb{Z}^{\op},\Sp^{BG})
\]
defined on an object $(i,j)$ by 
\begin{equation}\label{eq: Xij}
X_{i,j}=\begin{cases} 
                    \overline{E}^{\tensor i} \tensor X & \text{ for } i >0,j \le 0 \\
                    X & \text{ for } i\le 0,j\le 0 \\
                       \overline{E}^{\tensor i} \tensor F( \widetilde{E}G_+^{(j-1)},X) & \text{ for } i > 0, j > 0 \\
                 F(\widetilde{E}G_+^{(j-1)},X)  & \text{ for } i\le 0,j> 0. \\
                \end{cases}
\end{equation} 
Note that $X_{\bullet,\sbsq}$ is simply the exterior tensor product $\Gamma_E(\mb{S})_\bullet  \overline{\tensor}~~\Gamma^h(X)_{\sbsq}$. 

We then consider the Day convolution 
\[ 
	\Gamma^h_E(X)_{\bullet} \coloneqq \Gamma_E(\mb{S})_{\bullet}\tensor \Gamma^h(X)_\bullet.
\]
By definition, $\Gamma^h_E(X)_{\bullet}$ is the left Kan extension in the diagram 
\[
    \begin{tikzcd}
        \mathbb{Z}^{\op}\times \mathbb{Z}^{\op} \ar[d,"+",swap] \ar[r,"X_{\bullet,\sbsq}"] &  \Sp^{BG} \\
        \mathbb{Z}^{\op}\ar[ur,"\Gamma^h_E(X)_{\bullet}", swap] & 
    \end{tikzcd}
\]
and it is explicitly defined on objects by 
\[ 
	\Gamma^h_E(X)_s=(\Gamma_E(\mb{S})_{\bullet}\tensor \Gamma^h(X)_{\bullet})_s=\underset{a+b\ge s}{\colim}\thinspace \Gamma_{E}(\mb{S})_a\tensor \Gamma^{h}(X)_b.
\]
Unpacking this, we compute that 
\begin{equation}\label{eq:GammahEX}
    \Gamma^h_E(X)_{s}=\begin{cases}
     X & \text{ if } s  \le 0\\ 
                   \underset{a+b\ge  s}{\colim}\thinspace \overline{E}^{\tensor a} \tensor
                   F(\widetilde{E}G_+^{(b-1)},X) &                                  \text{ if } s >0 .
     \end{cases}
\end{equation}
Let $\mathbb{N}^{\op}$ denote the full subcategory of $\mathbb{Z}^{\op}$ generated by the nonnegative integers and let $\mathbb{N}^{\delta}$ denote the category whose objects are nonnegative integers and whose only morphisms are the identity maps. By Lemma \ref{lemma:E1} we can identify
\begin{equation}\label{eq:MANSS-E1} 
	\text{cofib}(\Gamma^h_{E}(X)_{s+1}\to \Gamma^h_E(X)_{s})=\bigoplus_{a+b=s}E \otimes \overline{E}^{\tensor a} \tensor F(EG^{(b)}_+/EG^{(b-1)}_+, X).
\end{equation}

\begin{notation}
For $\star\in RO(G)$ we write $\star=*+\diamond$ where $*$ is the virtual dimension of $\star$ (representing a sum of trivial representations), and $\diamond$ is a virtual representation of virtual dimension $0$.\\
\end{notation}

Using this notation, we have a $ \text{RO}(G)\times \mathbb{Z}$-graded spectral sequence with signature 
\[{}^{E \mh \textup{MASS}}\underline{E}_1^{s,n+\diamond}(X)=\bigoplus_{a+b=s}\underline{\pi}_{n+\diamond} (E \otimes \overline{E}^{\tensor a} \tensor F(EG^{(b)}_+/EG^{(b-1)}_+, X)) \implies \underline{\pi}_{n+\diamond}(X_{E}^{h}) \]
where 
\begin{equation}\label{eq: notation for Borel E-completion} 
	X^h_E \coloneqq \lim_s X/\Gamma^h_E(X)_s. 
\end{equation}
We shall refer to this spectral sequence as the \emph{modified $E$-based Adams spectral sequence} ($E$-MASS) for $X$. If $X\in \Sp_{G}$, then abuse terminology and refer to the $E$-MASS for $X^h$ as the $E$-MASS for $X$. 

\subsection{The $E$-MASS as a generalized Adams spectral sequence}\label{sec: MANSS as ANSS}
The following lemma will be useful for identifying certain filtration quotients. 
\begin{lemma} \label{lem:EG-quotient}
There is an equivalence 
\[ 
	F(E{G}^{(k)}_+/E{G}^{(k-1)}_+,\mathbb{S}) \hteq \mathbb{S}[G]\tensor\overline{\mathbb{S}[G]}^{\tensor k}
\]
where $\overline{\mathbb{S}[G]} \coloneqq \fiber(\mathbb{S}\to \mathbb{S}[G])$. 
\end{lemma}
\begin{proof}
By definition, $EG^{(k)} \hteq |\text{sk}_kB_\bullet(G,G,*)|$. 
Consequently, there is an equivalence
\[ 
	F(EG^{(k)}_+,\mathbb{S}) \hteq \Tot_k(F(B_\bullet(G,G,*)_+, \mathbb{S})). 
\]
By \cite[Proposition 2.14]{MNN17}, there is therefore an equivalence 
 \[ 
 	F(EG^{(k)}_+,\mathbb{S}) \simeq \cofib(\overline{\mathbb{S}[G]}^{\tensor k+1}\to \mathbb{S}).
\]

The lemma then follows from the diagram of fiber sequences 
\[
    \xymatrix{
        \overline{\mathbb{S}[G]}^{\tensor k+1}\ar[r]\ar[d] & \mathbb{S}\ar@{=}[d]\ar[r] & F(EG_+^{(k)},\mathbb{S})\ar[d] \\
	\overline{\mathbb{S}[G]}^{\tensor k}\ar[r]\ar[d] & \mathbb{S}\ar[r]\ar[d] & F(EG_+^{(k-1)},\mathbb{S})\ar[d]\\
	\overline{\mathbb{S}[G]}^{\tensor k}\tensor \mathbb{S}[G]\ar[r] &  0 \ar[r] & F(\Sigma^{-1} EG_+^{(k)}/EG_+^{(k-1)},\mathbb{S}).
	}
\]
\end{proof}
\begin{remark}
When $G=C_2$, we have $\overline{\mathbb{S}[C_2]}=\mathbb{S}^{-\sigma}\coloneqq F(\mathbb{S}^{\sigma},\mathbb{S})$. To see this, note that there is a fiber sequence 
\[
	\mathbb{S}^{-\sigma}\too{a}\mathbb{S}^0 \to Ca =\mathbb{S}[C_2]
\] 
so $\overline{\mathbb{S}[C_2]} = \mathbb{S}^{-\sigma},$ and the
map $\overline{\mathbb{S}[C_2]}^{\tensor k} = \mathbb{S}^{-k\sigma}\to \mathbb{S}$ is simply
$a^k$. Consequently, by Lemma \ref{lem:EG-quotient} we have 
 equivalences
\begin{align*}
	E{C_2}^{(k)}_+/E{C_2}^{(k-1)}_+ & \hteq \Sigma^{k\sigma}Ca
\end{align*}
for each $k\ge 0$. 
\end{remark}
Now we prove the main result from this section. 
\begin{proposition}[The $E$-MASS is a generalized Adams spectral sequence]\label{prop:MANSS=NE-Adams}
The modified $E$-based Adams spectral sequence for a finite group $G$
agrees with the $E[G]$-based Adams spectral sequence from the $E_2$-term onwards.
\end{proposition}
\begin{proof}
To prove this result, it suffices to prove that the filtered spectrum that we use to construct the modified $E$-based Adams spectral sequence is in fact an $E[G]$-based Adams resolution. 

Since $\widetilde{E}G^{(s)}$ is finite, the modified $E$-based Adams spectral sequence for $X$ is the spectral sequence associated to the Day convolution 
\[ 
	\Gamma^h_E(X)_\bullet \simeq \Gamma_E(\mb{S})_\bullet \otimes \Gamma^h(\mb{S})_\bullet \otimes X. 
\]
(see Lemma \ref{lem:EG-quotient}). 
We therefore wish to show that 
\[       
	 \dots  \to  (\Gamma_E(\mb{S})_\bullet \otimes \Gamma^h(\mb{S})_\bullet \otimes X)_1 \to  (\Gamma_E(\mb{S})_\bullet \otimes \Gamma^h(\mb{S})_\bullet \otimes X)_0 =  X 
\]
is an $E[G]$-injective resolution of $X$ in the sense of \cite{Mil81}. We need to check that
\begin{enumerate}
    \item there is an equivalence $\Gr^n(\Gamma_E(\mb{S})_\bullet \otimes \Gamma^h(\mb{S})_\bullet \otimes X)\simeq E[G]\tensor M$ for some spectrum $M$ and consequently, $\Gr^n(\Gamma_E(\mb{S})_\bullet \otimes \Gamma^h(\mb{S})_\bullet \otimes X)$ is $E[G]$-injective;
    \item the map 
    \[ 
        (\Gamma_E(\mb{S})_\bullet \otimes \Gamma^h(\mb{S})_\bullet \otimes X)_n\to \Gr^n(\Gamma_E(\mb{S})_\bullet \otimes \Gamma^h(\mb{S})_\bullet \otimes X)
    \]
    induces a monomorphism after applying $E[G]_*$ for all $n\ge 0$ (note that the case $n=0$ is clear because $(\Gamma_E(\mb{S})_\bullet \otimes \Gamma^h(\mb{S})_\bullet \otimes X)_1=\overline{E[G]}$).
\end{enumerate}
The first part, follows from \eqref{eq:MANSS-E1} since 
\begin{align}
  \bigoplus_{a+b=s} E\tensor \overline{E}^{\otimes a}\tensor \mathbb{S}[G]\tensor \overline{\mathbb{S}[G]}^{\tensor b} \tensor X 
   \simeq  \thinspace & E[G] \otimes ( \bigoplus_{a+b=s} \overline{E}^{\otimes a}\tensor \overline{\mathbb{S}[G]}^{\tensor b}) \otimes X
\end{align}
for each $s\ge 0$. 

To prove the second part, it suffices to show that the induced map
\[  
	E[G]_*((\Gamma_E(\mb{S})_\bullet \otimes \Gamma^h(\mb{S})_\bullet \otimes X)_{n+1})\to E[G]_*((\Gamma_E(\mb{S})_\bullet \otimes \Gamma^h(\mb{S})_\bullet \otimes X)_n)
\]
is the zero map.  To prove this, we first observe that the source and target can be identified with $E[G]_*$ applied a finite colimit. It therefore suffices to show that the maps 
\[ 
	E[G]_*(\overline{E}^{\otimes i} \otimes \overline{\mathbb{S}[G]}^{\otimes j})\to E[G]_*(\overline{E}^{\otimes i}\otimes \overline{\mathbb{S}[G]}^{\otimes j-1})
\]
and 
\[ 
	E[G]_*(\overline{E}^{\otimes i} \otimes \overline{\mathbb{S}[G]}^{\otimes j})\to E[G]_*(\overline{E}^{\otimes i-1} \otimes \overline{\mathbb{S}[G]}^{\otimes j})
\]
are the zero maps, but this is the case because $E^{\otimes i}\otimes \mathbb{S}[G]^{\otimes \bullet}$ is an $\mathbb{S}[G]$-Adams resolution for each $i$ and
$E^{\otimes \bullet}\otimes \mathbb{S}[G]^{\otimes j}$ is an $E$-Adams resolution for each $j$. 
\end{proof}

\begin{corollary}
The $E$-MASS for $X \in \Sp^{BG}$ conditionally converges if the underlying nonequivariant spectrum $X^e$ is $E^e$-nilpotent complete
\end{corollary}

\begin{proof}
It suffices to show that the map
$$ X \to \Tot E[G]^{\otimes \bullet + 1}\otimes X $$
is a Borel equivalence, which amounts to showing that it is an underlying equivalence.  This follows from the hypothesis that $X^e$ is $E^e$-nilpotent complete and the fact that 
$$ E[G]^e \simeq \bigoplus_G E^e. $$
\end{proof}

\begin{corollary}\label{cor:mult E}
If $X \in \Sp_{G}$ is a ring spectrum, the $E$-MASS for $X$ is a spectral sequence of algebras, and in general the $E$-MASS is a spectral sequence of modules over the $E$-MASS for $\mb{S}$.
\end{corollary}

\begin{proof}
This follows from the fact that $E[G] \simeq F(\Sigma^{\infty}G_+, E)$.
\end{proof}

\begin{corollary}\label{thm:MANSS=E/a-Adams}
The $E$-MASS for $G=C_2$ agrees with the $E/a$-based Adams spectral sequence from the $E_2$-term onwards.
\end{corollary}
\begin{proof}
This follows from the observation that $Ca\simeq \mathbb{S}[C_2]$ and Proposition \ref{prop:MANSS=NE-Adams}. 
\end{proof}

\subsection{The MANSS and the $C_2$-effective slice spectral sequence}\label{sec: ESSS vs MANSS}
We now introduce the main tool that we use this paper.
\begin{definition}\label{def: MANSS}
In the special case of Section \eqref{sec: MANSS as ANSS} where $E$ is the Borel spectrum $MU_p$ with trivial action,\footnote{Here we work in the category $\Sp^{BG}_{ip}$ so that all instances of $\otimes$ should be interpreted as $\otimes_{\mathbb{S}_{p}}$} we shall refer to the $MU_p$-MASS as the $p$-primary \emph{modified Adams--Novikov spectral sequence} (MANSS)
\[ 
	{}^{\MANSS}\underline{E}_1^{s,n+\diamond}(X)=\bigoplus_{a+b=s}\underline{\pi}_{n+\diamond} (MU_p \otimes \overline{MU}_p^{\tensor a} \tensor F(EG^{(b)}_+/EG^{(b-1)}_+, X)) \implies \underline{\pi}_{\star}(X_{MU_p}^{h}) 
\]
where $X_{MU_p}^{h}$ is defined as in \eqref{eq: notation for Borel E-completion}. 
\end{definition}

\noindent The following corollaries follow from the evident generalization of Proposition \ref{prop:MANSS=NE-Adams} to the category $\Sp^{BG}_{ip}$. 

\begin{corollary}[MANSS is $MU_p/a$-based ASS]\label{thm:MANSS=MU/a-Adams}
The $p$-primary MANSS for $G=C_2$ agrees with the ${MU}_{p}/a$-based Adams spectral sequence from the $E_2$-term onwards.
\end{corollary}
\begin{proof}{
This follows from Corollary~\ref{thm:MANSS=E/a-Adams} by taking $E=MU_p$. }
\end{proof}

\begin{corollary}\label{cor:mult}
If $X \in \Sp_{G,ip}$ is a ring spectrum, the $p$-primary MANSS for $X$ is a spectral sequence of algebras, and in general the $p$-primary MANSS is a spectral sequence of modules over the $p$-primary {MANSS} for $\mb{S}_p$.
\end{corollary}

\begin{proof}
This follows from Corollary~\ref{cor:mult E} by taking $E=MU_p$.
\end{proof}

The following corollary of Proposition \ref{prop:MANSS=NE-Adams} begins to explain the motivation behind Definition~\ref{def:main}.

\begin{corollary}
In the category $\Sp^{\Fil}_{G,ip}$ there is an equivalence
{\[ 
	\Sdef{G,p}\coloneqq \Dec^{\textup{ev.post}}(MU_p[G]^{\tensor \bullet+1}) \hteq \Dec^{\textup{ev.post}} \Gamma^h_{MU_p}\mb{S}_p. 
\]}
\end{corollary}

We already showed that the $2$-primary MANSS for the group $G=C_2$ agrees with the $MU_2/a$-based Adams spectral sequence in Corollary \ref{thm:MANSS=MU/a-Adams}. Therefore to compare the MANSS with the $C_2$-ESSS it suffices to compare the $MU_2/a$-based Adams spectral sequence to the $C_2$-ESSS. Since the MANSS is a spectral sequence in Borel $G$-spectra, it corresponds to an $a$-complete $p$-complete filtered object in the $C_2$-equivariant category under the equivalence
\[ 
	(\Sp_{C_2,ip})^\wedge_a \simeq \Sp_{ip}^{BG}.
\]
For the rest of this section we will work in the category $\Sp_{C_2,ip}$, and objects of $\Sp_{ip}^{BG}$ will be regarded as $a$-complete objects of $\Sp_{C_2,ip}$.

\begin{lemma}
\label{lem:moveYinside}
Let $X\in\Sp_{C_2,i2}$. 
For any $i\geq 1,$ there is a natural equivalence 
\begin{align}\label{eq:moveYinside}
(\tpost {n} X)\otimes (Ca)^{\tensor i} \to \tpost n (X\otimes (Ca)^{\tensor i}).
\end{align}
\end{lemma}

\begin{proof}
Since $\tpost n $ is closed under the operation $-\otimes Ca$, the existence of a map follows from the fact that $\tpost n$ is part of the data of a $t$-structure on $\Sp^{BC_2}$. 

For any $C_2$-spectrum $Z$, we have
\[
	\pi_k^e (\tpost {n} (Z\tensor Ca) )=\begin{cases}
								\pi^{e}_k (Z\tensor Ca), ~k\geq n \\
								0, \text{~otherwise}
								\end{cases},
\text{ and ~~}
	\underline \pi_k^{C_2} (\tau_{\ge n}(Z\tensor Ca) )=\begin{cases} 
	               									 \pi_k^e(Z), ~k\geq n \\
	               								 	  0, \text{~otherwise.}
	               								 	\end{cases}
\]
By induction, we compute the right hand side of \eqref{eq:moveYinside} 
\[
	\underline \pi_k \Big(\tpost {n} (X\tensor (Ca)^{\tensor i})\Big)=\begin{cases}
														\underline \pi_k ((\tpost {n}X)\tensor (Ca)^{\tensor i}), ~k\geq n = \begin{cases}
\pi_{k}^{e} (X\tensor (Ca)^{\tensor i-1}), \text{~at~} C_2/C_2,\\
\pi_{k}^e (X\tensor (Ca)^{\tensor i}), \text{~at~} C_2/e,
\end{cases} \\
0, \text{~otherwise}.
\end{cases}
\]
The Mackey functor valued homotopy of the left hand side equals
\[
\underline \pi_k \Big((\tpost {n} X)\tensor (Ca)^{\tensor i}\Big)=\begin{cases}
\pi_{k}^{e} \Big((\tau_{\ge n}X)\tensor (Ca)^{\tensor i-1}\Big), \text{~for~} C_2/C_2= \begin{cases}
\pi_{k}^e (X\tensor (Ca)^{\tensor i-1}), ~k\geq n\\
0, \text{~otherwise},
\end{cases}\\
\pi_{k}^e \Big((\tpost {n}X)\tensor (Ca)^{\tensor i}\Big), \text{~for~} C_2/e= \begin{cases}
\pi_{k}^e (X\tensor (Ca)^{\tensor i}), ~k\geq n\\
0, \text{~otherwise}.
\end{cases}
\end{cases}
\]
Since Mackey functor valued homotopy detects weak equivalences and these identifications are compatible with the Mackey functor structure maps, the result follows.
\end{proof}

\begin{theorem}[Comparison of $C_2$-ESSS and {$MU_2/a$-based ASS}]\label{ESSS is MANSS}
Let $X\in \Sp_{C_2,i2}$. There is an equivalence of filtered objects in $\Sp_{C_2,i2}$
\[
	\Big(\Dec^{\evslice}(MU_{\bR,2}^{\tensor \bullet+1}\tensor X)\Big)^\wedge_a\simeq \Dec^{\textup{ev.post}}\Big( \big((MU_2/a)^{\tensor \bullet+1} \tensor X\big) \Big).
\]
\end{theorem}
\begin{proof}
This follows by the composite of natural equivalences
\begin{align*}
	\Big( \lim_i \tslice {2n}(MU_{\bR,2}^{\tensor i+1}\tensor X)\Big)^\wedge_a & \stackrel{\text{commute limits}}{\simeq} \lim_i \Big( \big(\tslice{2n} (MU_{\bR,2}^{\tensor i +1} \tensor X)\big)^\wedge_a \Big)	\\
	& \stackrel{\text{Lem.} \ref{lem:post<slice}+ \text{Lem.} \ref{lem:acomp_eq}}{\simeq}  \lim_i  \Big( \big(\tpost{2n} (MU_{\bR,2}^{\tensor i+1} \tensor X)\big)^\wedge_a \Big)\\
	& \stackrel{\text{def}}{\simeq} \lim_j\lim_i \Big(\tpost{2n}(MU_{\bR,2}^{\tensor i+1} \tensor X)\tensor Ca^{\tensor j+1} \Big)\\
	& \stackrel{\text{diagonal}}{\simeq} \lim_i \Big(\tpost{2n}(MU_{\bR,2}^{\tensor i+1} \tensor X)\tensor Ca^{\tensor i+1} \Big)\\
	&\stackrel{\text{Lem.} \ref{lem:moveYinside}}{\simeq} \lim_i \Big(\tpost{2n}\big((MU_{\bR,2}/a) ^{\tensor i+1} \tensor X\big) \Big)\\
	&\stackrel{\text{shearing}}{\simeq} \lim_i \Big( \tpost{2n} \big((MU_2/a)^{\tensor i+1} \tensor X\big) \Big). 
\end{align*}
Here the shearing isomorphism is the equivalence
\[ 
	MU_{\R,2}/a \simeq MU_{2}/a 
\]
which follows from the fact that $MU_{2}$ and $MU_{\R,2}$ have equivalent underlying ring spectra.
The proof is then complete by Corollary \ref{thm:MANSS=MU/a-Adams}.
\end{proof}

Associated to the $a$-completion of the $C_2$-effective slice filtration is the $a$-\emph{completed $C_2$-effective slice spectral sequence}
\[ 
	{}^{C_2\mh \textup{ESSS}}\underline{E}^{f,n+\diamond}_1(X)^{\wedge}_a = \underline{\pi}_{n+\diamond}(\Gr^f \Gamma_\mb{R} X) \implies \underline{\pi}_{n+\diamond}(X^h_{MU}). 
\]
The following corollary follows from \eqref{eq:evDecSS} (note that the hypothesis that $MU_*(X^e)$ is even implies that the MANSS has no nontrivial $d_r$ differentials for $r$ even).

\begin{corollary}\label{cor:ESSS is MANSS}
Suppose that $(MU_{2})_*(X^e)$ is concentrated in even degrees.  Then for $r \ge 1$, there is a (page shifting) isomorphism of spectral sequences
\[ 
	{}^{\MANSS}\underline{E}^{s,n+\diamond}_{2r+1}(X) \cong {}^{C_2\mh \textup{ESSS}}\underline{E}^{\frac{n+s}{2},n+\diamond}_r(X)^{\wedge}_a. 
\]
\end{corollary}
We can now give the first evidence for why the category of artificial complete motivic spectra is a useful deformation of Borel $G$-spectra. 
\begin{theorem}[MANSS Galois reconstruction]\label{thm: MANSS Galois reconstruction}
{There are equivalences
\[ 
\widehat{\SH}_{C_2,2}^{\text{art}}\coloneqq 
\operatorname{Mod}((\Sp^{BC_2}_{i2})^{\Fil};\Sdef{G,2})
\hteq	\Mod((\Sp^{BG}_{ip})^{\Fil}; \Dec^{\textup{ev.post}} \Gamma^h_{MU_2}\mb{S}_2) \hteq (\SH(\R)^{\textup{AT}}_{i2})_{a}^{\wedge}. 
\]}
\end{theorem}
\begin{proof}
This follows from Proposition \ref{a-complete slice is a-complete Postnikov}, Theorem \ref{ESSS is MANSS}, and Theorem \ref{BHS Galois reconstruction}. 
\end{proof}

\section{Computing with the MANSS}\label{sec: computations}
We compute the $E_2$-term of ($p$-primary) MANSS (see Definition \ref{def: MANSS}) using double complex spectral sequences in Section \ref{sec: E2 MANSS}.
We then specialize to the case $G = C_2$, $p=2$, $G/H = C_2/C_2$ and write
$$ {}^{\MANSS}E_r(X) := {}^{\MANSS}\underline{E}_r(X)(C_2/C_2). $$
We express the $E_2$-term using $A^{s,2t}\coloneqq\Ext^{s,2t}((MU_{p})_*X^e)$ in two cases: 
\begin{enumerate}
	\item when $X$ is $MU_{\mathbb R}$-projective (Section \ref{sec: MUR proj case}), and
	\item when $X$ has trivial $C_2$-action and $(MU_{2})_*X^e$ is torsion free (Section \ref{sec: triv action case}).
\end{enumerate}
In the first case, we show the MANSS relates to the ANSS and the $C_2$-effective SS; in the second case, by studying a motivic version of MANSS, we relate it to the Atiyah--Hirzebruch spectral sequence.

Assume $x\in A^{s,2t}$ (see Notation \ref{not:braxangle}). In this section, we introduce the following notations.
\begin{center}
\begin{tabular}{ |c|c|c|c| } 
	\hline
	Notation & Degree & Definition & From\\
	\hline
	$\bar{x}$ & $(t,t\rho-s)$  & Notation \ref{not:xbar} & $C_2$-effective SS\\ 
	$\bar{ux}$ & $(t,t\rho-s+1-\sigma)$ & Notation \ref{not:xbar} & $C_2$-effective SS\\ 
 	$\bra{x}$  & $(s,t\rho-s)$ & Notation \ref{not:braxangle} & double complex SS\\ 
 	$\bra{ux}$  & $(s,t\rho-s)$ & Notation \ref{not:braxangle} & double complex SS\\ 
  	$[x]$  & $(s,2t-s)$ & Notation \ref{not:brax} & compare AHSS\\
    	$[ux]$ &  $(s,2t-s+1-\sigma,t)$ & Notation \ref{not:brax} & compare AHSS\\
 	$[x]$ & $(s,2t-s,t,t)$  & Proposition \ref{prop:motMANSSE2} &  motivic MANSS\\ 
 	$[ux]$ & $(s,2t-s+1-\sigma,t)$ & Proposition \ref{prop:motMANSSE2} & motivic MANSS\\ 
	\hline
\end{tabular}
\end{center}

\subsection{The $E_2$-term of the MANSS}\label{sec: E2 MANSS}
Let $X$ be an object of $\Sp_{ip}^{BG}$. Recall from Section \ref{sec: MANSS} that the $E_1$-term of the $p$-primary MANSS is
\[
	{}^{\MANSS}\underline{E}_1^{f,V}(X)=\bigoplus_{a+b=f}\underline{\pi}_{V}Z_{a,b}
\] 
where 
\[Z_{a,b}=F(EG^{(a)}/EG^{(a-1)}_+, X)\tensor MU_p \tensor\overline{MU}_p^{\tensor b}.\] 
In other words, the $E_1$-term $(\underline{E}_1^{\star,\ast}(X),d_1^{\MANSS})$ is the total complex of the double complex $\bigoplus_{a,b}\underline{\pi}_\star(Z_{a,b})$ whose horizontal differential 
$d_1^\ANSS\colon\thinspace \underline{\pi}_{\star}Z_{a,b}\to \underline{\pi}_{\star}Z_{a,b+1}$ is the Adams--Novikov $d_1$ and vertical differential $d_1^\hfpss\colon \thinspace \underline{\pi}_{\star}Z_{a,b}\to \underline{\pi}_{\star}Z_{a+1,b}$ is the homotopy fixed point spectral sequence $d_1$. We wish to compute the
cohomology $E_2^{\star,*}$ of the total complex.
There are two double complex spectral sequences 
\begin{align}
	\label{eq:double-complex1} {}_{I}\underline{E}_2^{a,b,\star}(X)= H^b(H^a(\dsums_{i,j} \underline{\pi}_\star(Z_{i,j}), d_1^\hfpss), d_1^\ANSS) \Longrightarrow {}^{\MANSS} \underline{E}_2^{a+b,\star}(X) \\
	\label{eq:double-complex2} {}_{\II}\underline{E}_2^{a,b,\star}(X)= H^a(H^b(\dsums_{i,j} \underline{\pi}_\star(Z_{i,j}), d_1^\ANSS), d_1^\hfpss) \Longrightarrow {}^{\MANSS} \underline{E}_2^{a+b,\star}(X)
\end{align}
which come from filtering the double complex by $i$ and $j$, respectively \cite[\S2.4]{McC01}.
We use these spectral sequences to describe the $E_2$-term of the MANSS. 

\begin{proposition}\label{prop:sseq-I}
There is an identification
\[
	{}_I\underline{E}_2^{a,b,n+\diamond}(X)(G/H) = \Ext^{b,a+b+n}(H^a(H;(MU_{p})_*(\mathbb{S}^{-i_H^*\diamond}\tensor X))) .
\]
In particular, 
the value of the $E_2$-term of the MANSS at $G/e$ is 
\[ 
	{}^{\MANSS}\underline{E}_2^{b,V}(X)(G/e) =\Ext^{b,b+|V|}((MU_{p})_*(X^e)).
\]
\end{proposition}

\begin{proof}
Since $EG^{(i)}/EG^{(i-1)}$ is dualizable, $(MU_{p})_*\bar{MU_{p}}$ is flat over $(MU_{p})_*$, and
$MU_p$ and $\overline{MU}_p$ have trivial $G$-action, we have:
\begin{align}
	\notag \underline{\pi}_{*+\diamond}(Z_{i,j}) & \cong \underline{\pi}_*(F(EG^{(i)}/EG^{(i-1)}, \mathbb{S}^{-\diamond}\otimes X)\tensor MU_{p}\tensor \overline{MU}_{p}^{\tensor j}) \\
	\label{eq:Z} & \isom \underline{\pi}_*(F(EG^{(i)}/EG^{(i-1)}, \mathbb{S}^{-\diamond}\tensor  X\tensor MU_p))\tensor_{(MU_{p})_*}(MU_{p})_*\overline{MU}_{p})^{\tensor_{(MU_{p})_*} j}.
\end{align}
At the $G$-set $G/H$, this is $(MU_{p})_*\overline{MU}_{p}^{\tensor j}$ tensored with the $E_2$-term of a
homotopy fixed points spectral sequence, and so we have
\[ 
	H^a(\underline{\pi}_{*+\diamond}(Z_{\bullet,j})(G/H),d_1^\hfpss) \isom H^a(H;MU_*(\mathbb{S}^{-i_H^*\diamond}\tensor X))\tensor_{(MU_{p})_*}(MU_{p})_*\overline{MU}_{p}^{\tensor_{(MU_{p})_*} j}. 
\]
This is isomorphic to the cobar complex 
\[ 
	C^{j,*}_{(MU_{p})_*MU_{p}}(H^a(H; (MU_{p})_*(\mb{S}^{-i_H^* \diamond}\otimes X)))
\]
for the $(MU_{p})_*MU_{p}$-comodule
\[
	H^a(H; (MU_{p})_*(\mb{S}^{-i_H^* \diamond}\otimes X))
\]
and we obtain ${}_I\underline{E}_2^{a,b,*+\diamond}(G/H)$ by applying $d_1^\ANSS$.
In the $H=e$ case, the double complex spectral sequence collapses at ${}_I\underline{E}_2$ because it is concentrated in degree $a=0$.
\end{proof}

\begin{corollary}\label{cor: MANSS to ANSS}
The restriction to the trivial subgroup $e \le G$ gives a map of spectral sequences from the value of the MANSS at $G/H$ to the classical Adams--Novikov spectral sequence.  The induced map on $E_2$-terms is 
\[ 
	{}^{\MANSS}\underline{E}_2^{b,V}(X)(G/H)\longrightarrow \Ext_{(MU_{p})_*MU_{p}}^{b,b+|V|}((MU_{p})_*(X^e)).
\]
\end{corollary}

\begin{proposition}\label{prop:sseq-II}
There is an identification
\[ 
	{}_{\II}\underline{E}_2^{a,b,n+\diamond}(X)(G/H) \isom H^a(H; \Ext^{b,a+b+n}((MU_{p})_*(\mathbb{S}^{-\diamond}\tensor X))). 
\]
Moreover, when $X^e$ is $MU_{p}$-nilpotent complete and the $E_2$-term of the Adams--Novikov spectral sequence for $X^e$ is concentrated on the zero line, then this spectral sequence collapses and we can identify the $E_2$-term of the MANSS with 
\[
	\underline{E}_2^{a,n+\diamond}(G/H)=H^a(H; \pi_n (\Sigma^{-\diamond} X) ).
\] 
\end{proposition}

\begin{proof}
Starting with \eqref{eq:Z}, we have
\small
\begin{align*}
	{}_{\II}\underline{E}_1^{a,b,n+\diamond}(X)(G/H) & \isom \Ext^{b,a+b+n}(\underline{\pi}_*(F(EG^{(a)}/EG^{(a-1)}, \mathbb{S}^{-\diamond}\tensor X\tensor MU_p)(G/H)))\\
	& \isom \big(\Ext^{b,a+b+n}(\pi_*(F(EG^{(a)}/EG^{(a-1)},\mathbb{S}^{-\diamond}\tensor X\tensor MU_p)^{H}))\big) \\ 
	& \isom \big(\Ext^{b,a+b+n}(\pi_*((\mathbb{S}[G]\tensor \overline{\mathbb{S}[G]}^{\tensor a}\tensor \mathbb{S}^{-\diamond}\tensor X \tensor MU_p)))\big)^{H}\\
	& \isom \big(\Ext^{b,a+b+n}((MU_{p})_*((\mathbb{S}[G]\tensor \overline{\mathbb{S}[G]}^{\tensor a} )\tensor_{(MU_{p})_*} \pi_*(\mathbb{S}^{-\diamond}\tensor X \tensor MU_p)))\big)^{H}\\
	& \isom \big(\Ext^{b,a+b+n}(H\mathbb{Z}_*(\mathbb{S}[G]\tensor \overline{\mathbb{S}[G]}^{\tensor a} )\tensor (MU_{p})_*(\mathbb{S}^{-\diamond}\tensor X ))\big)^{H}\\
	& \isom \big(H\mathbb{Z}_*(\mathbb{S}[G]\tensor \overline{\mathbb{S}[G]}^{\tensor a} )\tensor \Ext^{b,a+b+n}((MU_{p})_*(\mathbb{S}^{-\diamond}\tensor X )))\big)^{H}
\end{align*}
\normalsize
where the third line uses Lemma \ref{lem:EG-quotient}, the fact that the homotopy fixed points spectral sequence computing the homotopy groups of the Borel spectrum collapses by $H$-freeness, and the fact that $H$-fixed points is an exact functor on the category of free $H$-modules. The last line is precisely the cellular chain complex computing the cohomology of $BH$ with local coefficients in the $\Z[H]$-module given by the $\Ext$ group. Thus we have
\begin{align*}
	{}_{\II}\underline{E}_2^{a,b,n+\diamond}(G/H) & \isom H^a(H; \Ext^{b,a+b+n}((MU_{p})_*(\mathbb{S}^{-\diamond}\tensor X))). \qedhere
\end{align*}
\end{proof}

\begin{corollary}
The MANSS is independent of a choice of Adams--Novikov resolution. The MANSS is also independent of the choice of $\mathbb{S}[G]$-Adams resolution.
\end{corollary}
\begin{proof}
Given an Adams--Novikov resolution $L_{\bullet}$ of $X^e$ and a map of Adams--Novikov resolutions $\overline{MU}_{p}^{\bullet}\otimes X^e \to L_{\bullet}$, Proposition \ref{prop:sseq-II} produces a map of double complex spectral sequences that is an isomorphism on the $E_2$-terms and satisfies the conditions of \cite[Theorem 8.3]{Boa99}. The second statement follows from the same argument using Proposition \ref{prop:sseq-I} and noting that tower \eqref{eq:hfpss-tower} may be regarded as a $\mathbb{S}[G]$--Adams resolution of $X$. 
\end{proof}

The reader will notice that the $E_2$-terms in Propositions~\ref{prop:sseq-I} and \ref{prop:sseq-II} compare to the $E_2$-term of the homotopy fixed point spectral sequence.  The following proposition gives a relationship between the $E$-MASS and the homotopy fixed point spectral sequence. 

\begin{proposition}\label{prop: MANSS to HFPSS}
For a $G$-spectrum $X$,
let ${}^\hfpss \underline{E}_r$ denote the $E_r$-term of the homotopy fixed points spectral sequence, which is
obtained by applying $\underline{\pi}_\star$ to the filtration \eqref{eq:hfpss-tower}.
There is a map of spectral sequences given on $E_r$-terms by 
\begin{align*}
	 {}^{\text{hfpss}}\underline{E}_r^{s,V}(X)\to {}^{E\mh \textup{MASS}}\underline{E}_r^{s,V}(X) 
\end{align*}
\end{proposition}

\begin{proof}
For each $a,b\in \mathbb{Z}$ satisfying $a+b=n\ge 0$, there is a canonical map 
\[
	F(EG_+\otimes \widetilde{E}G^{(a)},X)\otimes \overline{E}^{b}\longrightarrow 
 	\colim_{a+b=n}F(EG_+\otimes \widetilde{E}G^{(a)},X)\otimes \overline{E}^{\otimes b}.
\]
Fixing $b=0$, this gives a map of filtered spectra producing the map of spectral sequences.
\end{proof}

\subsection{The case of $MU_\mb{R}$-projective spectra}\label{sec: MUR proj case}

We now focus on the case of $G = C_2$ and $p = 2$, and let
$$ {}^{\MANSS}E_r(X) := {}^{\MANSS}\underline{E}_r(X)(C_2/C_2). $$
%By Corollary~\ref{thm:MANSS=E/a-Adams} and the shearing isomorphism, the MANSS agrees with the $MU_{\R,2}$-based modified Adams spectral sequence starting at the $E_2$ page. However, the double complex spectral sequences may be different, and in this section only we redefine $MU_2$ to be the Borel spectrum with complex conjugation action.

If $X \in \Sp_{C_2,i2}$ is $MU_\R$-projective (Definition~\ref{defn:MURprojective}), then  Theorem~\ref{thm:C2ESSSE2} gives a complete computation of the $E_1$-term of the $C_2$-ESSS.  However, this does not immediately carry over to a description of the $E_1$-term of the $a$-completed $C_2$-ESSS, which by Corollary~\ref{cor:ESSS is MANSS} is the $E_2$-term of the MANSS.  In this section, we will instead 
argue that the second double complex spectral sequence 
\[
	{}_{\II}E^{a,b,n-j+j\sigma}_2(X) = H^a(C_2; \Ext^{b,b+n}((MU_{2})_*(\mathbb{S}^{j-j\sigma}\otimes X))) \implies {}^{\MANSS}E_2^{a+b,n-j+j\sigma}(X) 
\]
for such $X$ collapses.  This method of computing the $E_2$-term of the MANSS has the added benefit of giving explicit meaning to the generators.

It is noteworthy that the property of being $MU_\R$-projective is not preserved by $a$-completion.  Thus we are in the unusual situation of computing the $E_2$-term of the MANSS, an inherently Borel equivariant construction, by assuming the Borel equivariant spectrum has a genuine equivariant lift which is $MU_\R$-projective.

\begin{remark}
	A similar analysis to that performed in this section can be used to prove the other double complex spectral sequence ${}_IE^{*,*,\star}_r(X)$ also collapses at $E_2$ for $MU_\R$-projective $X$, but the resulting presentation of the $E_2$-term of the MANSS will be less useful for us.
\end{remark}

We first identify the $E_2$-term 
\[ 
	{}_{\II}E^{a,b,n -j +j\sigma}_2(X) = H^a(C_2; \Ext^{b,b+n}((MU_{2})_*(\mathbb{S}^{j-j\sigma}\otimes X))) 
\]
of the second double complex spectral sequence in the $MU_\R$-projective case.  The first task is to identify the action of $C_2$.

\begin{lemma}\label{lem:action}
	Suppose $X \in \Sp_{C_2,i2}$ is $MU_\R$-projective. Then the action of the generator of $C_2$ on 
	$$\Ext^{s,2t}((MU_{2})_*(\mathbb{S}^{j-j\sigma}\otimes X^e))$$
	is multiplication by $(-1)^{t+j}$.
\end{lemma}

\begin{proof}
Note that as $MU_2$ is given the trivial $C_2$-action, the action of $C_2$ on 
$$ \Ext^{s,2t}_{MU_*MU_2}((MU_2)_*,(MU_{2})_*(X^e)) $$
is induced from the action of $C_2$ on the cosimplicial resolution
\begin{equation}\label{eq:trivres}
 X^e \to MU^{\otimes \bullet + 1}_2 \otimes X^e
\end{equation}
arising by the action of $C_2$ on $X^e$.  By contrast, letting $MU^e_{\RR,2}$ denote the underlying spectrum of $MU_{\R,2}$ with its non-trivial $C_2$-action by conjugation,
the cosimplicial resolution
\begin{equation}\label{eq:conjres}
 X^e \to (MU^e_{\R,2})^{\otimes \bullet+1} \otimes X^e
\end{equation}
has a $C_2$-action coming from the \emph{diagonal} action of $C_2$ on $(MU^e_{\R,2})^{s \bullet+1} \otimes X^e$.
We will write
\begin{equation}\label{eq:conjE2}
 \Ext^{s,2t}_{(MU^e_\R)_*MU^e_\R}((MU^e_{\R,2})_*, (MU^e_{\R,2})_*(X^e))
 \end{equation}
for the $E_2$-term of the associated Bousfield Kan spectral sequence as a $C_2$-module. 

However, since the different $C_2$ actions on the two Adams-Novikov resolutions (\ref{eq:trivres}) and (\ref{eq:conjres}) cover the \emph{same} $C_2$-action on $X^e$, the induced actions on $E_2$-pages are the same.  We deduce that despite the fact that $(MU^e_{\R,2})_*X^e$ may not be equivariantly isomorphic to $(MU_2)_*X^e$, it is nevertheless the case that there \emph{is} an equivariant isomorphism
\begin{equation}\label{eq:trivconjiso}
 \Ext^{s,2t}_{MU_*MU_2}((MU_2)_*,(MU_{2})_*(X^e)) \cong \Ext^{s,2t}_{(MU^e_\R)_*MU^e_\R}((MU^e_{\R,2})_*, (MU^e_{\R,2})_*(X^e)).
 \end{equation}
Therefore, it suffices to show that the action of the generator of $C_2$ on 
$$ \Ext^{s,2t}_{(MU^e_\R)_*MU^e_\R}((MU^e_{\R,2})_*, (MU^e_{\R,2})_*(\mathbb{S}^{j-j\sigma}\otimes X^e)). $$
is multiplication by $(-1)^{t+j}$.
    This follows directly from the $MU_\R$-projective condition. To see this, note that $C_2$ acts on the underlying spectrum $MU_{2}$ of $MU_{\R,2}$ by complex conjugation. Therefore, the generator of $C_2$ acts on $\pi_{2t}(MU^e_{\R,2} \otimes \mathbb{S}^{j-j\sigma})$ by $(-1)^{t+j}$ (see \cite[Rmk.~5.14]{HHR16}). 
\end{proof}

\begin{rmk}
Lemma~\ref{lem:action}, and in particular isomorphism (\ref{eq:trivconjiso}), may seem surprising at first for, e.g., $X = \mathbb{S}$, for it implies that
$$ \Ext^{s,2t}_{MU_*MU_2}((MU_2)_*,(MU_{2})_*) $$
is an $\FF_2$-vector space for $t$ odd.  However, this is actually true (see \cite[Cor.~2.6]{Burklund}).
\end{rmk}

\begin{notation}
\label{not:braxangle}
For an abelian group $A$, we let 
\[ 
	A[2] \coloneqq \{ x \in A \: : \: 2x = 0\}. 
\]
Let $X \in \Sp_{C_2,i2}$ be $MU_\R$-projective, and define
\[ 
	A^{s,2t} = \Ext^{s,2t}((MU_{2})_*X^e). 
\]
Then by Lemma \ref{lem:action} we have 
\begin{equation}\label{eq:anglebrax}
 	A^{s,2t} \cong H^0(C_2; \Ext^{s,2t}((MU_{2})_*(\mathbb{S}^{t-t\sigma}\otimes X^e))) = {}_{\II}E^{0,s,t\rho-s}_2(X).
\end{equation}
Given $x \in A^{s,2t}$, we denote its image under \eqref{eq:anglebrax} by $\bra{x}$.
By Lemma \ref{lem:action},
there is also an isomorphism
\begin{align}\label{eq:anglebraux}
	A^{s,2t}[2] & \cong   H^0(C_2; \Ext^{s,2t}((MU_{2})_*(\mathbb{S}^{(t+1)-(t+1)\sigma}\otimes X^e)))  \\ 
	\nonumber &=  {}_{\II}E^{0,s,t\rho-s+1-\sigma}_2(X).
\end{align}
Given $x \in A^{s,2t}[2]$, we denote its image under \eqref{eq:anglebraux} by $\bra{ux}$.
\end{notation}

Recall (see \cite{HuKriz}) that there is a subring (the ``positive cone'')
\[ 
	\Z[u^2,a]/(2a) \subset \pi^{C_2}_\star H\underline{\Z}. 
\]
with $|a| = -\sigma$ and $|u| = 1-\sigma$.  The images of the elements $a$ and $u^2$ in $\pi_\star H\underline{\Z}^{\wedge}_a$ induce isomorphisms 
\begin{align*} 
	\pi^{C_2}_{*+\diamond} H\underline{\Z}^{\wedge}_a  \cong H^{-*}(\mathbb{S}^{\diamond}_{hC_2})  \cong H^{-*}(C_2; H^0(\mathbb{S}^\diamond))  \cong \Z[u^{\pm 2}, a]/(2a).
\end{align*}
There is an isomorphism 
\[
	\Z_2[u^2, a]/(2a) \cong H^*(C_2;\Ext^{0,0}_{(MU_{2})_*MU_{2}}((MU_{2})_*, (MU_{2})_*(\mathbb{S}^{-\diamond}_2))) = {}_{\II}E^{*,0,\diamond}_2(\mathbb{S}_2). 
\]
The elements $u^2$ and $a$ detect the corresponding elements in
\begin{align*} 
	\pi^{C_2}_{\ast+\diamond} H\underline{\Z_2} \cong \Ext^{0,0}((MU_2)_* \otimes \pi_{\ast+\diamond}^{C_2} H\underline{\Z}) = {}^{C_2-\textup{ESSS}}E^{0,\ast+\diamond}_1(\mathbb{S}_2)_a^{\wedge} \cong {}^{\MANSS}E^{-\ast,\ast+\diamond}_{3}(\mathbb{S}_2)
\end{align*}
in the second double complex spectral sequence. Therefore these elements are permanent cycles in the double complex spectral sequence, and the second double complex spectral sequence is a spectral sequence of modules over $\Z[u^{\pm 2}, a]/(2a)$. The following proposition is an elementary immediate consequence of Lemma~\ref{lem:action} using the well-known group cohomology of $C_2$ with arbitrary coefficients.

\begin{proposition}\label{prop:doublecomplexE2}
Let $X \in \Sp_{C_2,i2}$ be $MU_\R$-projective. As before, define
\[ 
	A^{s,2t} = \Ext^{s,2t}((MU_{2})_*X^e). 
\]
Then under the isomorphisms \eqref{eq:anglebrax} and \eqref{eq:anglebraux} we have
\[ 
	{}_{\II}E_2^{*,*,\star}(X) \cong \big(A^{*,*} \otimes \Z[u^{\pm 2},a]/(2a)\big) \oplus \big(A^{*,*}[2] \otimes \mb{F}_2[u^{\pm 2},a]\big). 
\]
\end{proposition}

\begin{corollary}\label{cor:doublecomplexE2}
	For $MU_\R$-projective $X$, every element of ${}_{\II}E_2^{*,*,\star}(X)$ is of the form $u^{2i}a^j\bra{x}$ or $u^{2i}a^j \bra{ux}$ for $i \in \Z$, $j \ge 0$, and $x \in A^{*,*}$ (respectively $x \in A^{*,*}[2]$). 
\end{corollary}

Our strategy for proving that the second double complex spectral sequence collapses will be to produce elements in $\E{\MANSS}{}^{*,\star}_2(X)$ which the generators in Corollary~\ref{cor:doublecomplexE2} converge to.  To do this, we will require the following ``universal coefficient theorem''.

\begin{lemma}\label{lem:splitting}
For any torsion-free $(MU_{2})_*MU_{2}$-comodule $N$, there is a non-canonically split short exact sequence
\[ 
	\xymatrix@C-1em{
		0 \ar[r] & \Ext^{*,*}(N) \otimes \mb{F}_2 \ar[r]_-{\alpha} & \Ext^{*,*}(N/2) \ar@/_1pc/[l]_p \ar[r]_-{\beta} & \Tor_1^\Z (\Ext^{*+1,*}(N), \mb{F}_2) \ar[r] \ar@/_1pc/[l]_\iota & 0. 
	}
\]
If $N$ is an algebra in the category of comodules, then the short exact sequence is a short exact sequence of modules over $\Ext^{*,*}(N)$, and the splitting is a splitting in the category of $\Ext^{*,*}(N)$-modules.
\end{lemma}

\begin{proof}
The short exact sequence is easily deduced from the long exact sequence in $\Ext$ coming from the short exact sequence
\[ 
	0 \to N \xrightarrow{\cdot 2} N \rightarrow N/2 \to 0. 
\]
The splitting map $p$ comes from choosing a splitting $p'$ of the short exact sequence of free abelian groups
\[
	\xymatrix{
		0 \ar[r] & Z^{*,*} \ar[r] & C^{*,*}(N) \ar[r]_d \ar@{.>}@/_1pc/[l]_{p'} & B^{*+1,*} \ar[r] & 0
	}
\]
where $C^{*,*}(N)$ is the cobar complex for $N$, and $Z^{*,*}$ and $B^{*,*}$ are the subcomplexes of cocycles and coboundaries, respectively.
The composite
\[ 
	p''\colon \thinspace C^{*,*}(N) \to Z^{*,*} \twoheadrightarrow \Ext^{*,*}(N) 
\]
is then verified to be a map of cochain complexes, where the target is given zero differential.
The splitting map $p$ is then obtained by applying $H^{*,*}(-)$ to the map of cochain complexes 
\[ 
	C^{*,*}(N) \otimes \mb{F}_2 \xrightarrow{p'' \otimes 1} \Ext^{*,*}(N) \otimes \mb{F}_2. 
\]

Now assume that $N$ is a comodule algebra.  Then the maps $\alpha$ and $\beta$ are easily checked to be maps of modules over $\Ext^{*,*}(N)$.  Because the module structure on 
$\Ext^{*,*}(N) \otimes \mb{F}_2$ and 
$\Ext^{*,*}(N/2)$ come from the maps of rings
\[
	\Ext^{*,*}(N) \to \Ext^{*,*}(N) \otimes \mb{F}_2 \xrightarrow{\alpha} \Ext^{*,*}(N/2), 
\]
the fact that $p$ is a map of $\Ext^{*,*}(N)$-modules follows from the commutativity of the following diagram
\[
	\xymatrix{
		&& \Ext^{*,*}(N/2) \ar[dd]^{p} \\
		\Ext^{*,*}(N) \ar[r] & \Ext^{*,*}(N) \otimes \mb{F}_2 \ar[ur]^\alpha \ar[dr]_= & \\
		&& \Ext^{*,*}(N) \otimes \mb{F}_2
	}
\]
\end{proof}

\begin{notation}\label{not:xbar}
Suppose that $X \in \Sp_{C_2,i2}$ is $MU_\R$-projective, and choose an element $x \in \Ext^{s,2t}((MU_{2})_*X^e)$.  Let $\bar{x}$ denote the image of this element under the map
\begin{equation}\label{eq:brmap}	
  	\Ext^{s,2t}((MU_{2})_*X^e \otimes \pi_0^{C_2}H\underline{\Z}) \to {}^{C_2\mh \textup{ESSS}}E_1^{t,t\rho-s}(X) 
 \end{equation}
induced by Theorem~\ref{thm:C2ESSSE2}. By abuse of notation, we will also use $\bar{x}$ to denote the image of this element under the map 
\begin{equation}\label{eq:C2ESSStoMANSS}	
  {}^{C_2\mh \textup{ESSS}}E_1^{t,t\rho-s}(X) \to {}^{C_2\mh \textup{ESSS}}E_1^{t,t\rho-s}(X)^{\wedge}_a \cong {}^{\MANSS}E_2^{s,t\rho-s}(X). 
\end{equation}
Now choose a splitting $\iota$ of the short exact sequence
 \begin{equation}\label{eq:bruxSES}
 	\xymatrix@C-1em{
 		0 \ar[r] & A^{s-1,2t} \otimes \mb{F}_2 \ar[r] & \Ext^{s-1,2t}((MU_{2})_*(X^e)/2) \ar[r] & \Tor_1^{\Z}(A^{s,2t}, \mb{F}_2) \ar[r] \ar@/_1.5pc/[l]_\iota &  0
	 }
 \end{equation}
 as in Lemma~\ref{lem:splitting} (note $(MU_{2})_*X^e$ is torsion-free by the $MU_\R$-projective hypothesis). For
 \[ 
 	x \in A^{s,2t}[2] = \Tor_1^\Z(A^{s,2t},\mb{F}_2), 
\]
 let $\overline{ux}$ denote the image of $\iota(x)$ under the maps
\[ 
	\Ext^{s-1,2t}((MU_{2})_*X^e \otimes \mb{F}_2\{a\} ) = \Ext^{s-1,2t}((MU_{2})_*X^e \otimes \pi_{-\sigma}H\underline{\Z}) \to {}^{C_2 \mh ESSS}E_1^{t,t\rho-s+1-\sigma}(X) 
\]
	and 
\[ 
	{}^{C_2\mh \textup{ESSS}}E_1^{t,t\rho-s+1-\sigma}(X) \to {}^{C_2\mh \textup{ESSS}}E_1^{t,t\rho-s+1-\sigma}(X)^{\wedge}_a \cong {}^{\MANSS}E_2^{s,t\rho-s+1-\sigma}(X). 
\]
\end{notation}

\begin{remark}\label{rmk:brux}
Given $x \in \Ext^{s,2t}((MU_{2})_*X^e)$ of order 2, the element $\overline{ux}$ defined above is not unique, as it depends on a choice of lift.  However, the short exact sequence (\ref{eq:bruxSES}) implies that any two lifts $\overline{ux}$ and $\overline{ux}'$ differ by an element of the form $a\bar{y}$ for $y \in \Ext^{s-1, 2t}((MU_{2})_*X^e)$.
\end{remark}

\begin{remark}\label{rmk:brmult}
If $X$ is a ring spectrum, the maps (\ref{eq:brmap}) and  (\ref{eq:C2ESSStoMANSS}) are maps of rings, and therefore we have
\[ 
	\bar{x} \cdot \bar{y} = \overline{xy}. 
\]
Moreover the map $\iota$ in (\ref{eq:bruxSES}) is a map of modules over $A^{*,*}$ by Lemma~\ref{lem:splitting}.  It follows that we have 
\[ 
	\bar{x} \cdot \overline{uy} = \overline{uxy}. 
\]
\end{remark}

\begin{proposition}\label{prop:E2MANSS}
	Let $X\in \Sp_{C_2,i2}$ be $MU_\R$-projective, and $A^{*,*}$ as above.  For each $x \in A^{*,*}$, the element $\bar{x}$ is detected in the second double complex spectral sequence by the element $\bra{x}$, and for each $x \in A^{*,*}[2]$, the element $\overline{ux}$ is detected by $\bra{ux}$.  Therefore, the double complex spectral sequence collapses, and every element of ${}^{\MANSS}E^{*,\star}_2(X)$ is of the form $u^{2j}a^i\bar{x}$ or $u^{2j}a^i\overline{ux}$ for $j \in \Z$, $i \ge 0$, and $x$ as above. 
\end{proposition}

\begin{proof}
	Using the fact that there is a $C_2$-equivariant isomorphism $(MU_{2})_*((X \otimes Ca)^e) \cong (MU_{2})_*(X^e[C_2])$, the $E_2$-term of the second double complex spectral sequence for $X \otimes Ca$ is given by 
\[ 
	H^\ell(C_2;\Ext^{*,*}((MU_{2})_*(\mathbb{S}^{-\diamond}\otimes (X \otimes Ca)^e))) \cong \begin{cases}
																	A^{*,*}\otimes \Z[u^{\pm}], & \ell = 0, \\
																	0, & \ell > 0.
																  \end{cases}
\]
Here we use the powers of $u$ as a placeholder for the fact that the computation is independent of $\diamond$, but this also reflects the 
module structure over $\Z[u^{\pm 2}, a]/(2a)$ (where $a$ acts as zero).
Therefore, the second double complex spectral sequence for $X \otimes Ca$ collapses.
Consider
 the map of second double complex spectral sequences
\[
	\xymatrix{
		H^*(C_2,\Ext^{*,*}((MU_{2})_*(\mathbb{S}^{-\diamond}\otimes X^e))) \ar@{=>}[r] \ar[d] & {}^{\MANSS}E_2^{*,*+\diamond}(X) \ar[d] \\
		A^{*,*}\otimes \Z[u^{\pm}] \ar@{=>}[r] &  {}^{\MANSS}E_2^{*,*+\diamond}(X \otimes Ca)
	}	
\]
induced by the map
\[ 
	X \to X \otimes Ca. 
\]
Under this map, for each element $x \in A^{s,2t}$, the elements $\bra{x}$ and $\bar{x}$ map to $u^{-t}x$, and for each element $x \in A^{s,2t}[2]$, the elements $\bra{ux}$ and $\overline{ux}$ map to $u^{-t+1}x$.  Note that by Remark~\ref{rmk:brux}, this is independent of the choice of $\overline{ux}$, because $a$ acts trivially on the second double complex spectral sequence for $X \otimes Ca$. In particular, the map
\[ 
	{}_{\II}E_2^{0,*,\star}(X) \to {}_{\II}E^{0,*,\star}_2(X \otimes Ca) 
\]
is an injection. We may therefore deduce that the elements $\bra{x}$ detect $\bar{x}$ and the elements $\bra{ux}$ detect $\overline{ux}$ in the double complex spectral sequence for $X$.  It follows that the elements $\bra{x}$ and $\bra{ux}$ in ${}_{\II}E^{*,*,\star}_2(X)$ are all permanent cycles in the second double complex spectral sequence, and thus, by Corollary~\ref{cor:doublecomplexE2}, the double complex spectral sequence $\{{}_{\II}E^{*,*,\star}_r(X)\}$ collapses at the $E_2$-term.
\end{proof}

\subsection{The case of a trivial $C_2$-action}\label{sec: triv action case}

In this section we continue to assume that $G = C_2$ and $p = 2$.  Instead of assuming that $X$ is $MU_\R$-projective, we simply assume that $X$ is a Borel $C_{2}$-spectrum with $(MU_{2})_*X$ torsion-free, but we make the additional assumption that the action of $C_2$ on $X$ is \emph{trivial}.

The next lemma says that the $C_2/C_2$ level of the MANSS, which 
arises from applying $\pi^{C_2}_{*,*}$ to the filtered Borel $C_2$-spectrum
$\Gamma^h_{MU_{2}}$, can also be viewed as a collection of non-equivariant spectral
sequences, one for each weight $j$. In Section \ref{sec:MANSS-AHSS}, we will
define a $\mathbb{C}$-motivic analogue of the MANSS which generalizes this
$C_2/C_2$ level version.

\begin{lemma} \label{lem:h-to-Pinfty}
Suppose $X\in \Sp^{BC_2}_{i2}$ has trivial action.
Fix $j\in \Z$ and let
\begin{equation}\label{eq:GammahMUXj}
	\Gamma^{hC_2}_{MU_{2}}(X)[j]_{s} \coloneqq\begin{cases}
     									X^{P_j^\infty} & \text{ if } s  \le 0 \\ 
                  							\underset{a+b\ge  s}{\colim}\thinspace (\overline{MU}_2^{\tensor a} \tensor X)^{P_{j+b}^\infty} & \text{ if } s >0 .
    			 					  \end{cases}
\end{equation}
The spectral sequence associated to this filtered object agrees with the $C_2/C_2$-level of the MANSS starting in the $E_2$-term.
\end{lemma}
This agrees on the level of filtrations with the $C_2$-fixed points of a variant of the MANSS which uses a different filtration of $EC_2$, namely the filtration of $S^{\infty\sigma}\simeq EC_2$ coming from $S^{k\sigma}$ rather than the skeletal filtration of the bar complex model of $EC_2$; this is reflected in Lemma \ref{lem:h-to-Pinfty}. Nonetheless, as in Proposition \ref{prop:MANSS=NE-Adams}, the resulting spectral sequences agree starting on the $E_2$-term.
\begin{remark}
    The reader may wonder why the tensor with $\overline{MU}_{2}$ happens inside the parentheses in (\ref{eq:GammahMUXj}), when this tensoring happens outside the parentheses in (\ref{eq:GammahEX}).  However, the tensor in (\ref{eq:GammahEX}) is taken in the category of Borel $G$-spectra, and therefore homotopy fixed points are taken after taking the tensor.
\end{remark}

Next we identify the associated graded of \eqref{eq:GammahMUXj} to give a new description of the $C_2/C_2$ level of the MANSS $E_1$ term (see Definition \ref{def: MANSS}).
\begin{lemma} \label{lem:C2/C2-MANSS-E1}
There is an isomorphism
\[
	{}^{\MANSS}E_1^{s,2t-s+j\sigma}(X)
	\cong \bigoplus_{a+b = s} C^b_{\textup{cell}}(\Sigma^{-j}P^\infty_{j}; \mathbb{Z}) \otimes C^{a,2t+j}_{(MU_{2})_*MU_{2}}((MU_{2})_*X).
\]
where $C^{*,*}_{(MU_{2})_*MU_{2}}$ denotes the reduced cobar complex and $C^*_{\textup{cell}}$ denotes the cellular cochain complex.
\end{lemma}
Thus, the MANSS $E_1$-term is the total complex of the double complex 
\[
C^{\sbsq}_{\textup{cell}}(\Sigma^{-j}P_j^\infty ;\mathbb{Z} ) \otimes  C^{\bullet,2*+j}_{(MU_{2})_*MU_{2}}((MU_2)_*X)).
\]
Intuitively, the MANSS is a combination of the ANSS for $X$ and the Atiyah--Hirzebruch spectral sequence (AHSS) for $X^{P^{\infty}_j}$.  The goal of this subsection and the next is to make this intuition precise.  

We have already related the MANSS to the ANSS in Corollary~\ref{cor: MANSS to ANSS}.  In this subsection we will show that the first double complex spectral sequence always collapses in the case at hand, giving us a computation of the $E_2$-term of the MANSS which has an Atiyah--Hirzebruch flavor to it.  In the next subsection we will then add some additional hypotheses on $X$ which will allow us to relate higher differentials in the MANSS to differentials in the ANSS for $X$ and the algebraic AHSS for $\Ext^{*,*}((MU_{2})_*(X^{P_j^\infty}))$.

As in the previous subsection, for brevity we adopt the notation
\[ 
	A^{s,2t} \coloneqq \Ext^{s,2t}((MU_{2})_*X^e). 
\]
Because both the cobar complex and the cellular chain complex are torsion-free, we have a non-canonically split K\"unneth short exact sequence
\begin{equation}\label{eq:Kunneth} 
	0 \to \bigoplus_{a+b = s} A^{a,2t+j} \otimes H^b(\Sigma^{-j}P^\infty_{j}) \xrightarrow{\alpha} {}^{\MANSS}E_2^{s,2t-s+j\sigma}\\ 
	\xrightarrow{\beta} \bigoplus_{a+b = s+1} \Tor_1^\Z(A^{a,2t+j}, H^b(\Sigma^{-j}P_j^\infty)) \to 0.
\end{equation}

\begin{notation}\label{not:brax}
Let $1 \in H^0(P_0^\infty)$ denote the unit.
Given $x \in A^{s,2t}$, we let 
\[ 
	[x] \in {}^{\MANSS}E_2^{s,2t-s}(X) 
\]
denote the image of $x \otimes 1$ under the first map in \eqref{eq:Kunneth}.  This gives a map
\begin{equation}\label{eq:brax}
A^{s,2t} \hookrightarrow {}^{\MANSS}E_2^{s,2t-s}(X).
\end{equation}
Furthermore, there is an isomorphism 
\[ 
	\Tor_1^\Z(A^{s,2t}, H^1(\Sigma^{1}P_{-1}^\infty)) \cong A^{s,2t}[2]. 
\]
By choosing a splitting of the second map of \eqref{eq:Kunneth}, we get a map
\begin{equation}\label{eq:braux}
 	A^{s,2t}[2] \hookrightarrow {}^{\MANSS}E_2^{s,2t-s+1-\sigma}(X).
 \end{equation}
If $x \in A^{s,2t}[2]$, then we will denote the image of $x$ under these maps by $[ux]$.
\end{notation}

Recall from the previous subsection that ${}^{\MANSS}E_2^{*,\star}(X)$ is a module over $\pi^{C_2}_\star H \underline{\Z}^{\wedge}_a = \Z[a, u^{\pm 2}]/(2a)$. By \cite[Theorem~V.2.4]{Hinfty} we have
\[ 
	\pi^{C_2}_{i+j\sigma}(H\underline{\Z}^{\wedge}_a) \cong H^{-i}(P^\infty_j),
\]
and under this isomorphism the element $a^k u^{2l} \in \pi^{C_2}_{2l-(k+2l)\sigma}(H\underline{\Z}^{\wedge}_a)$ corresponds to a generator of $H^{-2l}(P_{-k-2l}^\infty)$,
and the first map of \eqref{eq:Kunneth} is a map of $\pi_\star^{C_2}(H\underline{\Z}^{\wedge}_a)$-modules.

\begin{remark}\label{rmk:braux}
Let $x$ be an element in $A^{s,2t}[2]$. Analogous to the definition of $\overline{ux}$ \eqref{not:xbar}, the choice of the element $[ux]$ may not be unique, but may depend on the choice of splitting of \eqref{eq:braux}. However different lifts $[ux]$ coming from two choices of splitting differ by an element of the form
\[ 
	\sum_i \left( \frac{a^2}{u^2}\right)^i a[y_i] 
\]
for $y_i \in A^{s-1-2i,2t}$. \\ 
\end{remark}

\begin{rmk}\label{rmk:bramult}
If $X$ is a ring spectrum, the map $\alpha$ in (\ref{eq:Kunneth}) is a map of rings, and therefore we have
\[ 
	[x] \cdot [y] = [{xy}]. 
\]
Moreover the map $\beta$ in (\ref{eq:Kunneth}) is a map of modules over $A^{*,*}$.  It follows from an argument analogous to that of Lemma~\ref{lem:splitting} that any splitting of (\ref{eq:Kunneth}) is a splitting of modules over $A^{*,*}$, and therefore given $x \in A^{s,2t}$ and $y \in A^{s',2t'}[2]$, we have 
\[ 
	[x] \cdot [uy] = [uxy]. 
\]
\end{rmk}

The following proposition is an immediate consequence of the K\"unneth theorem \eqref{eq:Kunneth}.

\begin{proposition}\label{prop:MANSSE2triv}
For $X \in \Sp^{BG}_{i2}$ with trivial action and $MU_*X$-torsion-free, under the maps \eqref{eq:brax} and \eqref{eq:braux} there is an isomorphism
\[ 
	{}^{\MANSS}E^{*,\star}_2(X) \cong \big(A^{*,*} \otimes \Z[u^{\pm 2},a]/(2a)\big) \oplus \big(A^{*,*}[2]\otimes \mb{F}_2[u^{\pm 2},a]\big). 
\]
\end{proposition}

\begin{corollary}\label{cor:MANSSE2triv}
For $X$ as in Proposition~\ref{prop:MANSSE2triv}, the first double complex spectral sequence
\[ 
	{}_{I}E_2^{a,b,2t-a-b-j+j\sigma}(X) = \Ext^{a,2t}((MU_{2})_*X \otimes H^{b+j}(P_j^\infty)) \implies {}^{\MANSS}E_2^{a+b, 2t-a-b-j + j\sigma}(X) 
\]
collapses.
\end{corollary}

\begin{proof}
Assume $P_j^\infty$ is given a finite type cellular structure.
Let $B^k \subseteq Z^k \subseteq C^k_{\mit{cell}}(P_j^\infty)$ denote the subgroups of coboundaries and cocycles of the cellular cochain complex.  Then we have a resolution of free abelian groups
\[ 
	0 \to B^{k} \to Z^k \to H^k(P_j^\infty) \to 0. 
\]
Since $(MU_{2})_*(X)$ is torsion free, we have a short exact sequence of $(MU_{2})_*(MU_{2})$-comodules
\[ 
	0 \to (MU_{2})_*X \otimes B^{k} \to (MU_{2})_*X \otimes Z^k \to (MU_{2})_* \otimes H^k(P_j^\infty) \to 0, 
\]
and therefore an exact sequence on Ext: 
\begin{equation*}
	A^{a,2t} \otimes B^{k} \xrightarrow{(1)} A^{a,2t} \otimes Z^k \to \Ext^{a,2t}((MU_{2})_*X \otimes H^{k}(P_j^\infty)) \to A^{a+1,2t}\otimes B^k \xrightarrow{(2)} A^{a+1,2t} \otimes Z^k. 
\end{equation*}
It follows that we have a short exact sequence (``universal coefficient theorem'')
\begin{equation}\label{eq:UCTExt}
	 0 \to A^{a,2t} \otimes H^k(P_j^\infty) \to \Ext^{a,2t}( (MU_{2})_*X \otimes H^{k}(P_j^\infty))  \to \Tor_1^\Z(A^{a+1,2t},H^k(P_j^\infty)) \to 0
 \end{equation}
involving the cokernel of (1) and the kernel of (2).
Comparing the description of ${}_I E_2$ in \eqref{eq:UCTExt} with the description of ${}^{\MANSS}E_2$ in \eqref{eq:Kunneth}, we see that there can be no differentials in the first double complex spectral sequence.
\end{proof}

\begin{rmk}
A similar argument, using the classical Universal Coefficient Theorem for cohomology, shows that the second double complex spectral sequence also collapses.
\end{rmk}

\subsection{The $\mathbb{C}$-motivic MANSS and the algebraic AHSS}\label{sec:MANSS-AHSS}

We continue with the assumption in the last subsection that $X \in \Sp^{BC_2}_{i2}$ is a Borel spectrum with trivial $C_2$ action, and $(MU_{2})_*X$ is torsion-free.
We will now proceed with our goal of relating differentials in the MANSS to those in the Atiyah--Hirzebruch spectral sequence (AHSS).
Our strategy will be to use complex motivic homotopy as a means of intrinsically encoding the Adams--Novikov aspect of the MANSS.  For this purpose we will add the following assumption for the rest of the section.

\begin{assumption}
The graded group $(MU_{2})_*X$ is concentrated in even degrees.
\end{assumption}

\noindent Let
\[ 
	\nu_\mb{C} : \Sp_{i2} \to \SH(\mb{C})_{i2}
 \]
be the section of Betti realization $\Be: \SH(\mb{C})_{i2} \to \Sp_{i2}$ constructed in \cite{GIKR}.  Letting $\tau \in \pi^{\mb{C}}_{0,-1}\mathbb{S}_2$ denote the canonical generator, we let $C\tau$ denote the cofiber of $\tau$, and define for $Z \in \SH(\CC)$
\[ 
	Z/\tau \coloneqq Z \otimes C\tau.
\]

We shall define a complex motivic lift $\E{\MANSS}{\mot}^{*,\star,*}_r(\nu_\CC X)$ of the MANSS (the ``motivic MANSS'').  There will then be a zig-zag
\[ 
	\E{\MANSS}{}_r^{*,\star}(X) \cong \E{\MANSS}{\mot}^{*,\star,*}_r(\nu_\CC X[\tau^{-1}]) \leftarrow  \E{\MANSS}{\mot}^{*,\star,*}_r(\nu_\CC X) \rightarrow  \E{\MANSS}{\mot}^{*,\star,*}_r(\nu_\CC X/\tau)
\]
and we will show that the spectral sequence $\E{\MANSS}{\mot}^{*,*+j\sigma,*}_r(\nu_\CC X/\tau)$ is a variant of the algebraic AHSS for $\Ext^{*,*}((MU_{2})_*(X^{P_j^\infty}))$.

Because of our hypothesis that $MU_*X$ is even, we have
\begin{equation}\label{eq:cofibertau}
 	\pi^{\mb{C}}_{n,t}(\nu_\mb{C}X/\tau) \cong \Ext^{2t-n,2t}((MU_{2})_*X) = A^{2t-n,2t}
\end{equation}
and the $\tau$-Bockstein spectral sequence for $\nu_{\mb{C}}X$ is equivalent to the ANSS for $X$.  Specifically, we have isomorphisms
\[ 
	(MU_{2})_{2t}(X) \cong (MGL_{2})_{2t,t}(\nu_\CC X) 
\]
and with respect to these isomorphisms we have an isomorphism
\[ 
	(MGL_{2})_{*,*}(\nu_\CC X) \cong (MU_{2})_*(X)[\tau]. 
\]
The motivic ANSS (the $MGL_2$-based ASS in $\SH(\CC)_{i2}$) for $\nu_\CC X$ takes the form
\[ 
	\E{\ANSS}{\mot}^{s,n,w}_2(\nu_\CC X)  \Longrightarrow \pi^{\CC}_{n-s,w}\nu_\CC X. 
\]
We have
\[ 
	\E{\ANSS}{\mot}^{s,2w-s,w}_2(\nu_\CC X) \cong A^{s,2w} 
\]
and we denote the image of an element $x \in A^{s,2w}$ under the above isomorphism by
\[ 
	\tilde{x} \in \E{\ANSS}{mot}^{s,2w-s,w}_2(\nu_\CC X). 
\]
In the case where $X = \mathbb{S}_2$ the element $\tau \in \pi^\CC_{0,-1}\mathbb{S}_2$ lifts to an element which we also denote 
\[ 
	\tau \in \E{\ANSS}{\mot}^{-2,0,-1}_2(\mathbb{S}_2). 
\]
This gives $\E{\ANSS}{\mot}^{*,*,*}_2(\nu_\CC X)$ the structure of a $\Z[\tau]$-module, and we have
\[ 
	\E{\ANSS}{\mot}^{*,*,*}_2(\nu_\CC X) \cong A^{*,*}[\tau]. 
\]
A differential 
\[ 
	d_{2r+1}(x) = y 
\]
in the ANSS for $X$ corresponds to a differential
\[ 
	d_{2r+1}(\tilde{x}) = \tau^r\tilde{y} 
\]
in the motivic ANSS for $\nu_\CC X$.  

Define 
\[ 
	\nu_{\CC}(X,P_j^\infty) \coloneqq \lim_{m}\nu_\CC(X^{P_j^{2m}}). 
\]
We take the limit outside the $\nu_\CC$ because $\nu_\CC$ does not commute with limits in general. 
We will be interested in the complex motivic homotopy groups
\[ 
	\pi^{\mb{C}}_{i,w}\nu_\mb{C}(X,{P_j^\infty}) 
\]
as lifts of the equivariant homotopy groups
\[ 
	\pi^{C_2}_{i+j\sigma}X^h \cong \pi_i(X^{P_j^\infty}). 
\]

\begin{remark}
It may be the case that there is an equivalence
\[
	\nu_\CC (X, P_j^\infty) \simeq (\nu_{\CC}X)^{L_j^\infty} 
\]
where the spectra $L_j^\infty$ are the stunted versions of $B\mu_2$ studied by \cite{GregersenRognes}, \cite{Quigley},
but we are unable to determine this.
\end{remark}

\begin{lemma}\label{lem:nuXPcofibers}
The maps $P_j^{2m} \to P^{2m}_{j+1}$ induce cofiber sequences
\begin{align*}
	\nu_\CC(X^{P_{j+1}^{2m}}) &\to \nu_{\CC}(X^{P_j^{2m}}) \to \Sigma^{-j, -\lceil j/2 \rceil} \nu_\CC X
	\\
	\nu_\CC(X,{P_{j+1}^{\infty}}) &\to \nu_{\CC}(X,P_j^\infty) \to \Sigma^{-j, -\lceil j/2 \rceil} \nu_\CC X. 
\end{align*}
\end{lemma}

\begin{proof}
Following, for example, \cite[Sec.~8]{Behrensroot}, using the torsion-free hypothesis on $(MU_{2})_*X$, we have isomorphisms of $(MU_{2})_*$-modules
\begin{equation}
	\begin{split}\label{eq:MU*P} (MU_{2})_*(X^{P_{2k-1}^{2m}}) & \cong (MU_{2})_*(X) \otimes_{(MU_{2})_*} \frac{x^k ((MU_{2})_*[x]/(x^{m}))}{\langle x^{i-1}[2](x) \: : \: i \ge {k} \rangle}, \\
	(MU_{2})_*(X^{P_{2k}^{2m}}) & \cong (MU_{2})_*(X) \otimes_{(MU_{2})_*} \frac{x^k ((MU_{2})_*[x]/(x^{m}))}{\langle x^{i-1}[2](x) \: : \: i \ge {k+1} \rangle}
	\end{split}
\end{equation}
where $\abs{x} = -2$ and $[2](x) \in (MU_{2})_*[[x]]$ is the 2-series of the universal formal group law.  
Applying $\nu_\CC$ to the cofiber sequences
\begin{gather*}
	X^{P^{2m}_{2j+1}} \to X^{P^{2m}_{2j}} \to \Sigma^{-2j}X  \\
	\Sigma^{-2j}X \to  X^{P^{2m}_{2j}} \to X^{P^{2m}_{2j-1}}, 
\end{gather*}
the result then follows from \cite[Lem.~3.13, Prop.~3.18]{GIKR}.
\end{proof}

\begin{definition} \label{def:motivic-MANSS}
Consider the following motivic lift $\Gamma^{h\mu_2}_{MGL_{2}}(\nu_{\CC}X)[j] \in \SH(\CC)^{\Fil}$ of the MANSS filtration:
\[ 
    \Gamma^{h\mu_2}_{MGL_{2}}(\nu_\CC X)[j]_{s} \coloneqq\begin{cases}
     \nu_\CC(X,{P_j^\infty}) & \text{ if } s  \le 0 \\ 
                   \underset{a+b\ge  s}{\colim}\thinspace  \overline{MGL}_2^{\tensor a} \widehat{\tensor}
                   \nu_\CC(X,{P_{j+b}^{\infty}})  &                                  \text{ if } s >0 .
     \end{cases}
\]
where the $\widehat{\otimes}$ notation is defined by  
\[ 
	Z \widehat{\otimes} \nu_{\CC}(X,P_\ell^\infty) \coloneqq \lim_m Z \otimes \nu_\CC(X^{P_\ell^{2m}})
\]
for $Z \in SH(\CC)_{i2}$.
We define the \emph{motivic MANSS} to be the spectral sequence associated to this filtration.
\end{definition}
Calculating the associated graded of this object gives a motivic analogue of Lemma \ref{lem:C2/C2-MANSS-E1}:
\[ 
    \E{\MANSS}{\mot}_1^{s,2t-s+j\sigma, *}(\nu_\CC X) \cong \dsums_{a+b=s} C_{\textup{cell},\textup{mot}}^{b+j,*}(\nu_\CC P_{j}^{\infty}) \tensor C^{a,2t+j,*}_{(MGL_2)_{*,*}MGL_2}((MGL_2)_{*,*}\nu_\CC X).
\]
where we are able to commute the limit in $\widehat{\otimes}$ past the colimit since the former is sequential and the latter can be taken to be finite.
Non-finiteness issues involved in dualizing $\nu_\CC(\mathbb{S}, P^\infty_j)$ to $\nu_\CC(P^\infty_j)$ are explained in \cite[Lemma 3.24]{GIKR}. Analogous arguments to \eqref{eq:Kunneth} hold in the motivic case and we obtain the following.

\begin{lemma}\label{lem:motivic-MANSS-E1}
There is a (non-canonically split) short exact sequence
\begin{multline}\label{eq:motKunneth} 
0 \to \bigoplus_{{a+b = s}} A^{a,2t} \otimes_{\mathbb{Z}_{2}} H^{b+j,w}(\nu_\CC P^\infty_{j}) \to \E{\MANSS}{\mot}_2^{s,2t-s-j+j\sigma,t-w}(\nu_\CC X) \\ 
\to \bigoplus_{a+b = s+1} \Tor_1^{\Z_{2}}(A^{a,2t}, H^{b+j,w}(\nu_\CC P_j^\infty)) \to 0.
\end{multline}
\end{lemma}

The same argument as Corollary~\ref{cor:MANSSE2triv} gives the following.

\begin{proposition}\label{prop:motMANSSE2}
There is an isomorphism
\begin{equation}\label{eq:motMANSSE2}
 \E{\MANSS}{\mot}_2^{s,2t-s-j+j\sigma, t-w}(\nu_\CC X) \cong \bigoplus_{a+b = s}\Ext^{a,2t}((MU_{2})_*X \otimes_{\mathbb{Z}_{2}} H^{b+j,w}(\nu_\CC P_j^\infty)).
\end{equation}
\end{proposition}

In particular, we have
\begin{equation}\label{eq:MANSS-mot-E2}
	\E{\MANSS}{\mot}_2^{*,*,*}(\nu_\CC X) \cong {}^{\MANSS}E^{*,*}_2(X)\otimes \Z[\tau]. 
\end{equation}
All elements of the $E_2$-term of the motivic MANSS may be represented by elements of the form
\[ \tau^m a^k u^{2l}[x] \text{ or }  \tau^m a^k u^{2l}[uy] \]
for $k,m \ge 0$, $l \in \Z$, $x \in A^{s,2t}$, $y \in A^{s,2t}[2]$, and 
\begin{align*}
[x] & \in \E{\MANSS}{\mot}_2^{s,2t-s,t}(\nu_\CC X),  
& [uy] & \in \E{\MANSS}{\mot}_2^{s,2t-s+1-\sigma,t}(\nu_\CC X),
\\
 \tau & \in \E{\MANSS}{\mot}_2^{-2,0,-1}(\nu_\CC X),
& a & \in \E{\MANSS}{\mot}_2^{1,-\sigma,0}(\nu_\CC X),
\\
u^2 & \in \E{\MANSS}{\mot}_2^{0,2-2\sigma,1}(\nu_\CC X). & & 
\end{align*}

Betti realization induces an equivalence \cite{DuggerIsaksen}
\begin{equation}
\SH(\CC)^{\mr{T}}_{i2}[\tau^{-1}] \simeq \Sp_{i2}
\end{equation}
where $\SH(\CC)^{\mr{T}}_{i2}$ denotes the subcategory of $\SH(\CC)_{i2}$ of cellular objects.  The following is an immediate consequence of \eqref{eq:MANSS-mot-E2}.

\begin{proposition}
    The $\tau$-local motivic MANSS with $E_{2}$-term $\E{\MANSS}{\mot}_2^{*,\star,*}(\nu_\CC X)[\tau^{-1}]$ is isomorphic to the $2$-primary MANSS with $E_{2}$-term ${}^{\MANSS}E_r^{*,\star}(X)[\tau ,\tau^{-1}]$.
\end{proposition}

\begin{proof}
The canonical map induces an isomorphism on $E_2$-terms. 
\end{proof}

We will now endeavor to identify the spectral sequence we obtain by smashing the motivic MANSS filtration with $C\tau$.
Define
\[ 
	\Ext^{*,*}((\widehat{MU}_{2})_*(X^{P_j^\infty})) \coloneqq \lim_m \Ext^{*,*}((MU_{2})_*(X^{P_j^{2m}})) \]
so we have
\[ 
	\pi^{\CC}_{n,t}\nu_{\CC}(X, P_j^\infty) = \Ext^{t-n, t}((\widehat{MU}_{2})_*(X^{P_j^\infty})). 
\]
By (\ref{eq:MU*P}), there are 
short exact sequences
\begin{align}
\label{eq:MU*PSESs1}
0 \to (MU_{2})_*(X^{P^{2m}_{2k+1}}) \to & (MU_{2})_*(X^{P^{2m}_{2k}}) \to (MU_{2})_*(\Sigma^{-2k}X) \to 0
\\
\label{eq:MU*PSESs2} 0 \to (MU_{2})_*(\Sigma^{-2k}X) \xrightarrow{\cdot x^{k-1}[2](x)} & (MU_{2})_*(X^{P^{2m}_{2k}}) \to  (MU_{2})_*(X^{P^{2m}_{2k-1}}) \to 0. 
\end{align}
Applying $\Ext(-)$ produces long exact sequences. We observe that each of the terms in these long exact sequences as well as each of the kernels and cokernels of the maps in the long exact sequences, are Mittag--Leffler sequences in $m$ and consequently applying the functor $\lim_m\Ext(-)$ also produces a long exact sequence. These long exact sequences splice together to give us an \emph{algebraic Atiyah--Hirzebruch spectral sequence (AHSS)} (compare with \cite[Sec.~8]{Behrensroot}):
\[
	\E{\AHSS}{\alg}_1^{l,s,t}(X,P_j^\infty) \implies \Ext^{s,t}( (\widehat{MU}_{2})_*(X^{P^{\infty}_j}))   
\]
with 
\[
\E{\AHSS}{\alg}_1^{l,s,t}(X,P_j^\infty) =
\begin{cases}
\Ext^{s,t}((MU_{2})_*(\Sigma^{-2k} X)), & l = 2k \ge j, \\
\Ext^{s+1,t}((MU_{2})_*(\Sigma^{-2k} X)), & l = 2k-1 \ge j, \\
0, & l < j.
\end{cases}
\]
The differentials in the algebraic AHSS are determined by the effects of the attaching maps of cells in $P_j^\infty$ on $\Ext((MU_{2})_*X)$.

There is a variant of the algebraic AHSS which involves taking 2 cells at a time.  Specifically, applying $\lim_m\Ext(-)$ to the decreasing filtration
\[ (MU_{2})_*(X^{P^{2m}_{j}}) \leftarrow \cdots \leftarrow (MU_{2})_*(X^{P^{2m}_{2k-1}}) \leftarrow (MU_{2})_*(X^{P^{2m}_{2k+1}}) \leftarrow \cdots \]
of $(MU_{2})_*(X^{P^{2m}_{j}})$ when $j$ is odd and the decreasing filtration 
\[ (MU_{2})_*(X^{P^{2m}_{j}}) \leftarrow  (MU_{2})_*(X^{P^{2m}_{j+1}})  \leftarrow  \cdots \leftarrow (MU_{2})_*(X^{P^{2m}_{2k-1}}) \leftarrow (MU_{2})_*(X^{P^{2m}_{2k+1}}) \leftarrow \cdots \]
when $j$ is even gives rise to a variant of the algebraic AHSS which we will called the \emph{accelerated algebraic AHSS}:
\[ \E{\AHSS}{\alg, \acc}_1^{k,s,t}(X,P_j^\infty) \implies \Ext^{s,t}((\widehat{MU}_{2})_*(X^{P^{\infty}_j}))   \]
with 
\[
\E{\AHSS}{\alg, \acc}_1^{k,s,t}(X,P_j^\infty) = 
\begin{cases}
\Ext^{s,t}(MU_*(\Sigma^{-2k} X \otimes M(2))), & k > j/2, \\
\Ext^{s,t}(MU_*(\Sigma^{-j} X)), & \text{$j$ even and} \: k = j/2, \\
0, & \text{otherwise}.
\end{cases}
\]
Note that the same considerations as before apply so that applying $\lim_m$ to the long exact sequence produces another long exact sequence. 

We define the motivic MANSS mod $\tau$ to be the spectral sequence with signature 
\[
\E{\MANSS}{\mot}_1^{s,i+j\sigma,w}(\nu_\CC X/\tau) \implies \pi^\CC_{i,w}(\nu_\CC(X,P_j^\infty)/\tau) 
\]
associated to the filtered 
motivic spectrum $\Gamma^h_{MGL_{2}}(\nu_\CC X)[j]/\tau$.  
Since we have
\[
H^{b+j,w}(\nu_\CC((P_j^\infty)_{2})/\tau) = 
\begin{cases}
\Z_{2}, & b = 0, \text{$j$ even}, w = \frac{j}{2},\\
\Z/2, & b > 0, \text{$b+j$ even}, w = \frac{b+j}{2}, \\
0 & \text{otherwise},
\end{cases}
\]
by (\ref{eq:motMANSSE2}) we have
\[  \E{\MANSS}{\mot}_2^{s,i+j\sigma, w}(\nu_\CC X /\tau) = 
\begin{cases}
\Ext^{2w-i, b}((MU_{2})_*X \otimes_{\mathbb{Z}_{2}} H^{b-2w, \frac{b}{2}-w}(\nu_\CC ((P_j^\infty)_{2})/\tau)) & b=i+s+j \text{ even,}  \\
0 & \text{otherwise}.
\end{cases}
\]

The following proposition gives an isomorphism between the motivic MANSS mod $\tau$ and the accelerated algebraic AHSS.  In particular, it implies that mod $\tau$, the differentials in the motivic MANSS are given by the effects of attaching maps in $P_j^\infty$ on $\Ext$.

\begin{theorem}\label{thm:AHSS}
The isomorphisms
\begin{multline*} \Ext^{2w-i,i+s+j}((MU_{2})_*X \otimes_{\mathbb{Z}_{2}} H^{i+s+j-2w, \frac{i+s+j}{2}-w}(\nu_\CC P_j^\infty/\tau)) \\
 \cong 
\begin{cases} 
 \Ext^{2w-i,2w}((MU_{2})_*(\Sigma^{-j}X)), 
 & i+s-2w = 0, \text{$j$ even}, \\
 \Ext^{2w-i,2w}((MU_{2})_*(X \otimes \Sigma^{2w-i-s-j}M(2))), 
  & i+s-2w > 0, \\
  0, & \text{otherwise}
\end{cases}
\end{multline*}
induce an isomorphism of spectral sequences
\[
\xymatrix{
\E{\MANSS}{\mot}^{s,i+j\sigma,w}_{2r+1}(\nu_\CC X/\tau)  \ar@{=>}[r] \ar[d]_\cong &
\pi^\CC_{i,w}(\nu_\CC(X, {P_j^\infty})/\tau) \ar[d]^\cong
\\
\E{\AHSS}{\acc,\alg}^{\frac{i+j+s}{2}-w, 2w-i, 2w}_{r}(X, P_j^\infty) 
\ar@{=>}[r] & \Ext^{2w-i,2w}((\widehat{MU}_{2})_*(X^{P_j^\infty})).
}
\]
\end{theorem}

\begin{proof}
The fundamental observation is that if $Y \in \Sp_{i2}$ and $(MU_{2})_*Y$ even, then 
\begin{equation}\label{eq: formula for f=0}
\pi^{\CC}_{i,w}(\overline{MGL}_2^{\otimes a} \otimes \nu_\CC(Y)/\tau)
=
\begin{cases}
\Ext^{2w-i,2w}((MU_{2})_*Y), & i < 2w-a, \\
Z^{2w-i,2w}((MU_{2})_*Y), & i = 2w-a, \\
0, & i > 2w-s
\end{cases}
\end{equation}
where $Z^{2w-i,2w}((MU_{2})_*Y) \subseteq C^{2w-i,2w}_{(MU_{2})_*{MU_{2}}}((MU_{2})_*Y)$ denotes the subgroup of cocycles in the normalized cobar complex.
This is proven by induction on $a$ using the cofiber sequences 
\[
\overline{MGL}_2^{\otimes a+1} \otimes \nu_\CC(X)/\tau \to 
\overline{MGL}_2^{\otimes a} \otimes \nu_\CC(X)/\tau
\to 
MGL_2 \otimes \overline{MGL}_2^{\otimes a} \otimes \nu_\CC(X)/\tau. 
\]

We will prove the theorem by directly comparing the exact couple of the accelerated Atiyah--Hirzebruch spectral sequence starting from the $E_2$-term to the exact couple of the motivic MANSS, run at double speed, starting from the $E_1$-term. To this end, we will first fix integers $j$ and $s$ and define a spectral sequence to compute 
\[ 
	\pi^\CC_{i,w} \Gamma^{h\mu_2}_{MGL_{2}}(\nu_\CC X)[j]_s/\tau. 
\]
We define a tower of spectra 
\[
\begin{tikzcd}
Y_s \ar[r]\ar[d,"\beta_s"']
&Y_{s-1} \ar[r]\ar[d,"\beta_{s-1}"] & \cdots \ar[r] & Y_1 \ar[r]& Y_0 \ar[r,equal]\ar[d,"\beta_0"]&\Gamma^{h\mu_2}_{MGL_{2}}(\nu_\CC X)[j]_s/\tau\\
S_s&S_{s-1}\ar[ul,"\alpha_{s-1}",dashed]&&&S_0\ar[ul,"\alpha_0",dashed]&
\end{tikzcd}
\]
where 
\[
Y_\ell \coloneqq \underset{\substack{a+b\ge  s \\ a \le \ell}}{\colim}\thinspace (\overline{MGL}_2^{\tensor a} \widehat{\tensor}
                   \nu_\CC(X,{P_{j+b}^{\infty}}))/\tau 
\]
and 
\[
S_f=\begin{cases}
\text{cof}(Y_{f+1}\to Y_f), &\text{ if }0\leq f\leq s-1,\\
	Y_s, & \text{ if }f=s, \text{ and }\\
	0 & \text{otherwise}.
\end{cases}
\]
Observe that we have 
\[
\pi_{i,w}^{\mathbb C}S_f =\begin{cases}
(MGL_2 \otimes \overline{MGL}_2^{\otimes f}) \widehat{\otimes} \nu_\CC(X,P_{j+s-f}^\infty)/\tau & \text{ if } 1\le f \le s \text{ and }\\
\overline{MGL}_2^{\otimes s} \widehat{\otimes} \nu_\CC(X,{P_j^{\infty}})  & \text{ if }f=0.\end{cases}
\]
so when $1\le f\le s$ we have
\begin{equation}\label{eq: formula for 1 le f le s}\pi_{i,w}^{\mathbb C}S_f = \begin{cases}
 \widehat{C}^{f,2w}_{(MU_{2})_*MU_{2}}(\widehat{(MU_{2})_*}(X^{P_{j+s-f}^\infty})), & i = 2w-f, \\
 0, & \text{otherwise}.
 \end{cases}
 \end{equation}
Applying $\pi_{i,w}^\mathbb{C}$, we obtain a conditionally convergent spectral sequence with signature
\begin{equation}
\label{eq:ssforY0}
E_1^{i,f,w} =\pi_{i,w}S_{f}\implies Y_0,
\end{equation}
with
\[
E_1^{i,f,w}=\begin{cases}
\widehat{C}^{f,2w}_{(MU_{2})_*{MU_{2}}}((\widehat{MU}_{2})_*(X^{P_{j+s-f}^\infty})), &i=2w-f, 0\leq f< s\\
Z^{s,2w}((\widehat{MU}_{2})_*(X^{P^{\infty}_j})), &i= 2w-s, f=s\\
\Ext^{2w-i,2w}((\widehat{MU}_{2})_*(X^{P^{\infty}_j})), &i< 2w-s, f=s\\
0, &\text{otherwise}
\end{cases}
\]
by \eqref{eq: formula for 1 le f le s} and \eqref{eq: formula for f=0}. Note that in each stem $i$, there is at most one non trivial bidegree (at filtration $f=2w-i$ when $2w-s\le i\leq 2w$, or $f=s$ for $i<2w-s$). Therefore, we will have no extensions and the spectral sequence is strongly convergent. For degree reasons, the only possibly non-trivial differentials are
\begin{equation}\label{eq: non edge}
d_1\colon \thinspace \widehat{C}^{f,2w}_{(MU_{2})_*MU_{2}}((\widehat{MU}_{2})_*(X^{P_{j+s-f}^\infty}))\to \widehat{C}^{f+1,2w}_{(MU_{2})_*MU_{2}}((\widehat{MU}_{2})_*(X^{P_{j+s-f-1}^\infty})),
\end{equation}
and 
\begin{equation}
\label{eq:edge}
d_1\colon \thinspace \widehat{C}^{s-1,2w}_{(MU_{2})_*MU_{2}}((\widehat{MU}_{2})_*(X^{P_{j+1}^\infty}))\to Z^{s,2w}((\widehat{MU}_{2})_*(X^{P^\infty_j})).
\end{equation}
Now, assume that $s+i+j$ is odd. Then $j+s+i-2w+1$ is even so the map
\[ 
	\widehat{C}^{*,*}_{(MU_{2})_*MU_{2}}((\widehat{MU}_{2})_*X^{P^\infty_{j+s+i-2w+1}}) \twoheadrightarrow \widehat{C}^{*,*}_{(MU_{2})_*MU_{2}}((\widehat{MU}_{2})_*X^{P^\infty_{j+s+i-2w}}) 
\]
is a surjection by \eqref{eq:MU*PSESs2}, and the map 
\[ 
	\widehat{C}^{*,*}_{(MU_{2})_*MU_{2}}((\widehat{MU}_{2})_*X^{P^\infty_{j+s+i-2w}}) \hookrightarrow \widehat{C}^{*,*}_{(MU_{2})_*MU_{2}}((\widehat{MU}_{2})_*X^{P^\infty_{j+s+i-2w-1}}) 
\]
is an injection by \eqref{eq:MU*PSESs1}.  
We compute \eqref{eq: non edge} using the commutative diagrams 
\[
\xymatrix{
\pi_{i,w}^\mathbb{C}S_{2w-i}=\widehat{C}^{2w-i,2w}_{(MU_{2})_*MU_{2}}((\widehat{MU}_{2})_*X^{P^\infty_{j+s+i-2w}}) \ar[r]^-{\alpha_{2w-i+1}} \ar[d]_d \ar[rd]^{d_1}&
\pi^\CC_{i-1,w}Y_{2w-i+1} \ar@{^{(}->}[d]^{\beta_{2w-i+1}}
\\
\widehat{C}^{2w-i+1,2w}_{(MU_{2})_*MU_{2}}((\widehat{MU}_{2})_*X^{P^\infty_{j+s+i-2w}}) \ar@{^{(}->}[r] &
\pi_{i-1,w}^\mathbb{C}S_{2w-i+1}=\widehat{C}^{2w-i+1,2w}_{(MU_{2})_*MU_{2}}((\widehat{MU}_{2})_*X^{P^\infty_{j+s+i-2w-1}})
}
\]
and
\[
\xymatrix{
\pi_{i,w}^\mathbb{C}S_{2w-i-1}=\widehat{C}^{2w-i-1,2w}_{(MU_{2})_*MU_{2}}((\widehat{MU}_{2})_*X^{P^\infty_{j+s+i-2w+1}}) \ar[r]^-{\alpha_{2w-i}} \ar@{->>}[d]\ar[dr]^{d_1} &
\pi^\CC_{i,w}Y_{2w-i} \ar@{^{(}->}[d]^{\beta_{2w-i}} 
\\
\widehat{C}^{2w-i-1,2w}_{(MU_{2})_*MU_{2}}((\widehat{MU}_{2})_*X^{P^\infty_{j+s+i-2w}}) \ar[r]_d &
\widehat{C}^{2w-i,2w}_{(MU_{2})_*MU_{2}}((\widehat{MU}_{2})_*X^{P^\infty_{j+s+i-2w}}). 
}
\]
We can compute the differential \eqref{eq:edge} in a similar way. 
We therefore have
\begin{equation}\label{eq:piGammamodtau}
\pi^\mathbb{C}_{i,w}\Gamma^{h\mu_2}_{MGL_{2}}(\nu_\CC X)[j]_s/\tau \cong 
\Ext^{2w-i,2w}((\widehat{MU}_{2})_*X^{P^\infty_{i+s-2w+j}})
\end{equation}
for $i+s+j$ odd. By construction, these isomorphisms are compatible in the sense that 
\[
    \begin{tikzcd}
    \pi^\mathbb{C}_{i,w}\Gamma^{h\mu_2}_{MGL_{2}}(\nu_\CC X)[j]_s/\tau \ar[r,"\cong"]  \ar[d] &  
\Ext^{2w-i,2w}((\widehat{MU}_{2})_*X^{P^\infty_{i+s-2w+j}}) \ar[d] \\
  \pi^\mathbb{C}_{i,w}\Gamma^{h\mu_2}_{MGL_{2}}(\nu_\CC X)[j]_{s-2/\tau} \ar[r,"\cong"] &  \Ext^{2w-i,2w}((\widehat{MU}_{2})_*X^{P^\infty_{i+s-2-2w+j}}) 
    \end{tikzcd}
\]
is a commutative diagram. 
Note that since 
\[ 
	\E{\MANSS}{\mot}^{s,i+j\sigma,w}_{2}(\nu_\CC X/\tau) 
\]
is concentrated in degrees with $i+j+s$ even, it only has non-trivial $d_r$-differentials for $r$ odd, and there are isomorphisms
\[ 
	\E{\MANSS}{\mot}^{s,i+j\sigma,w}_{2r}(\nu_\CC X/\tau) \cong 
\E{\MANSS}{\mot}^{s,i+j\sigma,w}_{2r+1}(\nu_\CC X/\tau).
\]
To prove the theorem, we must show that under the isomorphism
\[ 
	\E{\MANSS}{\mot}^{s,i+j\sigma,w}_{2}(\nu_\CC X/\tau)  \cong 
\E{\AHSS}{\acc,\alg}^{\frac{i+j+s}{2}-w, 2w-i, 2w}_{1}(X, P_j^\infty) 
\]
the $d_{2r+1}$ differentials in the motivic MANSS for $\nu_\CC X/\tau$ coincide with the $d_r$ differentials in the accelerated algebraic AHSS for $X^{P_j^\infty}$.
This follows from the following commutative diagram
\[
\begin{tikzcd}
\pi^\CC_{i,w}\left(\frac{\Gamma_s}{\Gamma_{s+1}}\right)\ar[d]\ar[r,hookleftarrow] & Z_1^{s,i+j\sigma,w} \ar[d,two heads]&   \\
\pi^\CC_{i-1,w}(\Gamma_{s+1}) & \E{\MANSS}{\mot}_2^{s,i+j\sigma,w}\ar[r,"\simeq"] &\E{\AHSS}{\acc,\alg}_1^{k, a-1,2w}\ar[d]\\
\pi^\CC_{i-1,w}(\Gamma_{s+2})\ar[u] \ar[rr,"\simeq"]&&	\Ext^{a,2w}((\widehat{MU}_{2})_*X^{P^\infty_{2k+1}})\\
\pi^\CC_{i-1,w}(\Gamma_{s+2r})\ar[u] \ar[rr,"\simeq"]&& \Ext^{a,2w}((\widehat{MU}_{2})_{}*X^{P^\infty_{2(k+r)-1}})\ar[u]\ar[d]\\
\pi^\CC_{i-1,w}(\Gamma_{s+2r+1}) \ar[u]\ar[d]\ar[rd]&
\E{\MANSS}{\mot}_2^{s+2r+1, i-1+j\sigma, w}\ar[r,"\simeq"]&  \E{\AHSS}{\acc,\alg}_1^{k+r, a,2w}\\
\pi^\CC_{i-1,w}\left(\frac{\Gamma_{s+2r+1}}{\Gamma_{s+2r+2}})\right) \ar[r,hookleftarrow]  &  
Z_1^{s+2r+1, i-1+j\sigma, w}\ar[u,two heads]
&\\
\end{tikzcd}
\]
 where $i+j+s$ is even, $k \coloneqq \frac{i+j+s-2w}{2}$, $a \coloneqq 2w-i+1$, we write
 \[\Gamma_s=\Gamma^{h\mu_2}_{MGL_{2}}(\nu_\CC X/\tau)[j]_s\]
and $Z^{*,*,*}_1 \subseteq \E{\MANSS}{\mot}^{*,*,*}_{1}(\nu_\CC X/\tau)$ denotes the subgroup of $d_1$-cycles.
\end{proof}

\section{Applications}\label{sec:applications}
We end with some computations of the MANSS to illustrate our perspective.

\subsection{Example: $ko_{C_2}$ at $p=2$}\label{sec: koC2-example}
Using the ($2$-primary) MANSS (see Definition \ref{def: MANSS}), we revisit the computation of 
the $RO(C_{2})$-graded homotopy groups of the $C_2$-homotopy fixed points of the $2$-complete $C_2$-equivariant connective real $K$-theory spectrum $\pi^{C_2}_{**}(ko_{C_2})^\wedge_{2,a}$. By \cite[\S 10]{BHS}, this is the same information as the Artin--Tate homotopy groups of the $a$-complete very effective Hermitian $K$-theory $\pi^{\textup{AT}}_{***}(kq^\wedge_{2,a})$. 
The MANSS for $(\mathrm{ko}_{C_{2}})_{2,a}^{\wedge}$ has signature
\[ 
	\E{\MANSS}{}_2^{s,i+j\sigma}((ko_{C_2})^{\wedge}_{2,a}) \implies \pi^{C_2}_{i+j\sigma}((ko_{C_2})^\wedge_{2,a}).
\]
Recall that the underlying spectrum of $(ko_{C_2})^{\wedge}_{2,a}$ is the $2$-complete connective real $K$-theory spectrum $ko_2$, whose Adams--Novikov $E_2$-term is of the form 
\[ {}^{\ANSS}E_2^{*,*}(ko_2)=\mathbb Z_2[h_1, v_1^2]/2h_1\]
where the generators have degree $|h_1|={(1,2)}$ and $|v_1^2|={(0,4)}$. Here we use Adams grading convention. 

There are two avenues to computing the $E_2$-term of the MANSS, namely those
given in Sections~\ref{sec: MUR proj case} and \ref{sec: triv action case}.
\begin{lemma}
The results in Section \ref{sec: MUR proj case} apply to the spectrum $(ko_{C_2})_2^\wedge$.
\end{lemma}
\begin{proof}
The $MU_\R$-projective hypothesis was used in that section in three ways.
First, we use the fact that $(MU_2)_*X^e$ is torsion-free in \eqref{eq:bruxSES};
this holds for $X^e = ko^\wedge_2$.
Second, in Lemma~\ref{lem:action}, we use this hypothesis to
calculate the action of $C_2$ on $\pi_{2f}(MU_{\R,2}^e\tensor
\mathbb{S}^{j-j\sigma})$.
Consider the long exact sequence in $MU_{\R,2}$-homology associated to the
``Wood cofiber sequence''
\[ \Sigma^{\sigma} (ko_{C_2})_2^\wedge \xrightarrow{{\eta}^\wedge_2} (ko_{C_2})_2^\wedge \to k\R_2 \]
of \cite[Prop.~10.13]{GHIR}.
Since the Hurewicz image of $\br{\eta} \in \pi^{C_2}_\sigma (\mb{S})$ in $\pi_{\sigma}MU_\R$ is trivial,
the long exact sequence splits, and so $(MU_{\R,2})_{*,*}(ko_{C_2})^\wedge_2$
injects into $(MU_{\R,2})_{*,*}k\R_2$.
So it suffices to calculate the $C_2$-action on $(MU_{\R,2})_{*,*}k\R_2$.
Since $k\R_2$ is strongly even (for example, using the criterion in \cite[Lemma 1.2(ii)]{Gre18} and the calculation of $\pi_{*,*}^{C_2}k\R_2$ in \cite{Behrens-Shah}), there is a splitting $MU_\R\tensor k\R_2\hteq \bigvee_k \Sigma^{k\rho} k\R_2$ by \cite[Proposition 2.31]{AKKQ-syntomic}.
Since the action on $k\R$ comes from complex conjugation, the generator of $C_2$ acts on $\pi_*(k\R_2^e \tensor \mathbb{S}^\star)$  with the same formula as in Lemma~\ref{lem:action} for $MU_{\R,2}$.

%Third, the $MU_\R$-projective property is a hypothesis for Theorem \ref{thm:C2ESSSE2}. However, the proof of that
%theorem applies in greater generality than the statement:
%it suffices to show that \cite[Lemma 5.39]{BHS} holds for $X$,\mhana{I am not quite sure what it means for BHS Lemma 5.39 to hold for $X$? Do you mean replacing $Sp^{proj}_{C_2,i2}$ in BHS 5.39 with the full subcategory of $Sp_{C_2,i2}$ consisting of spectra that is strongly even after tensoring $MU_{\mathbb{R},2}$? The blue arguement looks good to me.} which is
%equivalent to the condition that $MU_{\R,2}\tensor X$ is strongly even (\cite[Definition 3.1]{HM17}).
%To show this holds for $(ko_{C_2})_2^\wedge$, by the above inclusion of $(MU_{\R,2})_{*,*}(ko_{C_2})^\wedge_2$
%into $(MU_{\R,2})_{*,*}k\R_2$ it suffices to show that $MU_{\R,2}\tensor
%k\R_2$ is strongly even, and that follows by the splitting above and the fact
%that $k\R_2$ is strongly even.

Third, the $MU_\R$-projective property is a hypothesis for Theorem
\ref{thm:C2ESSSE2}. However, we may check the result directly using knowledge of
the $C_2$-effective slice spectral sequence for $(ko_{C_2})^\wedge_2$
\cite{Kon20} and the fact
(\cite[Proposition 1.37]{BHS}) that the motivic effective slice spectral
sequence agrees with the $\ta$-Bockstein spectral sequence.
\end{proof}

Now it follows from Proposition~\ref{prop:E2MANSS} and Remark~\ref{rmk:brmult} that we have
\[ 
\E{\MANSS}{}_2^{*,\star}((ko_{C_2})^{\wedge}_{2,a}) \cong 
\Z_2[u^{\pm 2},a,\bar{h}_1, \overline{v^2_1}]/(2a, 2\bar{h}_1) \oplus \mb{F}_2[u^{\pm 2}, a, \bar{h}_1, \overline{{v}^2_1}]\{\overline{uh_1}\}.
\]
Here, labelling degrees as {(filtration, $RO(C_2)$-stem)}, we have
\begin{align*}
 |u^2| & = (0, 2-2\sigma), & 
 |a| & = (1, -\sigma), \\
 |\bar{h}_1| & = (1, \sigma), & 
 |\overline{v^2_1}| & = (0, 2\rho), \\
 |\overline{uh_1}| & = (1,1). 
 \end{align*}
Note that while in general the generator $\overline{uh_1}$ depends on a choice of a lift, by Remark~\ref{rmk:brux}, this lift is unique in this case because $\Ext^{0,2}(MU_*ko) = 0$.  

However, as the underlying spectrum $ko_2^\wedge$ of $(ko_{C_2})^{\wedge}_{2,a}$ has a trivial $C_2$-action, Proposition~\ref{prop:MANSSE2triv} gives a \emph{different} presentation of the MANSS $E_2$-term:
\[ 
\E{\MANSS}{}_2^{*,\star}((ko_{C_2})^{\wedge}_{2,a}) \cong 
\Z_2[u^{\pm 2},a,[h_1], [v^2_1]]/(2a, 2[h_1]) \oplus \mb{F}_2[u^{\pm 2}, a, [h_1], [v^2_1]]\{[uh_1]\}.
\]
with 
\begin{align*}
 |[h_1]| & = (1,1), & 
 |[v^2_1]| & = (0,4), \\
 |[uh_1]| & = (1,2-\sigma). 
 \end{align*}
Note that analogous to the case of $\overline{uh_1}$, the generator $[uh_1]$ has the potential to depend on a choice of a lift, but by Remark~\ref{rmk:braux}, because $\Ext^{0,2}((MU_{2})_*ko_{2}) = 0$, this lift is also unique in this case.

Finally, Proposition~\ref{prop:motMANSSE2} and Remark~\ref{rmk:bramult} gives the motivic MANSS $E_2$-term:
\[ 
	\E{\MANSS}{\mot}_2^{*,\star,*}((\nu_\CC ((ko)_2)) \cong 
\Z_2[u^{\pm 2},a,[h_1], [v^2_1], \tau]/(2a, 2[h_1]) \oplus \mb{F}_2[u^{\pm 2}, a, [h_1], [v^2_1], \tau]\{[uh_1]\}.
\]
Recording tridegrees as (filtration, $RO(C_2)$-stem,  motivic weight), we have
\begin{align*}
 |u^2| & = (0, 2-2\sigma, 1), & 
 |a| & = (1, -\sigma,0), \\
|[h_1]| & = (1,1,1), & 
|[v^2_1]| & = (0, 4,2), \\
|[uh_1]| & = (1, 2-\sigma,1), & 
|\tau| & = (-2, 0,-1).
\end{align*}
Theorem~\ref{thm:AHSS} implies that mod $\tau$ the motivic AHSS can be identified with the accelerated algebraic AHSS for $\Ext((\widehat{MU}_{2})_*ko_{2}^{P_j^\infty})$.  Under this correspondence, the class
$
	a^k u^{2l}[u^{\epsilon} x]
$
corresponds to the class 
\[ 
	x[-2l-\epsilon] \in \E{\AHSS}{\alg}^{*,*,*}_1(ko_2,P^\infty_{-k-2l-\epsilon}) 
\]
which represents ``$x$ on the $(-2l-\epsilon$)-cell of $P^\infty_{-k-2l-\epsilon}$''.

The following theorem explains the philosophy that in the case of a trivial action, MANSS differentials can be derived from ANSS differentials and algebraic AHSS differentials.

\begin{theorem}
	\label{thm:kodiff}
	In the MANSS for $(ko_{C_2})^{\wedge}_{2,a}$, the $d_3$-differentials are determined by  
	\begin{enumerate}
 \item \label{it: a is permanent cycle} $d_3(a) = 0$, 
\item \label{it d3 on v1squared} $d_3([v_1^2])= [h_1^3],$
   \item  \label{it d3 on usquared} $d_3(u^2) = a^2[h_1]$,
   	\item  \label{it d3 uh1} $d_3([uh_1])= \frac{a^2}{u^2}[uh_1^2]$.
	\end{enumerate}
For dimensional reasons we have $E_4 = E_\infty$.
\end{theorem}

\begin{proof}
To prove \eqref{it: a is permanent cycle}, note that the element $a$ is always a permanent cycle in the MANSS due to the fact that it comes from $\pi^{C_2}_{-\sigma}\mathbb{S}^0$.

	 To prove \eqref{it d3 on v1squared}, note that for dimensional reasons, we have
 \[ 
 	d_3([v_1^2]) = \alpha [h_1^3] + \beta \frac{a^2}{u^2} [v_1^2h_1] 
\]
 for $\alpha, \beta \in \mb{F}_2$.
 Consider the map from the MANSS to the ANSS from Corollary \ref{cor: MANSS to ANSS}. 
  In the Adams--Novikov spectral sequence, we have $d_3(v_1^2)=h_1^3$. Since $[v_1^2]$ maps to $v_1^2$ and $\alpha [h_1^3]+\beta \frac{a^2}{u^2}[v_1^2h_1]$ maps to $\alpha h_1^2$, we deduce that $\alpha = 1$. 
  In the motivic MANSS for $\nu_\CC (ko_2)$ we therefore have
\begin{align*}
 d_3([v_1^2]) & = \alpha \tau [h_1^3] + \beta \frac{a^2}{u^2} [v_1^2h_1] \\
  & = \tau [h_1^3] + \beta \frac{a^2}{u^2} [v_1^2h_1].
  \end{align*}
Therefore in the motivic MANSS for $\nu_\CC(ko_2)/\tau$ we have 
\[ 
d_3([v_1^2]) = \beta \frac{a^2}{u^2} [v_1^2h_1].
\]
Since the $0$-cell splits off of $P_0^\infty$, we have
\[ 
	d_1([v_1^2]) = 0 
\]
in the accelerated algebraic AHSS for $\Ext((\widehat{MU}_{2})_*ko_2^{P_0^\infty})$.  By Theorem~\ref{thm:AHSS} we deduce $\beta = 0$.

To prove \eqref{it d3 on usquared}, note that for dimensional reasons, in the MANSS we have 
 \[ 
 	d_3(u^2) = \alpha a^2[h_1]
\]
for $\alpha \in \mb{F}_2$.  This corresponds to a differential
\[ 
	d_3(u^2) = \alpha a^2[h_1]
\]
in the motivic MANSS for $\nu_\CC(ko_2)$, and hence a differential
\[ 
	d_3(u^2) = \alpha a^2[h_1]
\]
in the motivic MANSS for $\nu_\CC(ko_2)/\tau$.  However, since the 0-cell attaches to the -2-cell of $P^\infty_{-2}$ by an $\eta$-attaching map, we have
\[ d_1(u^2) = a^2[h_1] \]
in the accelerated algebraic AHSS for $\Ext((\widehat{MU}_{2})_*P^\infty_{-2})$.  It follows that $\alpha = 1$. (Note \eqref{it d3 on usquared} could also have been deduced from Proposition~\ref{prop: MANSS to HFPSS}.)

	To prove \eqref{it d3 uh1}, note that for dimensional reasons, in the MANSS we have 
 \[ d_3([uh_1]) = \alpha \frac{a^2}{u^2}[uh_1^2] \]
for $\alpha \in \mb{F}_2$.  This corresponds to a differential
 \[ d_3([uh_1]) = \alpha \frac{a^2}{u^2}[uh_1^2] \]
in the motivic MANSS for $\nu_\CC(ko_2)$, and hence a differential
 \[ d_3([uh_1]) = \alpha \frac{a^2}{u^2}[uh_1^2] \]
in the motivic MANSS for $\nu_\CC(ko_2)/\tau$.  However, since the 1-cell attaches to the -1-cell of $P^\infty_{-1}$ by an $\eta$-attaching map, we have
\[ d_1([uh_1]) = \frac{a^2}{u^2}[uh_1^2] \]
in the accelerated algebraic AHSS for $\Ext((\widehat{MU}_{2})_*P^\infty_{-1})$.  It follows that $\alpha = 1$.
\end{proof}

The element $[uh_1]$ is not a product in the $E_2$-term of the MANSS for $(ko_{C_2})^{\wedge}_{2,a}$ because $u$ does not exist on the $E_2$-term.  The following proposition is therefore needed to completely identify the multiplicative structure in $\E{\MANSS}{}((ko_{C_2})^{\wedge}_{2,a})$.

\begin{proposition}
	\label{prop:exoticmulti}
 In $\E{\MANSS}{}_2((ko_{C_2})^{\wedge}_{2,a})$ there is a multiplicative relation
	\[ 
		[uh_1]^2 = u^2 [h_1^2] + a^2 [v_1^2].
	\]
\end{proposition}

\begin{proof}
For degree reasons, we know that
	\[ [uh_1]^2=\alpha u^2 [h_1^2] + \beta a^2 [v_1^2] \]
for some $\alpha,\beta\in \mb{F}_2.$ 
	By the Leibniz rule, we have 
 \begin{align*}
	0= d_3([uh_1]^2) & = d_3(\alpha u^2 [h_1^2] + \beta a^2 [v_1^2]) \\
 & = \alpha a^2 [h^3_1] + \beta a^2 [h_1^3].
 \end{align*}
	Therefore, we have $\alpha =\beta$. Note that under the map of spectral sequences from the MANSS to the ANSS given by Corollary~\ref{cor: MANSS to ANSS}, we have $[u h_1]$ maps to $h_1$. Since $h_1^2 \neq 0$, we conclude that $\alpha=\beta =1.$
\end{proof}

Next we discuss how to compare the classes $[x]$ and $\bar{x}$ that appear in the two forms of 
$\E{\MANSS}{}_2^{*,\star}(X^e)$ when the results of both Sections \ref{sec: MUR proj case} and \ref{sec: triv action case} hold for $X$.
If $x \in \Ext^{s,2t}((MU_{2})_*X^e)$, then we have two different MANSS classes:
\begin{align*}
\bar{x} & \in \E{\MANSS}{}^{s,t-s+t\sigma}(X^\wedge_a), \: \text{and} \\
[x] & \in \E{\MANSS}{}_2^{s,2t-s}(X^{\wedge}_a). 
\end{align*}
By Lemma~\ref{lem:action}, we either have {$t = 2k$ or $2x = 0$}.  Write $t = 2k+\epsilon$.
By considering the map to the ANSS (Corollary \ref{cor: MANSS to ANSS}), we have 
\[ [x] \equiv u^{2k}\overline{u^\epsilon x} \pmod a. \]
For example, when $X = (ko_{C_2})^\wedge_{2,a}$, we deduce that 
\begin{align*}
[v^2_1] & = u^2\overline{v_1^2}, \\
[h_1] & = \overline{uh_1}, \text{ and }\\
[uh_1] & = u^2\bar{h}_1
\end{align*}
for dimensional reasons. 
In particular, squaring the last of these equations and using Proposition~\ref{prop:exoticmulti} gives
$$ u^4 \bar{h}_1^2 = u^2[h_1^2] + a^2[v_1^2]. $$
Dividing by $u^4$, we get
\begin{equation}\label{eq:h1^2bar}
\overline{h_1^2} = u^{-2}[h_1^2] + u^{-4}a^2[v_1^2].
\end{equation}
which gives an example of an $x \in \Ext$ where $[x] \ne u^t\bar{x}$.
More generally we have
\[ \overline{h_1^{2k}} = (u^{-2}[h_1^2] + u^{-4}a^2[v_1^2])^k. \]

Analogous to the case of $[uh_1]$, the element $\overline{uh_1}$ is not a product in the $E_2$-term of the MANSS for $(ko_{C_2})^{\wedge}_{2,a}$. We deduce the following exotic multiplicative relation, which was originally proven by the fourth author in \cite{Kon20} in the $C_2$-ESSS.

\begin{proposition}\label{prop: exotic mult rel}
    In $\E{\MANSS}{}((ko_{C_2})^{\wedge}_{2,a})$ there is a multiplicative relation
	$$\overline{uh_1}^2 = u^2 \overline{h_1^2} + a^2 \overline{v_1^2}.$$
\end{proposition}

\begin{proof}
Solving (\ref{eq:h1^2bar}) for $[h_1^2]$, we have
\begin{align*}
[h_1]^2 & = u^2 \overline{h_1^2} + u^{-2} a^2 [v_1^2] \\
& = u^2 \overline{h_1^2} + a^2 \overline{v_1^2}.
\end{align*}
The result follows from the fact that $[h_1] = \overline{uh_1}$.
\end{proof}

We deduce the following differentials in the $C_2$-ESSS naming convention, which correspond to the generating differentials in the $C_2$-ESSS for $ko_{C_2}$ as computed by the fourth author in \cite{Kon20}. These differentials are implied by the differentials computed in Theorem \ref{thm:kodiff} and the multiplicative relation from Proposition \ref{prop: exotic mult rel} so we only mention them in order to provide a dictionary between the computations above and the computations in \cite{Kon20}. 

\begin{corollary}
In the MANSS for $(ko_{C_2})^{\wedge}_{2,a}$, we have
\begin{align*}
d_3(u^2)  = a^2 \overline{uh_1}, \quad \text{and } \quad
d_3(\bar{v}_1^2)  = \overline{uh_1^3}.
\end{align*}
\end{corollary}

\begin{proof}
The first differential follows immediately from Theorem~\ref{thm:kodiff}\thinspace \eqref{it: a is permanent cycle}.  The second is computed from Theorem~\ref{thm:kodiff} and Proposition~\ref{prop:exoticmulti}.  On the one hand, we have
\[
	d_3(u^2\overline{v_1^2})= d_3([v_1^2]) = (\overline{uh_1})^3 \stackrel{\text{Thm.} \ref{thm:kodiff}}{=}[h_1]^3=\overline{uh_1} \cdot (u^2 \bar{h}_1^2 + a^2 \overline{v_1^2}).
\]
However, the Leibniz rule gives
\[
	d_3(u^2\overline{v_1^2}) = a^2\overline{uh_1} \cdot \overline{v_1^2} + u^2 d_3(\overline{v_1^2}). 
\]
We deduce that
\[ 
	d_3(\overline{v_1^2}) = \overline{uh_1} \cdot \bar{h}_1^2 = \overline{uh_1^3}. 
\]
\end{proof}

\subsection{Example: $\tmf(2)$ at $p=3$}
We now consider the 3-complete spectrum $\tmf(2)_{3}$ of topological modular forms for $\Gamma(2)$ studied in depth in \cite{HahnSengerWilson}.  Recall that a $\Gamma(2)$-structure on an elliptic curve $C$ is a choice of isomorphism
\[ 
	\FF_2 \oplus \FF_2 \xrightarrow{\cong} C[2] 
\]
where $C[2]$ denotes the $2$-torsion subgroup-scheme of $C$.  Associated to the Deligne-Mumford compactification of the moduli stack of elliptic curves with level 2 structure is a non-connective commutative ring spectrum $\Tmf(2)$ \cite[\S 5.3]{Stojanoska}.
The group $GL_2(\FF_2)\cong \Sigma_3$ acts transitively on the set of $\Gamma(2)$-structures, and this gives $\Tmf(2)$ the structure of a commutative Borel $\Sigma_3$-spectrum, which we may regard as a genuine commutative $\Sigma_3$-spectrum.
Following \cite[Notation 1.7]{HahnSengerWilson}, we define $\tmf(2) \in \Sp_{\Sigma_3}$ to be the $\Sigma_3$-equivariant connective cover of $\Tmf(2)$. 
The periodic version of $\tmf(2)_3$ is given by
\[ 
	\TMF(2)_3 \coloneqq \tmf(2)_3[\Delta^{-1}] \simeq (\tmf(2)_{E(2)})_3. 
\]
By \cite[Theorem 6.1]{Stojanoska}, we have
\[ 
	\TMF_{3} \simeq \TMF(2)_{3}^{h\Sigma_3}. 
\]  
Recall that there is a decomposition $\Sigma_3 \cong C_3 \rtimes C_2 $ where $C_2$ acts nontrivially on $C_3$. The HFPSS for the $C_2$-action on the parametrized $C_3$-homotopy fixed points $\TMF(2)_3^{h_{C_2}C_3}$ (see \cite{GM95} \cite[1.6. Definition]{QS21}) degenerates to give an isomorphism
\[ 
	\pi_* \TMF_3 \cong \pi_*(\TMF(2)_3^{h_{C_2}C_3})^{hC_2}. 
\]

We remark that $K(2)$-locally, there is a maximal finite subgroup $G_{12} = C_3 \rtimes C_4$ of the second Morava stabilizer group where $C_4$ acts nontrivially on $C_3$.  We have 
\[ 
	\Sigma_3 = G_{12}/C_2' = C_3 \rtimes (C_4/C_2') 
\]
where $C_2'$ is the subgroup of $C_4$ of order $2$.
There is an equivalence 
\[ 
	\TMF(2)_{K(2)} \simeq E_2^{hC_2' \times \mr{Gal}(\FF_9/\FF_3)} 
\]
where $E_2$ is the Morava $E$-theory spectrum associated to the formal group of the unique supersingular elliptic curve over $\FF_9$.  Thus we have (see \cite[Proposition 1.6.12]{BehrensTMF})
\[ 
	\TMF_{K(2)} \simeq E^{h(G_{12} \rtimes \rm{Gal})}_2 \simeq (E^{hC'_2 \rtimes \mr{Gal}}_2)^{h\Sigma_3} .
\]

With respect to the decomposition
$\Sigma_3 = C_3 \rtimes C_2 $, let $\gamma$ be a choice of generator of the subgroup $C_3$ and $\mu$ be the generator of the subgroup $C_2$.  
We have \cite[Sec.~7]{Stojanoska} 
\[ 
	\pi_* \tmf(2)_{3} \cong \ZZ_{3}[e_1, e_2]. 
\]
The $\Sigma_3$-action is encoded by
\begin{align*}
\gamma(e_1)  = e_2,~~ &\gamma(e_2)  = -e_1-e_2, \\
\mu(e_1)  = -e_2,~~  &\mu(e_2)  = -e_1.
\end{align*}
Define the element $\Delta^{1/2} \in \pi_{12}\tmf(2)_3$ to be 
\[ 
	\Delta^{1/2} \coloneqq 4e_1e_2(e_1+e_2) 
\]
so we have $(\Delta^{1/2})^2 = \Delta \in \pi_{24}\tmf(2)_3$ \cite[Eq.~(16)]{Stojanoska}.  
Note that
\[ 
 	\Delta^{1/2} = -4e_1 \cdot \gamma(e_1) \cdot \gamma^2(e_1). 
\]
We will now compute the ($C_3$-equivariant $3$-primary) MANSS for $\tmf(2)_3^{h}$ where 
\[
	\tmf(2)_3^{h}\coloneqq F({EC_{3}}_{+},\tmf(2)_3)
\]
is the Borel completion of  $\tmf(2)_3$. 
\begin{remark}
We warn the reader that the MANSS computation for $\tmf(2)_3$ in this subsection is  not particularly novel or deep.  We shall see that the MANSS in this case is isomorphic to the $RO(C_3)$-graded HFPSS.  This is expected given the fact that $E_2$ is $\mu_3$-oriented in the sense of \cite{HahnSengerWilson} and the MANSS should agree with the homotopy completion of the slice spectral sequence in this case. 

Indeed, the behaviour of the $RO(C_3)$-graded HFPSS spectral sequence of $E_2$ is well-known to experts, and is essentially described in \cite{HillHopkinsRavenelodd}, where it is related to a conjectural description of the odd primary slice spectral sequence.  This conjectural description is verified in \cite{Wilsonthesis}.  A detailed account of the behavior of the $RO(C_p)$-graded HFPSS for the Morava $E$-theory spectra $E_{n(p-1)}$ is the subject of a forthcoming paper of Andrew Senger. We are grateful to Senger for both generously sharing his results, and allowing us to use the $n=1$ case of his computations to illustrate this odd primary example of the MANSS.
\end{remark}

Write $RO(C_3) = \Z\{1, \lambda\}$, where $\lambda$ is the standard $2$-dimensional representation of $C_3$ given by rotation by $120$ degrees.  Note that the real regular representation $\rho = \R[C_3]$ of $C_3$ is given by $\rho = 1 +\lambda.$
Following
\cite{HillHopkinsRavenelodd,Wilsonthesis,HahnSengerWilson}, we shall find that the computation will be cleaner via the introduction of a spoke-enhanced representation grading.  Let
\[ 
	S^\spoke \subset S^\lambda 
\]
denote the unreduced suspension of the 3rd roots of unity in $\lambda$ (so that $S^\spoke$ is the $1$-skeleton of the usual $C_3$-CW complex structure of $S^\lambda$).  Since $S^\lambda$ is the $C_3$-analog of the $C_2$-representation sphere $S^{2\sigma}$, in some sense $S^\spoke$ plays the role of ``$S^{\lambda/2}$'', as it is the $C_3$ analog of $S^\sigma$.  For example, there is a cofiber sequence
\begin{equation}\label{eq:spokecofiber}
 (C_3)_+ \to S^0 \to S^\spoke.
 \end{equation}
However, we avoid taking this analogy too literally, as 
\[ 
	S^\spoke \otimes S^\spoke \simeq S^\lambda \vee (C_3)_+\otimes S^2 \not\simeq S^\lambda. 
\]
We define $S^{-\spoke} \coloneqq D(S^\spoke)$, but note that $S^\spoke$ is not invertible.  We do have
\[ 
	S^{\lambda - \spoke} \simeq S^\spoke 
\]
which implies that
\[ 
	S^{-\spoke} \simeq S^{\spoke}\otimes S^{-\lambda}.
\]
We therefore define
\[
	RO^\spoke(C_3) = \frac{RO(C_3) \oplus \ZZ\{\spoke\}}{(2\spoke = \lambda)}. 
\]
and we then have
\[ 
	-(\spoke) = \spoke - \lambda \in RO^\spoke(C_3). 
\]
Given $i+j\lambda + \epsilon\spoke$ (with $\epsilon \in \{0,1\}$), we define for $X \in \Sp_{C_3}$, 
\[
	\pi^{C_3}_{i+j\lambda + \epsilon\spoke}(X) \coloneqq [S^{i+j\lambda} \otimes (S^\spoke)^{\epsilon}, X]^{C_3}.
\]
In this subsection, we will use $\pi^{C_3}_\sstar$  to denote $RO^\spoke(C_3)$-graded homotopy groups.  If $R \in \Sp_{C_3}$ is a ring spectrum, then $\pi^{C_3}_\sstar R$ is a module over $\pi^{C_3}_{\star}R$.
Let
\[ 
	a_\spoke \colon \thinspace S^{-\spoke} \to S^0 
\]
denote the Spanier-Whitehead dual of the second map of (\ref{eq:spokecofiber}).  The $\lambda$-suspension of this map is the canonical inclusion
\[ 
	a_\spoke  \colon \thinspace S^\spoke \hookrightarrow S^\lambda. 
\]
If we abusively also use $a_\spoke$ to denote the Spanier-Whitehead dual of this map
\[
	a_\spoke  \colon \thinspace S^{-\lambda} \to S^{-\spoke} 
\]
then it makes sense to form $a^2_\spoke$, and we have
\[ 
	a_\spoke^2 = a_\lambda. 
\]
For any $X \in \Sp_{C_3}$, the $RO^\spoke(C_3)$-graded homotopy groups $\pi^{C_3}_\sstar$ are a module over $\ZZ[a_\spoke]$.

Our goal is to compute the $RO^\spoke(C_3)$-graded MANSS with signature
\[ 
	\E{\MANSS}{}^{s,\sstar}(\tmf(2)_3^{h}) \implies \pi^{C_3}_{\sstar}\tmf(2)^h_3. 
\]
Consider the second double complex spectral sequence
\[ 
	\E{}{II}^{a,b,\sstar}(\tmf(2)_3^{h}) = H^a(C_3; \Ext^{b,a+b}((MU_{3})_*(\tmf(2)_3 \otimes S^{-\sstar}))) \implies \E{\MANSS}{}^{a+b,\sstar}(\tmf(2)_3^{h}). 
\]
Since $\pi_*\tmf(2)_3$ is concentrated in even degrees, it is complex orientable, and therefore its ANSS collapses with $\Ext^{b,*}((MU_{3})_*\tmf(2)_3) = 0$ for $b > 0$.  In particular, the second double complex spectral sequence collapses to give an isomorphism
\begin{equation}\label{eq:MANSSE2tmf2}
	\E{\MANSS}{}_2^{s,\sstar}(\tmf(2)_3^{h}) \cong H^s(C_3; \pi_{s}(\tmf(2)_3 \otimes S^{-\sstar})).
\end{equation}
It follows from Proposition~\ref{prop: MANSS to HFPSS} that the MANSS is isomorphic to the HFPSS. We first identify the $E_2$-term of this spectral sequence.

\begin{proposition}
There are preferred generators $c_4$ and $c_6$ and an additive isomorphism
\begin{multline}\label{eq:MANSSE2tmf2comp}
\E{\MANSS}{}^{*,\sstar}(\tmf(2)_3^{h}) \cong \\ \frac{\pi^{C_3}_\sstar H\ZZ_3^h [c_4, c_6, \Delta^{1/2}]/(a_\spoke c_4, \: a_\spoke c_6, \: c_4^3 - c_6^2 -12^3(\Delta^{1/2})^2) \{\bar{v}_1, \bar{v}'_1, u_\spoke c_4, u'_\spoke c_4, u_\spoke c_6, u'_\spoke c_6  \}}{(c_4\cdot \bar{v}_1, \: c_6\cdot \bar{v}_1, \: a_\spoke\cdot \bar{v}_1', \: c_4 \cdot \bar{v}_1', \: {c_6} \cdot \bar{v}_1', a_\spoke \cdot u_\spoke c_4, \: a_\spoke \cdot u'_\spoke c_4, \:  a_\spoke \cdot u_\spoke c_6, \: a_\spoke \cdot u'_\spoke c_6, c_6 \cdot u_\spoke c_6)}.
\end{multline}
Furthermore the inclusion 
\[ 
	\pi^{C_3}_\sstar H\ZZ_3^h [c_4, c_6, \Delta^{1/2}]/(a_\spoke c_4, \: a_\spoke c_6, \: c_4^3 - c_6^2 -12^3(\Delta^{1/2})^2) \hookrightarrow \E{\MANSS}{}^{*,\sstar}(\tmf(2)_3^{h}) 
\]
 is an inclusion of rings, and the isomorphism \eqref{eq:MANSSE2tmf2comp}
is an isomorphism of modules over this subring.
\end{proposition}
\begin{proof}
Our first task will be to compute $\pi_*\tmf(2)_3$ as a $C_3$-module.  Let $\ZZ^\spoke$ denote the kernel of the augmentation
$ \ZZ[C_3] \to \ZZ.$
Note that we have an isomorphism of $C_3$-modules:
\[
	H_1(S^\spoke;\ZZ) \cong \ZZ^\spoke.
\]
As a $\ZZ[C_3]$-module we have
\begin{align*}
\pi_0\tmf(2)_3 & = \ZZ_3\{1\} \\
\pi_4\tmf(2)_3 & = \ZZ^\spoke_3\{e_1\} \\
\pi_8\tmf(2)_3 & = \ZZ_3[C_3]\{e_1^2\} \\
\pi_{12}\tmf(2)_3 & = \ZZ_3\{\Delta^{1/2}\} \oplus \ZZ_3[C_3]\{e^2_1e_2\} \\
\pi_{16}\tmf(2)_3 & = \ZZ^{\spoke}_3\{\Delta^{1/2}e_1\} \oplus \ZZ_3[C_3]\{e_1^4\} \\
\pi_{20}\tmf(2)_3 & = \ZZ_3[C_3]\{\Delta^{1/2}e^2_1, e_1^4e_2\} \\
& \vdots
\end{align*}
In general, we have an isomorphism of $C_3$-modules
\begin{equation}\label{eq:tmf2C3}
 \pi_*\tmf(2)_3 \cong \ZZ_3[\Delta^{1/2}]\otimes \left(\ZZ\{1\} \oplus \ZZ^\spoke\{e_1\} \oplus \ZZ[C_3]\{ e_1^{2i} e_2^{\epsilon} \: : \: i > 0, \epsilon \in \{0,1\}\} \right).
 \end{equation}

Note that we have (using the analog of \cite[Theorem~V.2.4]{Hinfty})
\begin{align*}
\pi^{C_3}_{i+j\lambda+\epsilon \spoke}H\ZZ_3^h & \cong H^{-i}((EC_3)_+ \otimes_{C_3} S^{j\lambda+\epsilon \spoke}; \ZZ_3) \\
  & \cong H^{-i-2j-\epsilon}(C_3; \ZZ_3^{\epsilon\spoke}).
  \end{align*}
It follows from elementary group cohomology computations that we have
\begin{equation}\label{eq:pistarHZC3}
 \pi_\sstar^{C_3} H\ZZ_3^h \cong \ZZ_3[u^{\pm}_\lambda, a_\spoke]/(3a_\spoke)
 \end{equation}
with
\begin{align*}
\abs{u_\lambda} = 2-\lambda, \text{ and }\abs{a_\spoke} = -\spoke.
\end{align*}

Observe that we have
\begin{align*} \E{\MANSS}{}^{0,\lambda+\spoke+1}(\tmf(2)_3)& = H^0(C_3; \pi_0(\tmf(2)_3 \otimes S^{-\lambda-\spoke-1}))\\
& \cong H^0(C_3; \ZZ_3^\spoke \otimes \ZZ^\spoke) \\
& \cong H^0(C_3; \ZZ_3) \oplus H^0(C_3; \ZZ_3[C_3]) \\
& \cong \ZZ_3 \oplus \ZZ_3.
\end{align*}
We let $\bar{v}_1$ denote a generator of the first summand, and $\bar{v}'_1$
denote a generator of the second summand. 
We let $c_4$ be a generator of 
\begin{align*}
\E{\MANSS}{}^{0,8}(\tmf(2)_3^{h}) & = H^0(C_3; \pi_{8}\tmf_0(2)_3) \\
& = H^0(C_3; \ZZ_3[C_3]\{ e_1^2 \}) \\
& = \ZZ_3, 
\end{align*}
and let $u_\spoke c_4$ and $u'_\spoke c_4$ denote the generators of 
\begin{align*}
\E{\MANSS}{}^{0,7+\spoke}(\tmf(2)_3^{h}) & = H^0(C_3; \pi_{7}\tmf_0(2)_3 \otimes S^{-\spoke}) \\
& = H^0(C_3; \ZZ_3[C_3]\{ e_1^2 \}\otimes \ZZ^\spoke) \\
& = \ZZ_3 \oplus \ZZ_3. 
\end{align*}
Similarly we take $c_6$ to be a generator of the second summand of 
\begin{align*}
\E{\MANSS}{}^{0,12}(\tmf(2)_3^{h}) & = H^0(C_3; \pi_{12}\tmf_0(2)_3) \\
& = H^0(C_3; \ZZ_3\{\Delta^{1/2}\} \oplus \ZZ_3[C_3]\{ e_1^2e_2 \}) \\
& = \ZZ_3 \oplus \ZZ_3 
\end{align*}
where the first summand is generated by $\Delta^{1/2}$.
We let $u_\spoke c_6$ and $u'_\spoke c_6$ be the generators of the first and second summands of
\begin{align*}
\E{\MANSS}{}^{0,11+\spoke}(\tmf(2)_3^{h}) & = H^0(C_3; \pi_{11}\tmf_0(2)_3 \otimes S^{-\spoke}) \\
& = H^0(C_3; \ZZ_3\{\Delta^{1/2}\}\otimes \ZZ^\spoke \oplus \ZZ_3[C_3]\{ e_1^2e_2 \}\otimes \ZZ^\spoke) \\
& = \ZZ_3 \oplus \ZZ_3. 
\end{align*}
If we take $c_4$ and $c_6$ to be the generators coming from the image of 
\[ 
	\pi_*\tmf_3 \to \pi_*\tmf(2)_3 
\]
then the result follows from (\ref{eq:MANSSE2tmf2}), (\ref{eq:tmf2C3}), and (\ref{eq:pistarHZC3}). 
\end{proof}

The $\mu_3$-orientation of $\tmf(2)$ in \cite[Theorem 1.9.]{HahnSengerWilson} gives rise to a $C_3$-equivariant map
\[ 
	\bar{v}_1  \colon \thinspace S^{\lambda+\spoke+1} \to \tmf(2)_3 
\]
which is an isomorphism on $\pi^e_4$ (after $3$-completion).
Since $\tmf(2)_3$ is a $C_3$-commutative ring spectrum, the norm-restriction adjunction \cite[Prop.~2.27]{HHR16} gives rise to a map
\[ 
	N_e^{C_3} \mr{Res}_e^{C_3} \tmf(2)_3 \to \tmf(2)_3. 
\]
Therefore, we may norm the element $e_1 \in \pi^e_{4}\tmf(2)_3$ to get an equivariant enhancement of $\Delta^{1/2}$:
\[ 
	\bar{\Delta}^{1/2} \coloneqq -4N^{C_3}(e_1) \in \pi^{C_3}_{4\rho}\tmf(2)_3. 
\]

\begin{proposition}
In the MANSS for $\mathrm{tmf}(2)_3$ there are differentials
\begin{align*}
d_5(u_\lambda^8 (\bar{\Delta}^{1/2})^2) & \doteq a^5_\spoke u_\lambda^7 (\bar{\Delta}^{1/2})^2, \\
d_5(u_\lambda) & \doteq a_\spoke^5 \bar{v}_1, \\
d_9(a_\spoke u^{17}_\lambda \bar{v}_1 (\bar{\Delta}^{1/2})^4) & \doteq a_\spoke^{10} u^{15}_\lambda (\bar{\Delta}^{1/2})^5  \\
d_9(u^2_\lambda \bar{v}_1) & \doteq a_\spoke^9 \bar{\Delta}^{1/2}.
\end{align*}
These generate all of the differentials in the MANSS under the Leibniz rule and in particular we find that 
\[ 
	u_\lambda^i u^\epsilon_\spoke c_4, \: u'_\spoke c_4, \: u_\lambda^i u_\spoke^\epsilon c_6, \: u'_\spoke c_6 , \: u_\lambda \bar{v}_1, \: \bar{v}'_1, u_\lambda^3 
\]
are permanent cycles. It follows that $(\Delta^{1/2})^{3} = u^{12}_{\lambda}(\bar{\Delta}^{1/2})^3$ is a permanent cycle, and the MANSS is therefore $(\Delta^{1/2})^{3}$-periodic.
\end{proposition}
\begin{proof}
The HFPSS/MANSS for $\tmf(2)_3$ maps to the HFPSS for $\TMF(2)_3$ converging to the homotopy groups of 
\[ 
	\TMF(2)_3^{h{C_3}} \simeq \TMF_3 \vee \Sigma^{36} \TMF_3.
\]
The structure of the HFPSS for $\TMF(2)_3$ in integer degrees is well known. It is a shifted direct sum of two copies of the HFPSS which converges to $\pi_*\TMF_3$ \cite{Bauer}. 
The $E_2$-term is given by  
\[ 
	\E{HFPSS}{}_2^{s,t}(\TMF(2)_3) \cong H^s(C_3; \pi_t \TMF(2)_3 ) \cong \frac{\ZZ_3[c_4, c_6, \Delta^{\pm 1/2}, \beta]\otimes E[\alpha]}{(c_4^3-c_6^2-12^3(\Delta^{1/2})^2, 3\alpha, 3\beta, \alpha c_4, \alpha c_6,  \beta c_4,  \beta c_6 )} 
\]
and under the map from the MANSS for $\tmf(2)_3$ to the HFPSS for $\TMF_0(2)$
we have
\begin{align*}
a_\spoke u_\lambda \bar{v}_1  \mapsto \alpha, \text{ and }
a_\spoke^2 u_\lambda^3 \bar{\Delta}^{1/2}  \mapsto \beta.
\end{align*}
In the HFPSS for $\TMF(2)_3$ there are differentials
\begin{align*}
d_5(\Delta)  \doteq \alpha \beta^2,  \text{ and }
d_9(\Delta^2 \alpha)  \doteq \beta^5.
\end{align*}
Here, we use $\doteq$ to indicate both sides of the equation are equal up multiplication by an element of $\ZZ^\times_3$.
This implies that in the MANSS for $\tmf(2)_3$, we have differentials
\begin{align*}
d_5(u_\lambda^8 (\bar{\Delta}^{1/2})^2) & \doteq a^5_\spoke u_\lambda^7 (\bar{\Delta}^{1/2})^2, \\
d_9(a_\spoke u^{17}_\lambda \bar{v}_1 (\bar{\Delta}^{1/2})^4) & \doteq a_\spoke^{10} u^{15}_\lambda (\bar{\Delta}^{1/2})^5.
\end{align*}
The element $a_\spoke$ is a permanent cycle because it is in the Hurewicz image, and the elements $\bar{v}_1$ and $\bar{\Delta}^{1/2}$ are permanent cycles because we constructed elements of $\pi^{C_3}_\sstar \tmf(2)_3^h$ which these converge to.  It follows that we have differentials
\begin{align*} 
d_5(u_\lambda)  \doteq a_\spoke^5 \bar{v}_1, \text{ and }
d_9(u^2_\lambda \bar{v}_1)  \doteq a_\spoke^9 \bar{\Delta}^{1/2}.
\end{align*}
By dimension considerations, we observe that these differentials generate all of the differentials in the MANSS under the Leibniz rule and the result follows.
\end{proof}

\bibliographystyle{alpha}
\bibliography{bib}

\end{document}